\documentclass[11pt,a4paper]{article}
\usepackage{ifthen,latexsym,amssymb,amsmath,bbm,fixmath}
\usepackage{bm,booktabs,float}
\usepackage{graphicx,tikz}
\usepackage{authblk}
\usepackage{enumitem}

\usepackage{graphicx}
\usepackage{multicol,multirow}
\usepackage{amsmath,amssymb,amsfonts}
\usepackage{mathrsfs}
\usepackage{amsthm}
\usepackage{rotating}
\usepackage{appendix}
\usepackage[numbers]{natbib}
\usepackage{ifpdf}
\usepackage[T1]{fontenc}
\usepackage{newtxtext}
\usepackage{newtxmath}
\usepackage{textcomp}
\usepackage[colorlinks]{hyperref}

\newtheorem{theorem}{Theorem}[section]
\newtheorem{lemma}[theorem]{Lemma}
\newtheorem{definition}[theorem]{Definition}
\newtheorem{corollary}[theorem]{Corollary}
\newtheorem{conjecture}[theorem]{Conjecture}
\newtheorem{proposition}[theorem]{Proposition}
\theoremstyle{definition}

\numberwithin{equation}{section}

\newtheorem{claim}{Claim}[theorem]
\newtheorem{claim2}{Claim}[section]

\newcommand{\bcpf}{\begin{proof}[Proof of Claim]}
\newcommand{\ecpf}{\end{proof} \medskip}

\usepackage{bm}
\usepackage{enumitem}

\begin{document}


\def\eps{{\varepsilon}}
\newcommand{\cP}{\mathcal{P}}
\newcommand{\cT}{\mathcal{T}}
\newcommand{\cL}{\mathcal{L}}
\newcommand{\ex}{\mathbb{E}}
\newcommand{\eul}{e}
\newcommand{\pr}{\mathbb{P}}
\newcommand{\supp}{\mathop{\mathrm{supp}}}
\newcommand{\Capa}{\mathop{\mathrm{Cap}}}
\newcommand{\opt}{\mathop{\textsc{opt}}}
\newcommand{\feas}{\mathop{\textsc{feas}}}
\newcommand{\pat}{\mathop{\textsc{pat}}}
\newcommand{\wt}{\mathop{\textsc{wt}}}
\newcommand{\coars}{\mathop{\mathrm{Coars}}}
\newcommand{\RL}{\mathop{\mathrm{RL}}}
\newcommand{\ext}{\mathop{\mathrm{ext}}}
\newcommand{\ba}{\bm{\alpha}}
\newcommand{\bb}{\bm{\beta}}

\title{Stability for the Erd\H{o}s-Rothschild problem}
\author[1]{Oleg Pikhurko}
\affil[1]{Mathematics Institute and DIMAP, University of Warwick, Coventry CV4 7AL, UK}
\author[2]{Katherine Staden}
\affil[2]{School of Mathematics and Statistics, The Open University, Walton Hall, Milton Keynes MK7 6AA, UK}
\maketitle

\begin{abstract}
Given a sequence $\bm{k} := (k_1,\ldots,k_s)$ of natural numbers and
a graph $G$, let $F(G;\bm{k})$ denote the number of colourings of the edges of $G$ with colours $1,\dots,s$ such that, for every  $c \in \{1,\dots,s\}$,
the edges of colour $c$ contain no clique of order $k_c$. 
Write $F(n;\bm{k})$ to denote the maximum of $F(G;\bm{k})$ over all graphs $G$ on $n$ vertices.
This problem was first considered by Erd\H{o}s and Rothschild in 1974,
but it has been solved only for a very small number of non-trivial cases.
In previous work with Yilma, we 
constructed a finite optimisation problem whose maximum is equal to the limit of $\log_2 F(n;\bm{k})/{n\choose 2}$
as $n$ tends to infinity
and proved a stability theorem for complete multipartite graphs $G$.

In this paper we provide a sufficient condition on $\bm{k}$ which guarantees a general stability theorem for \emph{any} graph $G$, describing the asymptotic structure of $G$ on $n$ vertices with $F(G;\bm{k}) = F(n;\bm{k}) \cdot 2^{o(n^2)}$ in terms of solutions to the optimisation problem.
We apply our theorem to systematically recover existing stability results as well as all cases with $s=2$. 
The proof uses a version of symmetrisation on edge-coloured weighted multigraphs.
\end{abstract}

\vspace*{14pt}
\section{Introduction}

Let a non-increasing sequence $\bm{k} = (k_1,\ldots,k_s) \in \mathbb{N}^s$ of natural numbers be given. By an \emph{$s$-edge-colouring} (or \emph{colouring} for brevity)
of a graph $G=(V,E)$ we mean a function $\chi:E\to [s]$, where we denote $[s]:=\{1,\dots,s\}$. Note that colourings do not have to be proper, that is, adjacent edges can have the same colour. A colouring $\chi$ of $G$ is called $\bm{k}$-\emph{valid} if, for every $c\in [s]$, the colour-$c$ subgraph $\chi^{-1}(c)$ contains no copy of $K_{k_c}$, the complete graph of order $k_c$. 
Write $F(G; \bm{k})$ for the number of $\bm{k}$-valid colourings of $G$.

In a previous paper with Yilma~\cite{psy}, we investigated  the \emph{Erd\H{o}s-Rothschild problem} of determining $F(n;\bm{k})$, the maximum of $F(G;\bm{k})$ over all graphs $G$ on $n$ vertices, and the \emph{$\bm{k}$-extremal} graphs, i.e.~order-$n$ graphs which attain this maximum. We assume throughout the paper, as we did there, that $s\ge 2$ and that $k_c \geq 3$ for all $c \in [s]$ (since $k_c=2$ just forbids
colour $c$ and the problem reduces to one with $s-1$ colours).

The case when $k_1 = \ldots = k_s =: k$, which we denote by $\bm{k}=(k_1,\ldots,k_s)=(k;s)$, was first considered by Erd\H{o}s and Rothschild in 1974 (see~\cite{ER,ER2}). 
A trivial lower bound on $F(n;(k;s))$ is obtained by taking the largest $K_k$-free graph on $n$ vertices, namely the Tur\'an graph $T_{k-1}(n)$ which is the complete $(k-1)$-partite graph with parts as equal as possible.
Any $s$-edge-colouring of this graph is $\bm{k}$-valid, so we have
\begin{equation}\label{trivial}
F(n;(k;s)) \geq s^{t_{k-1}(n)},
\end{equation}
where $t_{k-1}(n)$ is the number of edges in $T_{k-1}(n)$.
In particular, Erd\H{o}s and Rothschild conjectured that, when $\bm{k}=(3,3)$ and $n$ is sufficiently large, the trivial lower bound~(\ref{trivial}) is in fact tight and, furthermore,  $T_{2}(n)$ is the unique
$\bm{k}$-extremal graph on $n$ vertices.
The  conjecture was verified for all $n \geq 6$ by Yuster~\cite{yuster} (who also computed $F(n;(3,3))$ for smaller $n$).
Proving Yuster's extension of the conjecture,  Alon, Balogh, Keevash and Sudakov~\cite{abks} showed that an analogous result holds for two and three colours: for large $n$, the Tur\'an graph $T_{k-1}(n)$ is the unique $\bm{k}$-extremal graph for $\bm{k}=(k,k)$ and $\bm{k}=(k,k,k)$.
The proof of this result uses Szemer\'edi's Regularity Lemma, so the graphs to which it applies are very large indeed.
However, the assertions are not true for all numbers $n$ of vertices.
As was remarked in~\cite{abks}, the assertions do not hold when $k \leq n < s^{(k-2)/2}$, as in this case a random colouring of the edges of $K_n$ with $s$ colours contains no monochromatic $K_k$ with probability more than $\frac{1}{2}$.
Thus, for this range of $n$, we have $F(n;(k;s)) > s^{\binom{n}{2}}/2 \geq s^{t_{k-1}(n)}$.
H\`an and Jim\'enez~\cite{HJcontainer} used graph containers to obtain an essentially optimal lower bound for the order $n$ of graphs for which the trivial lower bound~(\ref{trivial}) for $s=2,3$ is tight. 

In this paper we are only interested in large $n$.
It was proved in~\cite[Proposition~5.1]{abks} that the limit 
 \begin{equation}\label{eq:lim}
 F(\bm{k}):=\lim_{n\to\infty} \frac{\log_2 F(n;\bm{k})}{n^2/2}
 \end{equation}
 exists (and is positive) when $\bm{k}=(k;s)$. 
It can be easily seen that the proof from~\cite{abks} extends to
an arbitrary fixed sequence $\bm{k}$.
The authors of~\cite{abks} noted that when more than three colours are used, the behaviour of $F(n;(k;s))$ changes, making its determination both harder and more interesting. Namely, it was shown in~\cite[page~287]{abks} 
that if $s \geq 4$ (and $k \geq 3$) then $F(n;(k;s))$ is exponentially in $n^2$ larger than $s^{t_{k-1}(n)}$. 
In particular, any extremal graph has to contain many copies of $K_k$.
The authors of~\cite{abks} determined $F(\bm{k})$ for $\bm{k}=(3,3,3,3)$ and $\bm{k}=(4,4,4,4)$,
where $T_4(n)$ and $T_9(n)$ respectively achieve the right exponent.
Pikhurko and Yilma~\cite{PY} were able to obtain an exact result for these cases: that these Tur\'an graphs are the unique respective extremal graphs.
Recently, Botler, Corsten, Dankovics, Frankl, H\`an, Jim\'enez and Skokan~\cite{skokan} announced the determination of $F(\bm{k})$ for $\bm{k}=(3;5),(3;6)$. For $s=6$ they proved that $T_8(n)$ is the unique extremal graph, and also proved a stability result.
For $s=5$ they uncovered new behaviour: for large $n$ there is an infinite family $\{S_{\alpha,\beta}(n): 0 \leq \alpha+\beta \leq \frac{1}{4}\} \cup \{T_{\alpha,\beta}(n): 0 \leq \alpha, \beta \leq \frac{1}{4}\}$ of asymptotically optimal graphs  with either $4$, $6$ or $8$ parts,
where $S_{\alpha,\beta}(n)$ denotes the complete partite graph with parts of size
$\frac{n}{4},\frac{n}{4},\alpha n,\alpha n,\beta n,\beta n,(\frac{1}{4}-\alpha-\beta) n,(\frac{1}{4}-\alpha-\beta) n$ and
$T_{\alpha,\beta}(n)$ denotes the complete partite graph with parts of size
$\alpha n, \alpha n, (\frac{1}{4}-\alpha)n, (\frac{1}{4}-\alpha)n, \beta n, \beta n, (\frac{1}{4}-\beta)n, (\frac{1}{4}-\beta)n$.
These are the only known results, asymptotic or exact.

\begin{table}
\centering
\caption{Known results}
\vspace{0.5cm}
\label{knownsolutions}
\begin{tabular}{rl|l|l|l}
$\bm{k}=$ & $(k;s)$        &                     $F(\bm{k})$        & extremal graph                           & citation \\ 
\hline
any $k$ & $s=2$   & $1-\frac{1}{k-1}$                     & $T_{k-1}(n)$                 &         \cite{abks}                 \\
 & $s=3$  & $(1-\frac{1}{k-1})\log_2 3$                            & $T_{k-1}(n)$              &     \cite{abks}      \\
$k=3$ & $s=4$           &  $\frac{1}{4}+\frac{1}{2}\log_2 3$                        & $T_{4}(n)$      &         \cite{abks,PY}  \\
 & $s=5$           &      $\frac{1}{2}+\frac{1}{2}\log_2 3$                 & $S_{\alpha,\beta}(n), T_{\alpha,\beta}(n)$  ~~~  (*)           &     \cite{skokan} \\      
 & $s=6$           &          $\frac{3}{4} + \frac{1}{2}\log_2 3$          & $T_{8}(n)$      &    \cite{skokan}           \\      
$k=4$ & $s=4$           &      $\frac{8}{9}\log_2 3$                       & $T_{9}(n)$                        &    \cite{abks,PY}                  
\end{tabular}

(*) These graphs are known to be asymptotically extremal only: they achieve the right exponent in $F(n;\bm{k})$.
\end{table}

Many other versions of the Erd\H{o}s-Rothschild problem have been studied, where the goal is to maximise the number of colourings of some discrete object when one forbids certain coloured substructures.
Erd\H{o}s and Rothschild themselves considered the generalisation where
one forbids a monochromatic graph $H$. In~\cite{abks} the authors showed that the trivial lower bound~(\ref{trivial}) is tight for large $n$
when $H$ is \emph{colour-critical}, that is, the removal of any edge 
from $H$ reduces its chromatic number.
(Note that every clique is colour-critical.)
In a further generalisation, Balogh~\cite{balogh} considered edge-colourings in which a specific colouring of a fixed graph $H$ is forbidden.
Other authors have addressed this question in the cases of forbidden monochromatic matchings, stars, paths, trees and some other graphs in~\cite{hkl2,hkl3}, matchings with a prescribed colour pattern in~\cite{hl}, rainbow stars in~\cite{hlos} and multicoloured cliques in~\cite{hln}.
A so-called `$q$-analogue' was addressed in~\cite{qana}, which considers a related problem in the context of vector spaces over a finite field $GF(q)$.
Alon and Yuster~\cite{alonyuster} studied a directed version of the problem, to determine the maximum number of $T$-free orientations of an $n$-vertex graph, where $T$ is a given $k$-vertex tournament. 
The problem of counting monochromatic $H$-free edge-colourings in hypergraphs was studied in~\cite{hkl,lprs,lps}.
Additive versions have also been studied, where an underlying group~\cite{HJ} or set~\cite{LSS} with addition is coloured, and monochromatic triples $(x,y,z)$ with $x+y=z$ are forbidden. 

\subsection{An optimisation problem}
This paper concerns the relation between the structure of almost extremal graphs for $F(n,\bm{k})$ and optimal solutions of a certain optimisation problem, Problem $Q_t$, which we now define.

\noindent\textbf{Problem~$Q_t$}:
 \it
 Given a sequence $\bm{k} := (k_1,\ldots,k_s) \in \mathbb{N}^s$ of natural numbers and $t\in\{0,1,2\}$, determine
\begin{equation}\label{Qdef}
Q_t(\bm{k}) := \max_{(r,\phi,\ba) \in \feas_t(\bm{k})} q(\phi,\ba),
\end{equation}
 the maximum value of
 \begin{equation}\label{q}
 q(\phi,\ba) := 2\sum_{\stackrel{1 \leq i < j \leq r}{\phi(ij) \neq \emptyset}} \alpha_i \alpha_j \log_2 |\phi(ij)|
 \end{equation}
 over the set $\feas_t(\bm{k})$ of \emph{feasible solutions}, that is, triples
$(r,\phi,\ba)$ such that
 \begin{itemize}
  \item $r \in \mathbb{N}$ and $\phi \in \Phi_t(r;\bm{k})$, where $\Phi_t(r;\bm{k})$ is the set of all
  functions $\phi : \binom{[r]}{2} \rightarrow 2^{[s]}$
  such that
  \[
  \phi^{-1}(c) := \left\{ ij \in \binom{[r]}{2} : c \in \phi(ij) \right\}
  \]
is $K_{k_c}$-free for every colour $c \in [s]$ and $|\phi(ij)|\ge t$
for all $ij\in \binom{[r]}{2}$;
  \item $\ba = (\alpha_1,\ldots,\alpha_r) \in \Delta^r$, where $\Delta^r$ is the set of all $\ba \in \mathbb{R}^r$ with $\alpha_i \geq 0$ for all $i \in [r]$, and $\alpha_1 + \ldots + \alpha_r = 1$.
  \end{itemize}
 \rm
 Note that for $t=1,2$, a triple $(r,\phi,\ba) \in \feas_t(\bm{k})$ necessarily has $r<R(\bm{k})$ 
where $R(\bm{k})$ is the \emph{Ramsey number} of $\bm k$ (i.e.\ the minimum $R$ such that $K_R$ admits no $\bm k$-valid $s$-edge-colouring).
Thus the maximum in~\eqref{Qdef} is attained for $t=1,2$ since $q(r,\phi,\cdot)$ is continuous for each of the finitely many pairs $(r,\phi)$, and $\feas_t(\bm{k})$ is a (non-empty) compact space. 
It is also attained for $t=0$ by~(\ref{Qt}) below.
We call $\phi \in \Phi_0(r;\bm{k})$ a \emph{colour pattern} and $\ba \in \Delta^r$ a \emph{vertex weighting}.
A triple $(r,\phi,\ba)$ is called \emph{$Q_t$-optimal} if it attains the maximum, that is, 
$(r,\phi,\ba)\in \feas_t(\bm{k})$ and $q(r,\phi,\ba)=Q_t(\bm{k})$.
One can easily show~\cite[Lemma~6]{psy} that
\begin{equation}\label{Qt}
Q(\bm{k}) := Q_2(\bm{k}) = Q_1(\bm{k})=Q_0(\bm{k}).
\end{equation}
Note that a $Q_0$-optimal triple can have $r$ arbitrarily large, by e.g.~adding vertices of weight 0 or splitting an existing vertex into two clones.
Let $\opt_t(\bm{k})$ be the set of $Q_t$-optimal triples $(r,\phi,\ba)$.
Let $\feas^*(\bm{k})$ be the set of $(r,\phi,\ba) \in \feas_2(\bm{k})$ with $\alpha_i>0$ for every $i \in [r]$.
Let $\opt^*(\bm{k})$ be the set of \emph{basic optimal solutions}, which are $Q_2$-optimal triples $(r,\phi,\ba)$ with $\alpha_i>0$ for every $i \in [r]$.

Given vectors $\bm{a}=(a_1,\ldots,a_s)$ and $\bm{b}=(b_1,\ldots,b_t)$, write $\bm{a}\leq\bm{b}$ if $a_i \leq b_i$ for all $i \leq \max\{s,t\}$ where $a_i:=0$ for all $i>s$ and $b_i:=0$ for all $i>t$.
We write $\| \bm{a}-\bm{b}\|_1 := \sum_{i \leq \max\{s,t\}}|a_i-b_i|$ for the \emph{$\ell^1$-distance} between $\bm{a}$ and $\bm{b}$.
In this paper we always take $\log$ to the base $2$; from now on we omit any subscript. 

One should think of feasible triples $(r,\phi,\ba)$ as vertex-weighted edge-coloured multigraphs.
It is not hard to show that $F(\bm{k})\geq Q(\bm{k})$.
Indeed, given $\bm{k}$, a $Q_1$-optimal triple $(r,\phi,\ba)$ and $n \in \mathbb{N}$, take the complete $r$-partite $n$-vertex graph $K_{\ba}(n)$ whose parts $X_1,\ldots,X_r$ satisfy $|\,|X_i|-\alpha_i n\,| \leq 1$ for all $i \in [r]$.
Consider those $s$-edge-colourings of $K_{\ba}(n)$ in which, for $x \in X_i$ and $y \in X_j$, we only allow colours in $\phi(ij)$ to be used on $xy$.
Every such colouring is $\bm{k}$-valid since every $\phi^{-1}(c)$ is $K_{k_c}$-free.
Clearly $F(n;\bm{k})$ is bounded below by the number of such colourings of $K_{\ba}(n)$, which is
$$
\prod_{ij\in\binom{[r]}{2}}|\phi(ij)|^{|X_i||X_j|} = 2^{q(\phi,\ba)n^2/2+O(n)} = 2^{Q(\bm{k})n^2/2+O(n)}.
$$ 
Taking the limit as $n\to\infty$, we have $F(\bm{k})\geq Q(\bm{k})$.
With Yilma, we proved the following results relating the determination of $F(n;\bm{k})$ to Problem $Q_1$, including a matching upper bound.

\begin{theorem}\cite{psy}\label{oldresults}
The following hold for every $s \in \mathbb{N}$ and $\bm{k} \in \mathbb{N}^s$.
\begin{enumerate}[label=(\roman*),ref=(\roman*)]
\item\label{oldresultsi}
For every $n \in \mathbb{N}$, at least
one of the 
$\bm{k}$-extremal graphs of order $n$ is complete multipartite.
\item\label{oldresultsii}
$F(n;\bm{k})=2^{Q(\bm{k}){n\choose 2}+o(n^2)}$, that is, $F(\bm{k})=Q(\bm{k})$.
\item\label{oldresultsiii}
For every $\delta>0$ there are $\eta>0$ and $n_0$
such that if $G=(V,E)$ is a complete multipartite graph of order $n\ge n_0$ with (non-empty) parts $V_1,\dots,V_r$ and $F(G;\bm{k})\ge 2^{(Q(\bm{k})-\eta)n^2/2}$ then there is a $Q_1$-optimal triple $(r,\phi,\bm{\alpha'})$ such that 
$\|\ba-\ba'\|_1 \leq \delta$, where
$\ba:=(\frac{|V_1|}{n},\dots,\frac{|V_r|}{n})$.
\end{enumerate}
\end{theorem}

Thus, by~(\ref{Qt}), to determine $F(\bm{k})$, it suffices to find the optimal solutions to Problem $Q_2$, which has the smallest feasible set among the problems $Q_t$.
Unfortunately this is difficult even when $\bm{k}$ is small.
Given a pair $(r,\phi)$, one can use the method of Lagrange multipliers to find a best possible $\ba$ for this pair; though the number of pairs $(r,\phi)$ is finite there are generally too many for a computer search.
Indeed, the upper bound of $R(\bm{k})$ for $r$ grows large very quickly, though we expect the optimal $r$ to be much smaller than $R(\bm{k})$.

\subsection{New results}

The main contribution of this paper is a general stability theorem that determines the structure of any $n$-vertex graph $G$ which is almost $\bm{k}$-extremal, i.e.~with $F(G,\bm{k})=F(n;\bm{k})\cdot 2^{o(n^2)}$.  This will show that the structure of any such graph is similar to an optimal solution to Problem $Q_0$, and almost all valid colourings almost follow an optimal colour pattern.
This stability result holds for all $\bm{k}$ satisfying a rather general condition which we call the \emph{extension property}. Given $\opt^*(\bm{k})$, one can easily check whether this condition holds.
Indeed, we show that in almost all instances for which $F(\bm{k})$ is known, $\bm{k}$ satisfies a strong version of this property.

\begin{definition}[Clones and extension property]\label{extprop}
\rm
Let $s \in \mathbb{N}$ and $\bm{k} \in \mathbb{N}^s$. 
Given $r \in \mathbb{N}$ and $\phi \in \Phi_0(r;\bm{k})$, we say that $i \in [r]$ is
\begin{itemize}
\item a \emph{clone of $j \in [r]\setminus \{ i \}$} \emph{(under $\phi$)} if $\phi(ik) = \phi(jk)$ for all $k \in [r] \setminus \{ i,j\}$ and $|\phi(ij)| \leq 1$;
\item a \emph{strong clone of $j$ (under $\phi$)} if additionally $\phi(ij) = \emptyset$.
\end{itemize}
We say that $\bm{k}$ has
\begin{itemize}
\item the \emph{extension property} if, for every $(r^*,\phi^*,\ba^*) \in \opt^*(\bm{k})$ and $\phi \in \Phi_0(r^*+1;\bm{k})$ such that $\phi|_{\binom{[r^*]}{2}} = \phi^*$ and $\sum_{i \in [r^*]: \phi(\{i,r^*+1\})\neq \emptyset}\alpha_i\log|\phi(\{ i,r^*+1 \})| = Q(\bm{k})$,
there exists $j \in [r^*]$ such that $r^*+1$ is a clone of $j$ under $\phi$;
\item the \emph{strong extension property} if in fact $r^*+1$ is a strong clone of $j$.
\end{itemize}
\end{definition}

We explain the intuition and motivation behind this property in Section~\ref{extsec}.
For now we remark that it is generally easy to check whether $\bm{k}$ has the extension property, when $\opt^*(\bm{k})$ is known.
We check it for some cases in Section~\ref{applicationsec}.
For $\bm{k}$ with the extension property, one can describe all solutions to $\opt_0(\bm{k})$ in terms of basic optimal solutions.
Every solution can be obtained by `blowing up' a basic optimal solution $(r^*,\phi^*,\ba^*)$; that is, taking arbitrarily many clones of the vertices
and modifying part sizes so that the sum of vertex weights of clones of $j$ equals $\alpha_j^*$ for every $j$; and then possibly adding colour $c$ edges between clones of each $j$, without creating a $c$-coloured copy of $K_{k_c}$, 
where $c$ is the colour with the largest forbidden clique. Without loss of generality,
we assume $k_1 \geq \ldots \geq k_s$, so that $c=1$.
If $k_1=k_2$ then one cannot add any colour $1$ edges between clones without
creating a forbidden clique.

\begin{lemma}\label{char}
Let $s \in \mathbb{N}$ and suppose that $\bm{k} \in \mathbb{N}^s$ with $k_1 \geq \ldots \geq k_s$ has the extension property.
Let $(r,\phi,\ba) \in \feas_0(\bm{k})$.
Then $(r,\phi,\ba) \in \opt_0(\bm{k})$ if and only if there exist $(r^*,\phi^*,\ba^*) \in \opt^*(\bm{k})$ 
and a partition $V_1 \cup \ldots \cup V_{r^*}$ of $\{i: \alpha_i>0\}$
such that the following hold.
\begin{enumerate}[label=(\roman*),ref=(\roman*)]
\item\label{chari} For all $j \in [r^*]$, $\alpha^*_j = \sum_{i \in V_j}\alpha_i$.
\item\label{charii} For all $ij \in \binom{[r^*]}{2}$, $i' \in V_i$ and $j' \in V_j$ we have that $\phi(i'j') = \phi^*(ij)$.
\item\label{chariii} For all $i \in [r^*]$ and distinct $i',j' \in V_i$ we have that $\phi(i'j') \subseteq \{ 1 \}$. Moreover, if at least one $\phi(i'i'')$ for distinct $i',i''\in V_i$ and $i \in [r^*]$ is non-empty, then $k_1>k_2$ and there is an integer vector $\bm{\ell} \in \mathbb{N}^{r^*}$ such that $\|\bm{\ell}\|_1 \leq k_1-1$ and $\phi^{-1}(1)[V_i]$ is $K_{\ell_i+1}$-free for all $i \in [r^*]$.
\end{enumerate}
\end{lemma}

Our main result is the following stability theorem.
The \emph{edit distance} $d_{\rm edit}(G,G')$ of two graphs $G,G'$ of the same order is the minimum number of edges that need to be added/removed to make $G'$ isomorphic to $G$.
Given graphs $G$ and $H$ and $\delta >0$, we say that $G$ is \emph{$\delta$-far from being $H$-free} if it has edit distance at least $\delta|V(G)|^2$ to every $H$-free graph with the same number of vertices (note that we only need to delete edges here).
Given disjoint $A,B \subseteq V(G)$ and $0 \leq d \leq 1$, we say that $G[A,B]$ is \emph{$(\delta,d)$-regular} if $d_G(A,B) := e_G(A,B)|A|^{-1}|B|^{-1} \in (d-\delta,d+\delta)$, and $|d_G(X,Y)-d_G(A,B)|<\delta$ for all $X\subseteq A$, $Y \subseteq B$ with $|X|/|A|, |Y|/|B| \geq \delta$.
We are now ready to state the stability theorem.
It says that, for any large $n$-vertex graph $G$ which has close to the maximum number of valid colourings; that is, $F(G;\bm{k})=F(n,\bm{k}) \cdot 2^{o(n^2)}$, 
for almost all of its valid colourings $\chi$, there is a solution to $\opt_0(\bm{k})$ which describes the structure of $\chi$: it looks like a `blow-up' of the solution.
Lemma~\ref{char} describes the structure of solutions in terms of basic optimal solutions, 
and therefore there is a basic optimal solution $(r^*,\phi^*,\ba^*)$ which describes the structure of $\chi$.
Part~\ref{stabilitysimpii} implies that, not only is there a partition of $V(G)$ such that $\chi$ has many edges of every colour $c \in \phi^*(ij)$ between the $i$-th and $j$-th parts (and few edges of any other colour), these edges are in fact well-distributed and of roughly equal densities between these parts.

\begin{theorem}[Stability]\label{stabilitysimp}
Let $s \in \mathbb{N}$ and suppose that $\bm{k} \in \mathbb{N}^s$ with $k_1\geq \ldots \geq k_s$ has the extension property.
Then for all $\delta > 0$ there exist $n_0 \in \mathbb{N}$ and $\eps > 0$ such that the following holds.
If $G$ is a graph on $n \geq n_0$ vertices such that
$$
\frac{\log F(G;\textbf{k})}{\binom{n}{2}} \geq Q(\textbf{k}) - \eps,
$$
then, for at least $(1-2^{-\eps n^2})F(G;\bm{k})$ colourings $\chi : E(G) \rightarrow [s]$ which are $\bm{k}$-valid, there are $(r^*,\phi^*,\ba^*) \in \opt^*(\bm{k})$ and a partition $Y_1 \cup \ldots \cup Y_{r^*}$ of $V(G)$ such that the following hold.
\begin{enumerate}[label=(\roman*),ref=(\roman*)]
\item\label{stabilitysimpi} For all $i \in [r^*]$, we have that $|\,|Y_i| - \alpha^*_in\,| < 1$.
\item\label{stabilitysimpii} For all $c \in \phi^*(ij)$ and $ij \in \binom{[r^*]}{2}$, we have that $\chi^{-1}(c)[Y_i,Y_j]$ is $(\delta,|\phi^*(ij)|^{-1})$-regular. In particular, $e_G(Y_i,Y_j) \geq (1-s\delta)|Y_i||Y_j|$.
\item\label{stabilitysimpiii} Suppose $\sum_{i \in [r^*]}e(G[Y_i]) > \delta n^2$. Then
$\bm{k}$ does not have the strong extension property, and
all but at most $\delta n^2$ edges in $\bigcup_{i \in [r^*]}G[Y_i]$ are coloured with $1$ under $\chi$.
Moreover if $\bm{\ell} := (\ell_1,\ldots,\ell_{r^*}) \in \mathbb{N}^{r^*}$ is such that $G[Y_i]$ is $\delta$-far from being $K_{\ell_i}$-free, then
$\|\bm{\ell}\|_1 \leq k_{1}-1$.
\end{enumerate}
\end{theorem}

Somewhat conversely, if~\ref{stabilitysimpi},~\ref{stabilitysimpiii} and $e_G(Y_i,Y_j)\geq (1-s\delta)|Y_i||Y_j|$ hold for some triple $(r^*,\phi^*,\ba^*) \in \opt^*(\bm{k})$, a partition $Y_1,\ldots,Y_{r^*}$ of an $n$-vertex graph $G$, and $\delta=o(1)$, and each $G[Y_i]$ is $K_{\ell_i+1}$-free for some vector $\bm{\ell}$ of $1$-norm at most $k_1-1$,
then $F(G;\bm{k})\geq 2^{(Q(\bm{k})-o(1))n^2/2}$.

One should note the similarities between the statements of Theorem~\ref{stabilitysimp} and Lemma~\ref{char}.
Indeed, this parallel shows that the gist of Theorem~\ref{stabilitysimp} is `near-extremal graphs look like blow-ups of solutions to $\opt^*(\bm{k})$'.
This is not quite true within parts, as here $G$ could be very far from a complete partite graph.
Note also that the partition $Y_1 \cup \ldots\cup Y_{r^*}$ may be different for different colourings $\chi$.

We illustrate these statements with the example $\bm{k} = (5,3)$. 
As it is not hard to show (or see Lemma~\ref{2col}), $\opt^*((5,3))$ consists of just one element, namely $(2,\phi^*,(\frac{1}{2},\frac{1}{2}))$, where $\phi^*(12):=\{1,2\}$.
Thus by Theorem~\ref{stabilitysimp} and the remark after it, the set of almost extremal graphs can be described as consisting of graphs that are $o(n^2)$-close to $T_2(n)$ with triangle-free graphs added into each part, or a $K_4$-free graph added into one part. Note that the partition $Y_1\cup Y_2$ of Theorem~\ref{stabilitysimp} may depend on the colouring.
For example, if $G$ is $T_4(n)$ with parts $V_1,\ldots,V_4$, then Theorem~\ref{stabilitysimp} gives that for a typical colouring, there are disjoint pairs $ij,\ell m$ such that almost all edges between $V_i$ and $V_j$  and between $V_\ell$ and $V_m$ are coloured with colour $1$; then $Y_1$ and $Y_2$ in Theorem~\ref{stabilitysimp} have to be $V_i\cup V_j$ and $V_\ell \cup V_m$, up to changing $o(n)$ vertices.

If $\bm{k}$ has the strong extension property, then $G$ is close to a complete multipartite graph by Theorem~\ref{stabilitysimp}.
So the pairs $(r^*,\ba^*)$ associated with the colourings specified by the theorem are close to each other up to a relabelling of colours.
We have the following corollary of Theorem~\ref{stabilitysimp} for $\bm{k}$ with the strong extension property.

\begin{corollary}\label{uniform}
Let $s \in \mathbb{N}$ and suppose that $\bm{k} \in \mathbb{N}^s$ with $k_1\geq\ldots \geq k_s$ has the strong extension property.
Then for all $\delta > 0$ there exists $n_0 \in \mathbb{N}$ and $\eps > 0$ such that the following holds.
Let $G$ be a graph on $n \geq n_0$ vertices such that
$$
\frac{\log F(G;\bm{k})}{\binom{n}{2}} \geq Q(\bm{k}) - \eps.
$$
Then there are $r^*,\ba^*$ and a partition V(G) = $V_1 \cup \ldots \cup V_{r^*}$ with $\left|\,|V_i|-\alpha^*_i n \,\right| < 1$ for all $i \in [r^*]$ such that
the edit distance between $G$ and $K[V_1,\ldots,V_{r^*}]$ is at most $\delta n^2$.
Moreover, for at least $(1-2^{-\eps n^2}) \cdot F(G;\bm{k})$ $\bm{k}$-valid $s$-edge-colourings $\chi$ of $G$, there exists $(r^*,\phi^*,\ba) \in \opt^*(\bm{k})$ where $\|\ba-\ba^*\|_1 \leq \delta$ such that $\chi^{-1}(c)[V_i,V_j]$ is $(\delta,|\phi^*(ij)|^{-1})$-regular for all $ij \in \binom{[r^*]}{2}$ and $c \in \phi^*(ij)$.
\end{corollary}

Recall from Theorem~\ref{oldresults}\ref{oldresultsi} that at least one extremal graph is complete multipartite.
The following conjecture was made in~\cite{psy}:

\begin{conjecture}\label{conj}
For every $\bm{k}$,
every $\bm{k}$-extremal graph is complete multipartite.
\end{conjecture}

In~\cite{ps3}, we will use Corollary~\ref{uniform} to prove an exact result for all $\bm{k}$ with the strong extension property: that for such $\bm{k}$ every large extremal graph is a complete multipartite graph with part ratios roughly $\alpha_1^*,\ldots,\alpha_{r^*}^*$ for some $\ba^*$ coming from a basic optimal triple $(r^*,\phi^*,\ba^*)$.
However it seems much harder to prove an exact result for $\bm{k}$ without the strong extension property, as Theorem~\ref{stabilitysimp} as well as the example $\bm{k}=(5,3)$ above show that in general there are many asymptotically extremal graphs which are far from complete multipartite.

\subsection{Applications}

Armed with Theorem~\ref{stabilitysimp}, to determine asymptotically $\bm{k}$-extremal graphs, one need only solve Problem $Q_2$ for $\bm{k}$, and then check for the extension property using the optimal solutions.

We apply our stability theorem to reprove stability for most of the cases in which $F(n;\bm{k})$ has already been determined, in a systematic fashion.
For this, it suffices to solve Problem $Q_2$ (which follows from these earlier works), and to prove the extension property. Proving the extension property is straightforward: there are $O(s^r)$ possible attachments of a new vertex to some $(r,\phi,\ba)$ so a computer could check these for reasonable $s,r$. Actually, when optimal solutions have a particular nice form, which they do in almost all known cases, one can reduce the problem to determining solutions to some simple exponential equation over the integers (see Lemma~\ref{numcheck}).
We cannot prove stability for $\bm{k}=(3,3,3,3,3)$, which \emph{does not} have the extension property, strong or otherwise.

\begin{theorem}\label{recover}
Each $\bm{k}$ among $(k,k)$, $(k,k,k)$, $(3,3,3,3)$, $(4,4,4,4)$ has the strong extension property, for all integers $k \geq 3$. Thus Theorem~\ref{stabilitysimp} applies to every such $\bm{k}$.
\end{theorem}

As an example, we solve the optimisation problem for $s=2$ and state a version of Theorem~\ref{stabilitysimp} to illustrate it.
Note that it would not be too difficult to prove stability for $s=2$ directly.

\begin{lemma}\label{2col}
Let $k \geq \ell \geq 3$ be positive integers, $\bm{k} := (k,\ell)$ and $\phi$ be the function on $\binom{[\ell-1]}{2}$ that assumes value $\{1,2\}$ for every pair. Then
$Q(\bm{k}) = 1 - \frac{1}{\ell-1}$ and
$\opt^*(\bm{k}) = \{ (\ell-1,\phi,\bm{u}) \}$ where $\bm{u}$ is uniform,
and $\bm{k}$ has the extension property. Moreover $\bm{k}$ has the strong extension property if and only if $k=\ell$.
\end{lemma}

We write $\omega(G)$ for the clique number of a graph $G$; that is, the size of its largest clique.

\begin{theorem}\label{2colcor}
Let $k \geq \ell \geq 3$ be integers.
For all $\delta>0$ there exist $n_0\in \mathbb{N}$ and $\eps>0$ such that the following holds.
Let $G$ be a graph on $n\geq n_0$ vertices such that $\log F(G;(k,\ell)) \geq (Q(k,\ell)-\eps)\binom{n}{2}$. Then there is a graph $G'$ which can be obtained from $G$ by modifying at most $\delta n^2$ adjacencies, and an equipartition $V(G')=A_1 \cup \ldots \cup A_{\ell-1}$ such that $G'[A_i,A_j]$ is complete for all distinct $i,j \in [\ell-1]$, and $\omega(G'[A_1])+\ldots+\omega(G'[A_\ell]) \leq k-1$.

Moreover, for at least $(1-2^{-\eps n^2})F(G;\bm{k})$ valid colourings $\chi$ of $G$, $\chi^{-1}(c)[A_i,A_j]$ is $(\delta,\frac{1}{2})$-regular for $c=1,2$ and all distinct $i,j \in [\ell-1]$ and all but at most $\delta n^2$ edges in $G[A_1]\cup \ldots \cup G[A_{\ell-1}]$ are coloured with colour $1$.
\end{theorem}

\subsection{A sketch of the proof of Theorem~\ref{stabilitysimp}}

Since the proof of Theorem~\ref{stabilitysimp} is quite involved, we give a fairly detailed sketch first.
Let $G$ be a graph on $n$ vertices such that $\log F(G;\bm{k}) \geq (Q(\bm{k})-\eps) \binom{n}{2}$.
There are several ingredients to the proof.

\noindent
\emph{Regularity lemma: $G$ is close to some nearly optimal $(r,\phi,\ba)$ in $\feas_0(\bm{k})$.}
The multicolour version of Szemer\'edi's regularity lemma was already used to prove the existing results on $F(n;\bm{k})$ in~\cite{abks,PY}.
(For definitions and statements related to the regularity lemma, see Section~\ref{regtools}.)
Given a $\bm{k}$-valid colouring $\chi$ of $G$, we obtain an equitable partition $U_1 \cup \ldots \cup U_r$ in which almost all pairs are regular in all colours in $[s]$.
Define a colour pattern $\phi : \binom{[r]}{2} \rightarrow 2^{[s]}$ by adding the colour $c$ to $\phi(ij)$ if $\chi^{-1}(c)[V_i,V_j]$ is regular, and has density that is not too small.
The embedding lemma (Lemma~\ref{embed}) implies that $\phi^{-1}(c)$ is $K_{k_c}$-free for all $c \in [s]$.
In this way, much of the information carried by $\chi$ is transferred to the tuple $\RL(\chi) := (r,\phi,\mathcal{U})$, where $\mathcal{U} := \{ U_1,\ldots,U_r \}$.

Of course, one still needs to prove that this process is in some sense reversible: that the structure of $G$ itself, as well as the structure of its colourings, can be recovered from $(r,\phi,\mathcal{U})$.
This may not always be the case: we could have chosen some pathological $\chi$ to generate $(r,\phi,\mathcal{U})$.
For example, in the case $\bm{k} = (3,3)$, the unique extremal graph is $T_2(n) = K_{\lfloor n/2 \rfloor, \lceil n/2 \rceil}$, but we could have chosen $\chi$ which colours every edge with colour $1$.
Then we cannot recover many further colourings from $(r,\phi,\mathcal{U})$.

For this reason, we only wish to consider tuples $(r,\phi,\mathcal{U})$ which are the image of many colourings; that is, some non-trivial proportion of all colourings.
Such a tuple is called \emph{popular}; and we think of colourings $\chi$ which map to this tuple as being good representatives of the set of all colourings of $G$.
Since, as we show in Proposition~\ref{mainlypop}, almost every colouring maps to a popular tuple, it suffices to fix a popular tuple $(r,\phi,\mathcal{U})$ and only consider colourings which map to this tuple.
Intuitively, all such colourings should be similar.

Let $\alpha_i := 1/r$ for all $i \in [r]$.
Then $(r,\phi,\ba) \in \feas_0(\bm{k})$.
So the regularity lemma allows us to pass from $G$ to a feasible solution to Problem~$Q_0$.
It turns out that since $(r,\phi,\mathcal{U})$ is popular, we have that
$$
q(\phi,\ba) \geq Q(\bm{k}) - 2\eps,
$$
and moreover that $G[U_i,U_j]$ is almost complete for all $ij \in \binom{[r]}{2}$ (see Claim~\ref{largecon}).
Since we can choose $r$ large (but still bounded), the number of edges of $G$ within any $U_i$ can be made very small.
Therefore the structure of $G$ can be recovered from $(r,\phi,\mathcal{U})$.

\emph{Symmetrisation: from $(r,\phi,\ba) \in \feas_0(\bm{k})$ to some $(r^*,\phi^*,\ba^*) \in \opt^*(\bm{k})$.}
This is the main part of the proof (Lemma~\ref{bulk}), and in it we forget about $G$ entirely, and instead concentrate on $(r,\phi,\ba)$.
We think of this object as a vertex-weighted coloured multigraph: the weights are given by $\ba$, and the coloured edges by $\phi$.
We will apply a version of symmetrisation to $(r,\phi,\ba)$.
Symmetrisation was originally used in (ordinary) graphs by Zykov~\cite{zykov} to give a new proof of Tur\'an's theorem.
In its most basic form, it is the process of considering two non-adjacent vertices $x$ and $y$ in a graph $G$, and replacing $x$ by a clone of $y$, i.e.\ a vertex $y'$ whose neighbourhood is the same as that of $y$.
With Yilma~\cite{psy} we used symmetrisation to modify any $\bm{k}$-extremal graph into one which is both extremal and complete multipartite (Theorem~\ref{oldresults}\ref{oldresultsi}).
Here, we use a version of symmetrisation as follows.
Suppose that there is some $ij \in \binom{[r]}{2}$ such that $|\phi(ij)| \leq 1$.
Then we create a new feasible solution on $r$ parts by making vertex $j$ a clone of vertex $i$, or vice versa.
One of these choices will be such that the new solution $(r,\phi',\ba)$ satisfies $q(\phi',\ba) \geq q(\phi,\ba)$.
At the end of this process, we will obtain $(r,\phi_f,\ba) \in \feas_0(\bm{k})$ (where $f$ is for \emph{final}) such that $|\phi_f(ij)| \geq 2$ whenever $i,j$ are not clones and
$$
q(\phi_f,\ba) \geq Q(\bm{k}) - 2\eps.
$$
(In fact we split each step into small steps to get slowly evolving colour patterns $\phi=\phi_0,\ldots,\phi_f$.)
The solution $(r,\phi_f,\ba)$ corresponds to a smaller solution $(r_f,\psi_f,\ba_f)$ in which all clones are merged, so $(r_f,\psi_f,\ba_f) \in \feas_2(\bm{k})$ (and $q(\psi_f,\ba_f)=q(\phi_f,\ba)$ is near-optimal).
A compactness argument (Lemma~\ref{Sfnear}) shows that there is some basic optimal solution $(r^*,\phi^*,\ba
^*)$ which is very close to the near-optimal $(r_f,\phi_f,\ba_f)$ in a very strong sense: $\ba^*$ and $\ba_f$ are close in $\ell^1$-distance, and $\phi_f$ equals $\phi^*$ between pairs with non-negligible weights.

\emph{The extension property: $(r,\phi,\ba)$ and $(r^*,\phi^*,\ba^*)$ are close.}
Here, we mean `close' in the sense of Lemma~\ref{char}.
So we would like to show that when we merge pairs $ij$ with $|\phi(ij)| \leq 1$ in $(r,\phi,\ba)$ we obtain a weight vector $\bm{u}$ which is close to $\ba^*$ in $\ell^1$-distance, and $\phi \subseteq \phi^*$ on pairs of non-negligible size. (This turns out to be a simplification; see below.)
It is of course far from clear that we have not changed $(r,\phi,\ba)$ drastically to obtain $(r^*,\phi^*,\ba^*)$.
That this is not so is effectively a consequence of the extension property.

Having obtained $(r^*,\phi^*,\ba^*)$ from symmetrisation, we can now, with the benefit of hindsight, follow the procedure backwards.
At each stage, we check whether the attachments of each vertex are large; if not we sequentially put bad vertices into an exceptional set, called $U_{i}^0$ at the $i$-th step.
At each stage, the extension property implies that every vertex $x$ is either in $U_{i}^0$ or it corresponds to some vertex $k$ in the optimal solution $(r^*,\phi^*,\ba^*)$; that is, $\phi_i(xy) \subseteq \phi^*(kj)$ for all $j \in [r^*] \setminus \{ k \}$ and all vertices $y$ corresponding to $j$.
Let $U_{i}^k$ be the set of vertices corresponding to $k$ at Step~$i$.
Note that $|\phi_i(xy)| \leq 1$ if $x,y \in U_{i}^k$.

So, going back through the procedure, at each step there are a small fraction of vertices $x$ which were exceptional but are no longer exceptional, and these are moved to the $U_{i}^k$ for which $x$ corresponds to $k$;
and there are a small fraction of vertices moved into $U_{i}^0$.
When we return to Step~0, we want to show that
$\phi \subseteq \phi^*$ between any two classes $U_0^j$ and $U_0^k$, even though
the attachment between $U_{0}^0$ and the rest of the vertices, as well as the colour pattern within a class $U_{i}^k$ may have changed.
Therefore, if we can show that $U_{0}^0$ is small (Claim~\ref{cl-Ui0}) and every $U_{0}^k$, $k>0$, is about the same size throughout, then the procedure did not really change $(r,\phi,\ba)$ much at all.
Furthermore, we have a partition $U_{0}^0,U_{0}^1,\ldots,U_{0}^{r^*}$ of $[r]$ such that $\phi \subseteq \phi^*$ between $U_{0}^i,U_{0}^j$, $i,j \neq 0$; and $|\phi(ij)| \leq 1$ within a class.
In fact, the above sketch is a simplification, and though the exceptional set is always small, the other parts could change significantly in size. Nevertheless, we show that for each $i$ there are $\tilde{r}_i$ nonzero nonexceptional parts with part ratios given by some $\tilde{\ba}_i$, and $(\tilde{r}_i,\tilde{\phi}_i,\tilde{\ba}_i) \in \opt^*(\bm{k})$, where $\phi_i\subseteq\tilde{\phi}_i$ between pairs and $|\phi_i(jk)| \leq 1$ within a class.

\emph{Recovering $G$ from $(r^*,\phi^*,\ba^*)$.}
We now transfer the information we have gleaned about $(r,\phi,\ba)$ back to $G$ itself.
The partition of $[r]$ induces a partition $X_0,X_1,\ldots,X_{r^*}$ on $V(G)$.
For every $\chi \in \RL^{-1}((r,\phi,\mathcal{U}))$, we almost always have $\chi(xy) \in \phi^*(ij)$ when $ij \in \binom{[r^*]}{2}$ and $x \in X_i, y \in X_j$.
In Claims~\ref{badprimer} and~\ref{bad}, we show for almost all these $\chi$, that $\chi^{-1}(c)[X_i,X_j]$ is regular.
Indeed, if $\chi^{-1}(c)[X_i,X_j]$ is not regular, then Lemma~\ref{badred} implies that there is $c^* \in [s]$ and large sets $X \subseteq X_i$ and $Y \subseteq X_j$ between which the density of colour $c^*$ edges is significantly less than it is between $X_i$ and $X_j$.
But Corollary~\ref{binbound} implies that very few $s$-edge colourings of $G[X_i,X_j]$ are such that some colour class ($c \in \phi^*(ij)$) has size much less than the average $|\phi^*(ij)|^{-1}$.

The final step is to look inside the classes.
We discretise by applying the regularity lemma again within each class (which is large by Lemma~\ref{solutions}).
This allows us to apply Lemma~\ref{nocap} to prove~\ref{stabilitysimpiii}.

\subsection{Organisation}
Most of the paper concerns optimal triples for Problem $Q_t$ rather than graphs.
In Section~\ref{optsols} we collect some tools concerning these optimal triples, and some consequences of the extension property.
The main result in Section~\ref{staboptsols} is the stability of optimal solutions given the extension property which is the key component of the proof of our graph stability result, Theorem~\ref{stabilitysimp}.
In Section~\ref{staboptgraphs} we transfer statements on optimal triples to graphs, and prove Theorem~\ref{stabilitysimp}, after stating and proving some tools concerning the regularity lemma.
Then we prove Corollary~\ref{uniform}.
Section~\ref{applicationsec} concerns applications of Theorem~\ref{stabilitysimp}, in particular it contains the derivations of Theorems~\ref{recover} and~\ref{2colcor}.
We finish with some concluding remarks in Section~\ref{conclude}.

\section{Optimal solutions of Problem $Q_t$}\label{optsols}

We will now explore some properties of optimal solutions, which will be useful in later sections.
The following proposition states that, in an optimal solution $(r,\phi,\ba) \in \opt_0(\bm{k})$, every vertex $i$ of positive weight \emph{attaches optimally}, in that the normalised contribution $q_i(\phi,\ba)$ to the sum
$$
q(\phi,\ba) = \sum_{i \in [r]}\alpha_i q_i(\phi,\ba)\quad\text{where}\quad
q_i(\phi,\ba) := \sum_{\substack{j \in [r] \setminus \{ i \}\\ \phi(ij)\neq\emptyset}}\alpha_j\log|\phi(ij)|
$$
is equal to $Q(\bm{k})$.

\begin{proposition}\label{lagrange}
Let $s \in \mathbb{N}$ and $\bm{k} \in \mathbb{N}^s$ and suppose that $(r,\phi,\ba) \in \opt_0(\bm{k})$.
For every $i \in [r]$ with $\alpha_i > 0$, we have that
$
q_i(\phi,\ba) = Q(\bm{k})
$.
\end{proposition}

\begin{proof}
Without affecting the statement, we can remove all indices $i$ with $\alpha_i=0$.
So assume that $\alpha_i>0$ for all $i \in [r]$.
We use the method of Lagrange multipliers.
Recall that the constraint is that $g(\ba) = 0$, where $g(\ba) := \|\ba\|_1 - 1$.
Fix $\phi \in \Phi_0(r;\bm{k})$ and let
$$
\mathcal{L}(\ba,\lambda) := q(\phi,\ba) - \lambda g(\ba).
$$
Since the optimal vertex weighting $\ba$ is in the interior of $\Delta^r$ (and $\mathcal{L}$ is continuously differentiable there), there is $\lambda$ such that $(\ba,\lambda)$ is a critical point of $\mathcal{L}$. Thus
$$
\frac{\partial \mathcal{L}}{\partial \alpha_i} = 2q_i(\phi,\ba) - \lambda = 0\ \ \text{ and }\ \ \frac{\partial \mathcal{L}}{\partial \lambda} = \|\ba\|_1-1 = 0\quad\text{for all }i \in [r].
$$
We see that $q_i(\phi,\ba)$, $i \in [r]$, are equal to each other and the common value is $\lambda/2=\|\ba\|_1\lambda/2 = \sum_{i \in [r]}\alpha_iq_i(\phi,\ba)=q(\phi,\ba)=Q(\bm{k})$.
Therefore every $(r,\phi,\ba) \in \opt_0(\bm{k})$ satisfies the equation in the statement of the proposition.
\end{proof}

The next proposition shows that the objective function $q(\phi,\cdot)$ is Lipschitz continuous.

\begin{proposition}\cite[Proposition 11]{psy}\label{continuity}
Let $s,r \in \mathbb{N}$ and $\bm{k} \in \mathbb{N}^s$.
Let $\phi \in \Phi_0(r;\bm{k})$ and $\ba,\bb \in \Delta^r$.
Then
$$
|q(\phi,\ba)-q(\phi,\bb)| < 2(\log s)\|\ba-\bb\|_1.
$$
\end{proposition}

The next lemma states that whenever we have a feasible solution $(r,\phi,\ba) \in \feas_1(\bm{k})$ which is almost optimal, there is some vertex weighting $\ba^*$ which is close to $\ba$ such that $(r,\phi,\ba^*)$ is an optimal solution.

\begin{lemma}\cite[Claim~15]{psy}\label{Sfnear}
Let $s \in \mathbb{N}$ and $\bm{k} \in \mathbb{N}^s$.
For all $\delta > 0$ there exists $\eps > 0$ such that the following holds.
Let $(r,\phi,\ba) \in \feas_1(\bm{k})$ be such that $q(\phi,\ba) \geq Q(\textbf{k}) - \eps$.
Then there exists $\ba^* \in \Delta^{r}$ such that $\| \ba-\ba^* \|_1 < \delta$ and $(r,\phi,\ba^*) \in \opt_1(\bm{k})$.
\end{lemma}

Consider the following definition.

\begin{definition}[Capacity]\label{cap}
\rm
Given $r \in \mathbb{N}$, a graph $G$ with vertex set $[r]$, and positive integers $\ell_1,\ldots,\ell_r$, we write $(\ell_1,\ldots,\ell_r)G$ to denote the graph obtained from $G$ by, for each $i \in [r]$, replacing vertex $i$ with a clique $K_{\ell_i}$, and joining every vertex in $K_{\ell_i}$ to all vertices in the
cliques $K_{\ell_j}$ such that $j$ is a neighbour of $i$ in $G$.

Let $\Capa(G,k)$ be the set of those $(\ell_1,\ldots,\ell_r) \in \mathbb{N}^r$ for which $(\ell_1,\ldots,\ell_r)G$ is $K_k$-free.
\end{definition}

Since $(1,\ldots,1)G = G$, certainly $\Capa(G,k) \neq \emptyset$ whenever $G$ is $K_k$-free.
It is easy to see that, if $\bm{b} \in \mathbb{N}^r$ satisfies $\bm{b} \in \Capa(G,k)$, then $\bm{a} \in \Capa(G,k)$ for all $\bm{a} \in \mathbb{N}^r$ with $\bm{a} \leq \bm{b}$.
We think of $\Capa(G,k)$ as being a measure of how far a $K_k$-free graph $G$ is from containing a copy of $K_k$.
For example, $(1,2,2) \in \Capa(K_3,6)$ but $\Capa(K_3,4) = \{ (1,1,1) \}$.

We now prove some facts about the capacity of $(\phi^*)^{-1}(c)$ in a basic optimal solution $(r^*,\phi^*,\ba^*)$,
i.e.~$q(\phi^*,\ba^*)=Q(\bm{k})$, $|\phi^*(ij)| \geq 2$ for all $ij \in \binom{[r^*]}{2}$, and $\alpha_i^*>0$ for all $i \in [r^*]$.
The most important part is~\ref{nocapii}, which is essentially of a consequence of the fact that maximally $K_k$-free graphs on $r$ vertices have non-trivial capacity if and only if $r<k-1$.
This lemma will be important in proving Theorem~\ref{stabilitysimp}\ref{stabilitysimpiii}.

\begin{lemma}\label{nocap}
Let $s \in \mathbb{N}$ and $\bm{k} \in \mathbb{N}^s$ with $k_1\geq \ldots \geq k_s$.
Let
$(r^*,\phi^*,\ba^*) \in \opt^*(\bm{k})$ and define $J_c := ([r^*],(\phi^*)^{-1}(c))$.
Then the following statements hold.
\begin{enumerate}[label=(\roman*),ref=(\roman*)]
\item\label{nocapi} $J_c$ is maximally $K_{k_c}$-free.
\item\label{nocapii}
With $\bm{1} \in \mathbb{N}^{r^*}$ denoting the all-$1$ vector, we have
$$
\Capa(J_c,k_c) = \begin{cases} 
\{ \bm{1} \} \cup \{ \bm{\ell} \in \mathbb{N}^{r^*} : \|\bm{\ell}\|_1 \leq k_c-1 \} &\mbox{if } c = 1 \text{ and } k_1 > k_2 \\
\{ \bm{1} \} &\mbox{otherwise;}
\end{cases}
$$
furthermore, if $\Capa(J_1,k_1) \neq \{ \bm{1} \}$, then $J_1 \cong K_{r^*}$.
\item\label{nocapiii} $r^* \geq k_{2}-1$.
\end{enumerate}
\end{lemma}

\begin{proof}
Suppose that~\ref{nocapi} does not hold for some $c \in [s]$.
Then there exist distinct $i',j' \in [r^*]$ such that $i'j' \notin E(J_c)$, and $J_c \cup \{ i'j' \}$ is $K_{k_c}$-free.
Let $\phi' : \binom{[r^*]}{2} \rightarrow 2^{[s]}$ be defined by setting $\phi'(ij) := \phi^*(ij)$ whenever $ij \neq i'j'$; and setting $\phi'(i'j') := \phi^*(ij) \cup \{ c \}$.
By construction, $\phi' \in \Phi_2(r^*;\bm{k})$.
So $(r^*,\phi',\ba^*) \in \feas_2(\bm{k})$. But
$$
q(\phi',\ba^*) - q(\phi^*,\ba^*) = 2\alpha^*_{i'}\alpha^*_{j'}\log\left(1+\frac{1}{|\phi^*(i'j')|}\right) \geq 2\alpha^*_{i'}\alpha^*_{j'}\log \left( 1 + \frac{1}{s} \right) > 0,
$$
a contradiction. This proves~\ref{nocapi}.

We now prove~\ref{nocapii}.
To do this, we need the following simple claim.

\begin{claim}\label{maxfree}
Let $k \in \mathbb{N}$ and let $H$ be maximally $K_k$-free. Then every $x \in V(H)$ lies in a copy of $K_{k-1}$ if and only if $|V(H)| \geq k-1$.
\end{claim}

\bcpf
To prove the claim, note that the `only if' direction is trivial.
Suppose now that the other direction does not hold; that is, $r := |V(H)| \geq k-1$ and there exists $x \in V(H)$ which does not lie in a copy of $K_{k-1}$.
First consider the case when $N_H(x) = V(H) \setminus \{ x \}$.
Then $H-x$ is a $K_{k-2}$-free graph on at least $k-2$ vertices.
Thus there is a non-edge $e$ in $H-x$. Let $H' := H \cup \{ e \}$.
Then $H'-x$ is $K_{k-1}$-free, and since $x$ not does lie in a $K_{k-1}$ in $H$, $x$ does not lie in a copy of $K_k$ in $H'$.
Therefore $H'$ is $K_k$-free, a contradiction.

Consider now the case when there is some $y \in V(H-x)$ such that $xy \notin E(H)$.
Then $H'' := H \cup \{ xy \}$ is $K_k$-free, since any clique which lies in $H''$ but not $H$ must contain $x$.
Again, this is a contradiction, proving the claim.
\ecpf

\medskip
\noindent
Suppose that $\Capa((\phi^*)^{-1}(c),k_c) \neq \{ \bm{1} \}$.
Then there exists $j \in [r^*]$ such that $\bm{1}+\bm{e}_j \in \Capa((\phi^*)^{-1}(c))$.
Observe that the graph $(\bm{1}+\bm{e}_j)J_c$ is obtained from $J_c$ by inserting a twin $j'$ of $j$ and adding the edge $jj'$.
If $r^* \geq k_c-1$, then $j$ lies in a copy of $K_{k_c-1}$ in $J_c$ by Part~\ref{nocapi} and the claim.
So $j'$ together with the vertices in this copy form a copy of $K_{k_c}$ in $(\bm{1}+\bm{e}_j)J_c$, a contradiction.
So $r^* \leq k_c-2$.

Suppose instead that $\Capa((\phi^*)^{-1}(c),k_c) = \{ \bm{1} \}$.
Then $(\bm{1}+\bm{e}_j)J_c$ contains a copy of $K_{k_c}$ for every $j \in [r^*]$, and this copy necessarily contains $j$ (since $J_c$ itself is $K_{k_c}$-free).
Trivially $r^* = |V(J_c)| \geq k_c-1$.
We have proved that
\begin{equation}\label{capstate}
\Capa((\phi^*)^{-1}(c),k_c) = \{ \bm{1} \} \ \text{ if and only if } \ r^* \geq k_c-1.
\end{equation}

Let $C := \{ c \in [s] : r^* \leq k_c-2 \}$.
If $C = \emptyset$, then $\Capa((\phi^*)^{-1}(c),k_c) = \{ \bm{1} \}$ for all $ c \in [s]$.
Note also that $\{ \bm{\ell} \in \mathbb{N}^{r^*} : \|\bm{\ell}\|_1 \leq k_c - 1 \} \subseteq \{ \bm{1} \}$.
So we are done in this case, and may assume that $C \neq \emptyset$.

By Part~\ref{nocapi}, for all $c \in C$ we have $(\phi^*)^{-1}(c) \cong K_{r^*}$ since this is the unique maximally $K_{k_c}$-free graph on $r^*$ vertices.
Define a new solution on $r^*+1$ vertices as follows.
Suppose, without loss of generality, that $\alpha_{r^*} \geq \alpha_i$ for all $i \in [r^*]$.
Let
\begin{align*}
\phi(ij) &:= \begin{cases} 
\phi^*(ij) &\mbox{if } ij \in \binom{[r^*]}{2} \\
\phi^*(ir^*) &\mbox{if } i \in [r^*-1], j = r^*+1\\
C &\mbox{if } \{ i,j \} = \{ r^*,r^*+1 \}
\end{cases}
\end{align*}
and
\begin{align*}
\alpha_i &:= \begin{cases} 
\alpha^*_i &\mbox{if } i \in [r^*-1] \\
\alpha^*_{r^*}/2 &\mbox{if } i \in \{ r^*,r^*+1 \}.
\end{cases}
\end{align*}
Then $\phi^{-1}(c)$ is $K_{k_c}$-free for all $c \in [s]$, so $(r^*+1,\phi,\ba) \in \feas_1(\bm{k})$.
Furthermore,
$$
0 \geq q(\phi,\ba)-q(\phi^*,\ba^*) = 2\left(\frac{\alpha_{r^*}^*}{2}\right)^2 \log |C|
$$
and so $|C|=1$.
Let $C = \{ c^* \}$.
Then $k_{c^*}-2 \geq r^* \geq k_c -1$ for all $c \in [s] \setminus \{ c^* \}$.
That is, $c^*=1$ and $k_1 > k_2$.
Suppose that $(\ell_1,\ldots,\ell_{r^*}) \in \mathbb{N}^{r^*}$.
Since $(\phi^*)^{-1}(1) \cong K_{r^*}$, the graph $(\ell_1,\ldots,\ell_{r^*})(\phi^*)^{-1}(1)$ is a clique of order $r^* + \sum_{i \in [r^*]}(\ell_i-1) = \sum_{i \in [r^*]}\ell_i$.
Therefore $\Capa((\phi^*)^{-1}(1),k_1) = \{ \bm{\ell} \in \mathbb{N}^{r^*} : \|\bm{\ell}\|_1 \leq k_1-1 \}$ (note that this set contains $\bm{1}$).
This completes the proof of~\ref{nocapii}.

Finally, Part~\ref{nocapiii} is an immediate consequence of Part~\ref{nocapii} and (\ref{capstate}).
\end{proof}

\subsection{The extension property}\label{extsec}

In this section, we explore the consequences of the extension property.
Recall the equality from Proposition~\ref{lagrange} which is necessary for all vertices of positive weight in an optimal solution.
Suppose we wish to extend an optimal solution by adding a new vertex of zero weight.
The following proposition shows that the normalised contribution of this new vertex cannot be more than $Q(\bm{k})$.
Given $\phi\in \Phi_0(r+1,\bm{k})$ and $\ba \in \Delta^r$, let
$$
\ext(\phi,\ba) := q_{r+1}(\phi,(\alpha_1,\ldots,\alpha_r,0)) = \sum_{\substack{i \in [r]\\ \phi(\{i,r+1\}) \neq \emptyset}}\alpha_i\log|\phi(\{ i,r+1 \})|,
$$
which is the ``normalised contribution'' of the zero-weighted vertex $r+1$ to $q(\phi,(\alpha_1,\ldots,\alpha_r,0))$.

\begin{proposition}\label{extendbd}
Let $s \in \mathbb{N}$ and $\bm{k} \in \mathbb{N}^s$.
Suppose that $(r,\phi',\ba) \in \opt^*(\bm{k})$.
Let $\phi \in \Phi(r+1,\bm{k})$ be such that $\phi|_{\binom{[r]}{2}} = \phi'$.
Then $\ext(\phi,\ba) \leq Q(\bm{k})$.
\end{proposition}

\begin{proof}
Suppose not.
We will show that we can transfer a small amount of weight from $[r]$ to $r+1$ and in so doing increase $q(\phi,\cdot)$.
Let $\gamma>0$ satisfy $\ext(\phi,\ba) = (1+\gamma)Q(\bm{k})$.
Let $\eps \in (0,2\gamma/(2\gamma+1))$.
Define $\bb \in \mathbb{R}^{r+1}$ by setting
$$
\beta_i := \begin{cases} 
(1-\eps)\alpha_i &\mbox{if } i \in [r] \\
\eps &\mbox{if }  i = r+1.
\end{cases}
$$
Then $\bm{\alpha} \in \Delta^r$ implies that $\bb \in \Delta^{r+1}$.
Now,
\begin{align*}
q(\phi,\bb) -q(\phi,\ba)
&= \eps(\eps-2) q(\phi',\ba) + 2(1-\eps)\eps \cdot \ext(\phi,\ba) = \eps(2\gamma-\eps(1+2\gamma))Q(\bm{k}) > 0,
\end{align*}
a contradiction.
\end{proof}

This motivates the \emph{extension property}, which we repeat for the reader's convenience:

\begin{definition}[Clones and extension property]\label{extprop2}
\rm
Let $s \in \mathbb{N}$ and $\bm{k} \in \mathbb{N}^s$. 
Given $r \in \mathbb{N}$ and $\phi \in \Phi_0(r;\bm{k})$, we say that $i \in [r]$ is
\begin{itemize}
\item a \emph{clone of $j \in [r]\setminus \{ i \}$} \emph{(under $\phi$)} if $\phi(ik) = \phi(jk)$ for all $k \in [r] \setminus \{ i,j\}$ and $|\phi(ij)| \leq 1$;
\item a \emph{strong clone of $j$ (under $\phi$)} if additionally $\phi(ij) = \emptyset$.
\end{itemize}
We say that $\bm{k}$ has
\begin{itemize}
\item the \emph{extension property} if, for every $(r^*,\phi^*,\ba^*) \in \opt^*(\bm{k})$ and $\phi \in \Phi_0(r^*+1;\bm{k})$ such that $\phi|_{\binom{[r^*]}{2}} = \phi^*$ and $\ext(\phi,\ba^*) = Q(\bm{k})$;
there exists $j \in [r^*]$ such that $r^*+1$ is a clone of $j$ under $\phi$;
\item the \emph{strong extension property} if in fact $r^*+1$ is a strong clone of $j$.
\end{itemize}
\end{definition}

The extension property says that if we extend any basic optimal solution by adding an infinitesimal part with optimal contribution $Q(\bm{k})$, then the new vertex clones an existing one (with perhaps one colour on the pair spanned by the two clones).
Assuming that $\bm{k}$ has the extension property, we can prove some properties of elements in $\opt^*(\bm{k})$, including a uniform lower bound for vertex weightings in $Q^*$-optimal solutions.

\begin{lemma}\label{solutions}
Let $s \in \mathbb{N}$ and suppose that $\bm{k} \in \mathbb{N}^s$ has the extension property.
Then there exists $\mu >  0$ such that $\alpha_i^* > \mu$ for all $(r^*,\phi^*,\ba^*) \in \opt^*(\bm{k})$ and $i \in [r^*]$.
\end{lemma}

\begin{proof}
Suppose not; then, for all $n \in \mathbb{N}$, there exists $(r^*_n,\phi^*_n,\ba^*_n) \in \opt^*(\bm{k})$ and $i_n \in [r^*_n]$ such that $\alpha^*_{i_n} < 1/n$.
By passing to a subsequence, since $r^*_n < R(\bm{k})$, we may assume that $r^*_n \equiv r$ and $\phi^*_n \equiv \phi$ and without loss of generality that $i_n \equiv r$.
Since $\ba^*_n \in \Delta^{r}$ and the simplex is closed and bounded, the Heine-Borel theorem implies that $\ba^*_1,\ba^*_2,\ldots$ has a convergent subsequence $\ba^*_{i_1},\ba^*_{i_2},\ldots$, with limit $\bb \in \Delta^{r}$.
Observe that $\beta_{r} = 0$.
Without loss of generality, assume that $\bb = (\beta_1,\ldots,\beta_t,0,\ldots,0)$, where $t \in [r-1]$ and $\beta_j > 0$ for all $j \in [t]$.
By continuity (Proposition~\ref{continuity}), 
$q(\phi,(\beta_1,\ldots,\beta_t))=q(\phi,\bb)=Q(\bm{k})$, so
$(t,\phi,(\beta_1,\ldots,\beta_t)) \in \opt^*(\bm{k})$.
Recall that $\alpha^*_{i_m,r} > 0$ for all $m \in \mathbb{N}$.
By continuity and Proposition~\ref{lagrange},
$$
\sum_{j \in [t]}\beta_j\log|\phi(rj)| = 
\ext(\phi,(\beta_1,\ldots,\beta_{r-1})) 
= \lim_{m\to\infty}q_r(\phi,(\alpha^*_{i_m,1},\ldots,\alpha^*_{i_m,r})) = Q(\bm{k}). 
$$
The extension property implies that there exists $i \in [t]$ such that $\phi(rj)=\phi(ij)$ for all $j \in [t] \setminus \{ i \}$ and $|\phi(ir)| \leq 1$, a contradiction to $\phi \in \Phi_2(r;\bm{k})$.
This completes the proof of the lemma.
\end{proof}

Next we prove that the strong extension property implies that optimal colour patterns have trivial capacity.

\begin{lemma}\label{strong}
Let $s \in \mathbb{N}$ and suppose that $\bm{k} \in \mathbb{N}^s$ with $k_1 \geq \ldots \geq k_s$ has the extension property.
\begin{enumerate}[label=(\roman*),ref=(\roman*)]
\item\label{strongi} If $\bm{k}$ has the strong extension property, then for every $(r^*,\phi^*,\ba^*) \in \opt^*(\bm{k})$ and $c \in [s]$, we have that $\Capa((\phi^*)^{-1}(c),k_c) = \{ \bm{1} \}$.
\item\label{strongii} If $k_1=k_2$, then $\bm{k}$ has the strong extension property.
\item\label{strongiii} If $(r^*,\phi^*,\ba^*) \in \opt^*(\bm{k})$ and $\phi \in \Phi_0(r^*+1;\bm{k})$ is such that $\phi|_{\binom{[r^*]}{2}}=\phi^*$ and $r^*+1$ is a clone of $i \in [r^*]$ under $\phi$, then $\phi(\{ i,r^*+1 \}) \subseteq \{ 1 \}$.
\end{enumerate}
\end{lemma}

\begin{proof}
For $c \in [s]$, write $C(c):=\Capa((\phi^*)^{-1}(c),k_c)$ as shorthand.
We use the following claim to prove all three parts.

\begin{claim}\label{capaclaim}
Let $(r^*,\phi^*,\ba^*) \in \opt^*(\bm{k})$, $c \in [s]$ and $j \in [r^*]$.
Define $\ba \in \Delta^{r^*+1}$ by setting $\alpha_i := \alpha^*_i$ for all $i \in [r^*]$, and $\alpha_{r^*+1} := 0$.
Define $\phi \in \Phi_1(r^*+1;\bm{k})$ by setting $\phi|_{\binom{[r^*]}{2}} := \phi^*$, and $\phi(\{ i,r^*+1 \}) := \phi^*(ij)$ for all $i \in [r^*]\setminus \{ j \}$, and $\phi(\{ j,r^*+1 \}) := \{ c \}$.
Then $\bm{1} + \bm{e}_j \in C(c)$ if and only if
$(r^*+1,\phi,\ba) \in \opt_1(\bm{k})$. 
\end{claim}

\bcpf
We need to show that $\bm{1}+\bm{e}_j \in C(c)$ if and only if both $(r^*+1,\phi,\ba)\in \feas_1(\bm{k})$ and $q(\phi,\ba)=Q(\bm{k})$.
Firstly, we have $(r^*+1,\phi,\ba) \in \feas_1(\bm{k})$ if and only if $\phi^{-1}(c)$ is $K_{k_c}$-free for all $c \in [s]$.
For $c'\in[s]\setminus\{c\}$, $\phi^{-1}(c')$ is obtained from $(\phi^*)^{-1}(c')$ by cloning vertex $r^*$, so is $K_{k_{c'}}$-free.
By definition, $\phi^{-1}(c)$ is $K_{k_c}$-free if and only if
$\bm{1} + \bm{e}_j \in C(c)$.
Secondly, $q(\phi,\ba)=q(\phi^*,\ba^*)=Q(\bm{k})$.
This proves the claim.
\ecpf

For~\ref{strongi}, suppose that $\bm{k}$ has the strong extension property but there is $(r^*,\phi^*,\ba^*)\in\opt^*(\bm{k})$ and $c \in [s]$ for which $C(c) \neq \{\bm{1}\}$.
Then there is $j \in [r^*]$ such that $\bm{1}+\bm{e}_j \in C(c)$. 
Let $\ba$ and $\phi$ be defined as in Claim~\ref{capaclaim}.
By the claim, $(r^*+1,\phi,\ba) \in \opt_1(\bm{k})$. So
$$
\sum_{i \in [r^*]}\alpha_i \log|\phi(\{ i,r^*+1\})| = \sum_{i \in [r^*]\setminus \{ j \}}\alpha^*_i\log|\phi^*(ij)| = Q(\bm{k})
$$
by Proposition~\ref{lagrange} applied to $(r^*,\phi^*,\ba^*)$.
But $r^*+1$ is not a strong clone of any $j' \in [r^*]$ under $\phi$ since $|\phi(\{r^*+1,j\})|=1$ (and $|\phi(\{r^*+1,i\})| = |\phi(ji)| \geq 2$ for all $i \in [r^*]\setminus\{j\}$).
So $\bm{k}$ does not have the strong extension property, a contradiction.

Next we prove~\ref{strongii}.
So suppose that $\bm{k}$ does not have the strong extension property. Then there is some $(r^*,\phi^*,\ba^*) \in \opt^*(\bm{k})$ and an extension $\phi\in\Phi_0(r^*+1;\bm{k})$ such that $\phi|_{\binom{[r^*]}{2}}=\phi^*$, $\ext(\phi,\ba^*)=Q(\bm{k})$, and
$r^*+1$ is a clone of some $j \in [r^*]$ under $\phi$, but not a strong clone.
So $\phi(\{i,r^*+1\})=\phi^*(ij)$ for all $i \in [r^*]\setminus \{j\}$, and $\phi(\{j,r^*+1\})=\{c\}$ for some $c \in [s]$.
Note that $(r^*+1,\phi,\ba)\in\opt_1(\bm{k})$, where $\ba$ is defined as in the claim. Thus, by the claim, $\bm{1}+\bm{e}_j \in C(c)$.
By Lemma~\ref{nocap}\ref{nocapii}, this implies $c=1$ and $k_1>k_2$.
This also gives Part~\ref{strongiii}.
\end{proof}

\subsection{The proof of Lemma~\ref{char}}

Recall that Lemma~\ref{char}, informally speaking, enables us to characterise all solutions to Problem $Q_0$ in terms of the basic optimal solutions $\opt^*(\bm{k})$.

\begin{proof}[Proof of Lemma~\ref{char}]
Note that the `if' direction is trivial so it remains to prove the `only if' direction.
Let $(r,\phi,\ba) \in \opt_0(\bm{k})$.
We can assume that $\alpha_i>0$ for all $i \in [r]$.

It is convenient to consider triples $(A,\phi,\ba)$ which are as feasible solutions $(r,\phi,\ba)$ except $A$ is a set of $r$ vertices (as opposed to $[r]$), $\phi : \binom{A}{2}\to 2^{[s]}$ and $\ba \in \Delta^A := \{(\alpha_i: i \in A): \alpha_i>0 \text{ for all } i \in A \text{ and }\sum_{i \in A}\alpha_i=1\}$.
Given $x,y \in A$, define a new vertex weighting $\overline{\ba} \in \Delta^{A\setminus\{y\}}$, the \emph{$(x,y)$-merging} of $\ba$,
 by setting $\overline{\alpha}_x := \alpha_x+\alpha_y$ and $\overline{\alpha}_z := \alpha_z$ for all $z \in [r] \setminus \{ x,y \}$.
Suppose $|\phi(xy)| \leq 1$. Then $(A\setminus\{y\},\phi',\overline{\ba}) \in \opt_0(\bm{k})$ where $\phi':=\phi|_{\binom{A\setminus\{y\}}{2}}$, and
\begin{align}
\label{merge}
q_y(\phi,\overline{\ba}) &= \sum_{z \in [r]\setminus \{ x,y \}}\alpha_z\log|\phi(zy)| + (\alpha_x+\alpha_y)\log|\phi(xy)|
= q_y(\phi,\ba)= Q(\bm{k}),
\end{align}
where the last equality follows from Proposition~\ref{lagrange}.

Consider the following claim.

\begin{claim2}
$\{ ij \in \binom{[r]}{2} : |\phi(ij)| \leq 1\}$ is a disjoint union of cliques.
\end{claim2}

\bcpf
Suppose for a contradiction to the claim that, without loss of generality, there is some $ij \in \binom{[r-1]}{2}$ such that $|\phi(ir)|,|\phi(jr)|\leq 1$ but $|\phi(ij)| \geq 2$.
Suppose first that there exists $i'j' \in \binom{[r]}{2}\setminus \{ir,jr\}$ such that $|\phi(i'j')| \leq 1$.
At least one of $i',j'$ is not in $\{i,j,r\}$, say $j'$.
Take the $(i',j')$-merging $\overline{\ba}$ of $\ba$.
By the above observations, $([r]\setminus\{j'\},\phi',\overline{\ba}) \in \opt_0(\bm{k})$, where $\phi':=\phi|_{\binom{[r]\setminus\{j'\}}{2}}$.
Note that $ir,jr,ij\in\binom{[r]\setminus\{j'\}}{2}$, and $\phi'(xy)=\phi(xy)$ for all $xy\in\binom{[r]\setminus\{j'\}}{2}$.

Do this repeatedly until the only pairs $i'j'$ with $|\phi(i'j')| \leq 1$ among the set $A$ of remaining vertices are $ir$ and $jr$.
Let $\bb$ be the weight function and $\psi := \phi|_A$ the colour pattern. We have $(A,\psi,\bb)\in\opt_0(\bm{k})$.
Now obtain the $(i,r)$-merging $\overline{\bb}$ of $\bb$ and let $A':=A\setminus\{r\}$.
By the above, $(A',\psi',\overline{\bb}) \in \opt_0(\bm{k})$ and $\psi,\overline{\bb},r$ satisfy~(\ref{merge}).
Further, $|\psi'(xy)|  = |\phi(xy)| \geq 2$ for every $xy \in \binom{A'}{2}$ and $\overline{\beta}_x>0$ for every $x \in A'$, so in fact $(A',\psi',\overline{\bb}) \in \opt^*(\bm{k})$. 
Since $\bm{k}$ has the extension property, there exists $y \in A'$ such that $\psi'(r\ell) = \psi'(y\ell)$ for all $\ell \in A'$, and $|\psi'(yr)|\leq1$.
In particular, $|\psi'(rj)| = |\psi'(yj)| \geq 2$, a contradiction to our assumption.
\ecpf

\medskip
\noindent
Proposition~\ref{lagrange} implies that, for every $i \in [r]$, we have that $q_i(\phi,\ba) = Q(\bm{k})$.
By the claim, there is a (unique up to relabelling) partition $[r] = V_1 \cup \ldots \cup V_{r^*}$ such that
\begin{equation}\label{almosteq}
\left\{ ij \in \binom{[r]}{2} : |\phi(ij)| \leq 1 \right\} = \bigcup_{j \in [r^*]}\binom{V_j}{2}
\end{equation}
(where a vertex $i'$ is the only member of some $V_{i}$ if and only if $|\phi(i'j)| \geq 2$ for all $j \in [r]\setminus \{ i' \}$). 
Assume without loss of generality that $i \in V_i$ for all $i \in [r^*]$.
Let $\ba^* \in \Delta^{r^*}$ such that $\alpha^*_i = \sum_{i' \in V_i}\alpha_{i'}$, and set $\phi^* := \phi|_{\binom{[r^*]}{2}}$.

\begin{claim2}\label{almost}
We have the following:
\begin{enumerate}[label=(\alph*),ref=(\alph*)]
\item\label{almostii} $(r^*,\phi^*,\ba^*) \in \opt^*(\bm{k})$.
\item\label{almostiii} For all $i \in [r^*]$ and $i' \in V_i$, we have that $\phi(i'j) = \phi(ij)$ for all $j \in [r^*]\setminus \{ i \}$.
\end{enumerate}
\end{claim2}

\bcpf
Let
$$
K := \bigcup_{j \in [r^*]}\{ (j,j') : j' \in V_j\setminus\{j\} \}
$$ 
 and $t := |K|$.
So $K$ is a union of spanning stars in the $V_j$'s.
We will form a new solution by transferring the total weight from $V_j$ to $j$.

Let $\ba_0 := \ba$, $\phi_0 := \phi$ and $A_0 := [r]$.
Order the elements $(j_1,x_1), \ldots, (j_t,x_{t})$ of $K$, and, for each $\ell \geq 1$, let $\ba_\ell$ be the $(j_\ell,x_\ell)$-merging of $\ba_{\ell-1}$ and $A_\ell := A_{\ell-1}\setminus\{x_\ell\}$ and $\phi_\ell := \phi|_{\binom{A_\ell}{2}}$.
Precisely as in (\ref{merge}), we have that $(A_\ell,\phi_\ell,\ba_\ell) \in \opt_0(\bm{k})$, and
\begin{equation}\label{iterate}
\sum_{k \in \binom{A_\ell}{2}}\alpha_{\ell,k}\log|\phi(x_\ell k)| = Q(\bm{k}).
\end{equation}
By construction, $A_t=[r^*]$ and $\alpha_{t,i} > 0$ for all $i \in [r^*]$.
Let $\ba^* := (\alpha_{t,1},\ldots,\alpha_{t,r^*})$ and $\phi' := \phi_t|_{\binom{[r^*]}{2}} = \phi|_{\binom{[r^*]}{2}}$.
Then $(r^*,\phi^*,\ba^*) \in \opt^*(\bm{k})$.
Moreover, by (\ref{iterate}), we have that
$
\sum_{i \in [r^*]}\alpha^*_i\log|\phi(\{ i,r^*+j \})|= Q(\bm{k})
$
for all $j \in [r-r^*]$, and
$
\alpha^*_i = \sum_{i' \in V_i}\alpha_{i'}
$.
So $\ba^*$ and $\phi^*$ satisfy~\ref{almostii}.

It remains to prove~\ref{almostiii}.
For each $i' \in \{ r^*+1,\ldots,r \}$, let $i \in [r^*]$ be such that $i' \in V_i$.
Apply the extension property to $(r^*,\phi^*,\ba^*)$ with $i'$ playing the role of the additional vertex, whose colour pattern is given by $\phi$.
So there is some $x_{i'} \in [r^*]$ which is a clone of $i'$ under $\phi$ in $[r^*]$.
But, by the definition of $V_i$, $|\phi(\{ k,i' \})| \leq 1$ if and only if $k \in V_i$.
But there is a unique member of $V_{i}$  which lies in $[r^*]$, namely $i$.
Certainly $i$ is a clone of itself.
So for all $i \in [r^*]$ and $i' \in [r]\cap V_i$, we have that $i'$ is a clone of vertex $i$ under $\phi$ in $[r^*]$.
So~\ref{almostiii} holds, completing the proof of the claim.
\ecpf

So Part~\ref{chari} of Lemma~\ref{char} holds by (\ref{almosteq}) and Claim~\ref{almost}\ref{almostii}. For~\ref{charii}, we need to prove that
$\phi(i'j') = \phi(ij)$ for all $i' \in V_i$ and $j' \in V_j$, whenever $i \neq j$.
Suppose that there is some $c \in \phi(i'j') \setminus \phi(ij)$.
Thus $c \notin \phi(ij)=\phi^*(ij)$, and
by Lemma~\ref{nocap}\ref{nocapi}, $(\phi^*)^{-1}(c)$ is maximally $K_{k_c}$-free, so there are vertices $\{ x_1,\ldots,x_{k_c-2} \} \in [r^*]\setminus \{ i,j \}$ which, together with $i,j$, span a copy of $K_{k_c}$ in $(\phi^*)^{-1}(c)$.
But Claim~\ref{almost}\ref{almostiii} implies that, for all $\ell \in [k_c-2]$, we have $c \in \phi^*(ix_\ell) = \phi(i'x_\ell)$ and $c \in \phi^*(jx_\ell) = \phi(j'x_\ell)$.
Therefore $x_1,\ldots,x_{k_c-2},i',j'$ span a copy of $K_{k_c}$ in $\phi^{-1}(c)$, a contradiction.
So $\phi(i'j') \subseteq \phi(ij)$.
Using Proposition~\ref{lagrange} and the fact that $|\phi(i'i'')| \leq 1$ for all $i'' \in V_i$, we have that
\begin{align*}
Q(\bm{k}) &= 
q_{i'}(\phi,\ba) = \sum_{j \in [r^*]\setminus \{ i \}}\sum_{j' \in V_j}\alpha_{j'}\log|\phi(i'j')| \leq \sum_{j \in [r^*]\setminus \{ i \}}\sum_{j' \in V_j}\alpha_{j'}\log|\phi(ij)|\\
&= 
q_i(\phi^*,\ba^*) = Q(\bm{k}).
\end{align*}
Therefore we have equality everywhere, and so $\phi(i'j') = \phi(ij) = \phi^*(ij)$ for all $ij \in \binom{[r^*]}{2}$ and $i' \in V_i$, $j' \in V_j$.
This completes the proof of~\ref{charii}.

For~\ref{chariii}, let $c \in [s]$, $i \in [r^*]$ and $i'i'' \in \binom{V_i}{2}$ with $c \in \phi(i'i'')$.
Then $\bm{1} + \bm{e}_i \in \Capa((\phi^*)^{-1}(c),k_c)$.
Lemma~\ref{nocap}\ref{nocapii} implies that $c = 1$ and $k_1>k_2$.
Now, for each $i \in [r^*]$, let $\ell_i$ be the size of the largest clique in $\phi^{-1}(1)[V_i]$.
By definition, $\bm{\ell} := (\ell_1,\ldots,\ell_{r^*}) \in \Capa((\phi^*)^{-1}(1),k_1)$, and so $\|\bm{\ell}\|_1 \leq k_1-1$ by Lemma~\ref{nocap}\ref{nocapii}.
This complete the proof of~\ref{chariii} and hence of the lemma.
\end{proof}

\subsection{Non-optimal attachments}

We derive a further quantifiable consequence of the extension property in the following lemma, which shows that if a basic optimal solution is extended by an infinitesimal part, if it is not a clone of an existing vertex then the deficit of its contribution is bounded away from zero.

\begin{lemma}\label{littlec}
Let $s \in \mathbb{N}$ and let $\bm{k} \in \mathbb{N}^s$ have the extension property.
Then there exists $\eta > 0$ such that the following holds.
Let $(r^*,\phi^*,\ba^*) \in \opt^*(\bm{k})$ and $\phi \in \Phi_0(r^*+1,\bm{k})$ such that $\phi|_{\binom{[r^*]}{2}}=\phi^*$ and $r^*+1$ is not a clone of any $i \in [r^*]$ under $\phi$.
Then
$\ext(\phi,\ba^*) \leq Q(\bm{k}) - \eta$.
\end{lemma}

\begin{proof}
Suppose that the statement of the lemma does not hold.
Then for all $n \in \mathbb{N}$, there exist $(r^*_n,\phi^*_n,\ba^*_n) \in \opt^*(\bm{k})$
and $\phi_n \in \Phi_0(r^*_n+1,\bm{k})$ such that $\phi_n|_{\binom{[r^*_n]}{2}} = \phi^*_n$ and $r^*_n+1$ is not a clone of any $i \in [r^*_n]$ under $\phi_n$, but
$
\ext(\phi_n,\ba_n^*) > Q(\bm{k}) - \frac{1}{n}
$.
By passing to a subsequence (since $r_n^*<R(\bm{k})$), we may assume that $r^*_n \equiv r$; $\phi^*_n \equiv \phi^*$ and $\phi_n \equiv \phi$, so
\begin{equation}\label{exteq}
\ext(\phi,\ba^*_n) > Q(\bm{k})-\frac{1}{n}
\end{equation}
and $r+1$ is not a clone of any $i \in [r]$ under $\phi$.
Since $\Delta^{r}$ is compact, we may choose a convergent subsequence $\ba^*_{i_1},\ba^*_{i_2},\ldots$ of $\ba^*_1,\ba^*_2,\ldots$, with limit $\bb$.
Now, since $q(\phi^*,\cdot)$ is continuous (by Proposition~\ref{continuity}),
$$
q(\phi^*,\bb) = \lim_{n \rightarrow \infty}q(\phi^*,\ba^*_{i_n}) = Q(\bm{k}),
$$
and Lemma~\ref{solutions} implies that $\beta_j =\lim_{n \rightarrow\infty}\alpha_{i_n,j}^* >0$ for every $j \in [r]$.
So $(r,\phi^*,\bb) \in \opt^*(\bm{k})$.
But taking the limit in (\ref{exteq}) implies that 
$\ext(\phi,\bb) = Q(\bm{k})$.
Now the extension property implies that there is some $i \in [r]$ such that $r+1$ is a clone of $i$ under $\phi$, a contradiction.
\end{proof}

Given colour patterns $\psi \in \Phi_0(r;\bm{k})$ and $\psi' \in \Phi_0(r';\bm{k})$ and a partition $\mathcal{V}=\{V_1,\ldots,V_r\}$ of $[r']$, we will say $\psi' =_{\mathcal{V}} \psi$ if $\psi'(i'j')=\psi(ij)$ for all $ij \in \binom{[r]}{2}$, $i' \in V_i$ and $j' \in V_j$, and $\psi(i'i'')=\emptyset$ for all $i \in [r]$ and $i'i'' \in \binom{V_i}{2}$.
Similarly, given $\ba \in \Delta^r$ and $\ba' \in \Delta^{r'}$ and a partition $\mathcal{V}=\{V_1,\ldots,V_r\}$ of $[r']$, we will say $\ba' =_{\mathcal{V}} \ba$ if $\sum_{j \in V_i}\alpha_j' = \alpha_i$ for all $i \in [r]$.

Let $(r^*,\phi^*,\ba^*) \in \opt^*(\bm{k})$, let $r \in \mathbb{N}$ and let $[r]$ have partition $\mathcal{V}=\{V_1,\ldots,V_{r^*}\}$.
Let $\phi \in \Phi_0(r+1;\bm{k})$ be such that $\phi|_{\binom{[r]}{2}}=_{\mathcal{V}} \phi^*$ and $\ba =_{\mathcal{V}} \ba^*$.
For $i \in [r^*]$, let
$$
d_i := \sum\{\alpha_{j'} : \phi(\{r+1,j'\}) \neq \phi^*(ij), j' \in V_j, j \in [r^*]\setminus\{i\}\}+\sum\{\alpha_{i'}: |\phi(\{r+1,i'\})| \geq 2, i' \in V_i\}
$$
be the minimum weight of edits of pairs at $r+1$ needed to make $r+1$ a clone of $i$.
If $d_i \leq \delta$, we say that $r+1$ is \emph{$\delta$-close to being a $\phi^*$-clone of $i$}, otherwise $r+1$ is \emph{$\delta$-far from being a $\phi^*$-clone of $i$}.

The next lemma extends the previous one by allowing an arbitrary feasible attachment to a (blow-up of a) basic optimal solution and supposing the new part is far from being a clone.

\begin{lemma}\label{nearidenticalii}
Let $s \in \mathbb{N}$ and let $\bm{k}$ have the extension property.
Then there exists $\eta>0$ such that the following holds.
Let $\delta>0$ and $(r^*,\phi^*,\ba^*) \in \opt^*(\bm{k})$, let $r \in \mathbb{N}$ and let $[r]$ have partition $\mathcal{V}=\{V_1,\ldots,V_{r^*}\}$.
Let $\phi \in \Phi_0(r+1;\bm{k})$ and $\ba \in \Delta^r$ be such that $\phi|_{\binom{[r]}{2}}=_{\mathcal{V}} \phi^*$ and $\ba =_{\mathcal{V}} \bb$ for some $\bb \in \Delta^{r^*}$ with $\|\bb - \ba^*\|_1 \leq \eta\delta$.
Suppose that $r+1$ is $\delta$-far from being a $\phi^*$-clone of any $i \in [r^*]$.
Then $\ext(\phi,\ba) \leq Q(\bm{k})-\eta \delta$.
\end{lemma}

\begin{proof}
We will derive the lemma from the following claim.

\begin{claim}
There exists $\eta>0$ such that when additionally $\bb=\ba^*$, we have $\ext(\phi,\ba) \leq Q(\bm{k})-2\eta\delta \log s$.
\end{claim}

Suppose that the claim holds and we wish to prove the lemma.
Let $\ba' \in \Delta^{r^*}$ be such that $\ba' =_{\mathcal{V}} \ba^*$ and $\|\ba-\ba'\|_1 \leq \eta\delta$.
Such an $\ba'$ exists: for example, for all $j \in [r^*]$ and all $i \in V_j$, take $\alpha'_i := \alpha_i\alpha^*_j/\beta_j$.
Since $\ext(\phi,\ba') \leq Q(\bm{k})-2\eta\delta\log s$, we have
$$
\ext(\phi,\ba) \leq \ext(\phi,\ba') + \log s\cdot \|\ba-\ba'\|_1 \leq Q(\bm{k})- \eta\delta,
$$
as required. So it remains to prove the claim.

\bcpf
Let $\eta'>0$ be the constant obtained from Lemma~\ref{littlec}.
We will show that we can take $\eta := \eta'/(2\log s)$.

It will be convenient to write $0$ instead of $r+1$ for the attachment.
So we require an upper bound for $q_0(\phi,\ba)$.
Let $1/n \ll 1/r,\delta,\eta'$ and for each $j \in [r^*]$,
subdivide the parts in each $V_j$ to get a total of $n$ subparts,
so that as many of these subparts as possible have the same size.
We may assume that in fact every subpart of parts in $V_j$ have the same size $\alpha^*_j/n$,
since the total size of smaller parts is at most $r/n$ which is negligible compared to $\eta'$ and $\delta$.
So, relabelling, we have a partition $\mathcal{U} := \{U_1,\ldots,U_r\}$ of $[r^*n]$ and $\ba_n \in \Delta^{r^*n}$
such that $\ba_{n,k} = \alpha_j^*/n$ for every $k \in \mathcal{U}_j := \{U_{j'} : j' \in V_j\}$, and $|\mathcal{U}_j|=n$.
Write $\mathcal{U}_j := \{x_{j,1},\ldots,x_{j,n}\}$.

For all $\ell \in [n]$, let $T_\ell := \{x_{1,\ell},\ldots,x_{r^*,\ell}\}$ be the $\ell$-th transversal, and let $\phi_\ell := \phi|_{\{0\} \cup T_\ell}$. Recall that $\phi|_{T_\ell} = \phi^*$.
For each $j \in [r^*]$, let $C_j$ be the set of all $\ell \in [n]$ such that $0$ is a clone of $j$ under $\phi_\ell$.
So $C := C_1\cup\ldots\cup C_{r^*}$ is a disjoint union.
By rearranging the transversals, we are going to make all sets $C_j$ empty except at most one.
For this, partition $C$ into pairs $\{\ell_1,\ell_2\}$ and a set $C_{0}$ such that in every pair $\{\ell_1,\ell_2\}$ we have $\ell_1 \in C_{j_1}$ and $\ell_2 \in C_{j_2}$ for distinct $j_1,j_2 \in [r^*]$, and there is at most one $j \in [r^*]$ such that $C_0 \cap C_j \neq \emptyset$.
For all pairs $\{\ell_1,\ell_2\}$, swap the labels of $x_{j_1,\ell_1}$ and $x_{j_1,\ell_2}$.
Update $C$.
Notice that now, $C$ is our previous $C_0$, as neither $\ell_1$ nor $\ell_2$ gives a transversal where $0$ is a clone (since $\phi_{\ell_1}$ has size two on every pair in $\{0\} \cup T_{\ell_1}$, and $\phi_{\ell_2}$ has size at most $1$ on exactly two pairs in $\{0\} \cup T_{\ell_2}$).

Let $\Phi$ be the set of all $\phi_{\ell}$. 
Let $\Phi_{\rm clone} \subseteq \Phi$ be such that $\phi_\ell \in \Phi_{\rm clone}$ if and only if $0$ is a clone under $\phi_\ell$ of some $j \in [r^*]$.
By construction, every such $\ell$ lies in $C$ and there is a unique such $j=j^* \in [r^*]$.

We can make edits of weight at most $1-|C|/n$ to make $0$ a clone of $j^*$ under $\phi$.
Indeed, we can edit each $\phi_\ell$ with $\ell \in [n]\setminus C$, requiring edits to parts of size $\alpha^*_1/n,\ldots,\alpha^*_{r^*}/n$ of total size $1/n$. Thus our hypothesis implies that
$$
1-|C|/n \geq \delta.
$$
 
Lemma~\ref{littlec} implies that $q_0(\psi,\ba^*) \leq Q(\bm{k}) - \eta'$ whenever $\psi \in \Phi\setminus\Phi_{\rm clone}$.
Therefore, using Proposition~\ref{extendbd},
\begin{align*}
q_0(\phi,\ba)  &= \sum_{\ell \in [n]}q_0(\phi_\ell,\ba_n) = \sum_{\ell \in C}q_0(\phi_\ell,\ba^*)/n + \sum_{\ell \in [n]\setminus C}q_0(\phi_\ell,\ba^*)/n\\
&\leq |C|Q(\bm{k})/n  + (n-|C|)((Q(\bm{k})-\eta')/n = Q(\bm{k})-(1-|C|/n)\eta' \leq Q(\bm{k})-\eta'\delta,
\end{align*}
as required.
\ecpf

This completes the proof of the lemma.
\end{proof}

The final lemma in this subsection considers an arbitrary \emph{not necessarily feasible} attachment. Now, either the new part is far from being a clone, or it lies in many forbidden monochromatic cliques.
This is the key tool in the proof of Lemma~\ref{bulk}.

\begin{lemma}\label{nearidenticaliii}
Let $s \in \mathbb{N}$ and let $\bm{k}=(k_1,\ldots,k_s)$ have the extension property, where $k_1 \geq \ldots \geq k_s$.
There exists $\eta>0$ such that the following hold.
Let $\delta>0$ and let $(r^*,\phi^*,\ba^*) \in \opt^*(\bm{k})$, let $r \in \mathbb{N}$ and let $[r]$ have partition $\mathcal{V}=\{V_1,\ldots, V_{r^*}\}$.
Let $\phi : \binom{[r+1]}{2} \to 2^{[s]}$ and $\ba \in \Delta^r$ be such that $\phi|_{\binom{[r]}{2}}=_{\mathcal{V}} \phi^*$
and $\ba =_{\mathcal{V}} \bb$ for some $\bb \in \Delta^{r^*}$ with $\|\bb- \ba^*\|_1 \leq \eta\delta$,
and $\ext(\phi,\ba) \geq Q(\bm{k})-\eta\delta$.
Then one of the following hold.
\begin{itemize}
\item There exists $\ell \in [r]$ such that $r+1$ is $\delta$-close to a $\phi^*$-clone of $\ell$.
\item Let $L$ be the set of sets $\{x_1,\ldots,x_{k_1-1}\} \in \binom{[r]}{k_1-1}$ such that
$\phi^{-1}(c)[\{r+1,x_1,\ldots,x_{k_1-1}\}] \supseteq K_{k_c}$ for some $c \in [s]$. Then $\sum_{(x_1,\ldots,x_{k_1-1}) \in L}\alpha_{x_1}\ldots \alpha_{x_{k_1-1}} \geq \eta$.
\end{itemize}
\end{lemma}

\begin{proof}
Let $\eta'>0$ be the constant obtained from Lemma~\ref{nearidenticalii}.
Suppose for a contradiction that 
for all $n \in \mathbb{N}$, there exist $(r^*_n,\phi^*_n,\ba^*_n) \in \opt^*(\bm{k})$,
$r_n$, $\mathcal{V}_n$, $\phi_n : \binom{[r_n+1]}{2} \to 2^{[s]}$, $\bb_n$ such that $\ba_n =_{\mathcal{V}_n} \bb_n$ and $\|\bb_n-\ba^*_n\|_1 \leq \eta'\delta/2$, $\phi_n |_{\binom{[r_n]}{2}} =_{\mathcal{V}_n} \phi^*_n$, $\ext(\phi_n,\ba_n) \geq Q(\bm{k})-\eta'\delta/2$,
$r_n+1$ is $\delta$-far from being a clone of any $\ell \in [r_n]$ under $\phi_n$, and
defining the set of tuples $L_n$ as in the statement of the lemma, we have $\sum_{(x_1,\ldots,x_{k_1-1}) \in L_n}\alpha_{n,x_1}\ldots \alpha_{n,x_{k_1-1}} < \frac{1}{n}$.
Note that we may assume that $r_n \leq 2^{s r^*_n}$ as $r_n+1$ has at most $2^s$ different attachments to parts in each $V_{n,i}$ in $\mathcal{V}_n$, so if $|V_{n,i}| > 2^{s r^*_n}$ at least two of its parts are clones under $\phi^*_n$ and we can merge them.
As usual, we may assume that $r_n^*=r^*$ and $\phi^*_n=\phi^*$, and thus also $r_n=p$, $\phi_n=\psi$ and $\mathcal{V}_n = \mathcal{V}$.
Choose a convergent subsequence $\ba^*_{i_1},\ba^*_{i_2},\ldots$ of $\ba^*_1,\ba^*_2,\ldots$ with limit $\bb^*$, and
a convergent subsequence $\ba_{i_{j_1}},\ba_{i_{j_2}},\ldots$ of $\ba_{i_1},\ba_{i_2},\ldots$ with limit $\bb$.
The function $\ext(\psi,\cdot)$ is continuous, so $\ext(\psi,\bb) \geq Q(\bm{k})-\eta'\delta/2$.
Also, writing $\bb_{\mathcal{V}}$ to be such that $\bb =_{\mathcal{V}} \bb_{\mathcal{V}}$, we have $\|\bb_{\mathcal{V}}-\bb^*\|_1 \leq \eta'\delta/2$. 
Let $L_c$ be the set of sets $(x_1,\ldots,x_{k_1-1}) \in \binom{[p]}{k_1-1}$ such that $\psi^{-1}(c)[\{p+1,x_1,\ldots,x_{k_1-1}\}] \supseteq K_{k_c}$.
Note that by our assumption above, $L_c$ does not change with $n$.
We have $\sum_{(x_1,\ldots,x_{k_1-1}) \in L_c}\beta_{x_1}\ldots \beta_{x_{k_1-1}} = 0$.
Thus the density of $K_{k_c}$ containing $p+1$ is $0$, and we can remove parts of size $0$ from $[p]$ to obtain a set (without loss of generality $[p']$) such that $\psi^{-1}(c)[\{p+1\} \cup [p']]$ is $K_{k_c}$-free for all $c \in [s]$.
Thus $\psi \in \Phi_0(\{p+1\} \cup [p'];\bm{k})$, i.e. the attachment of $p+1$ under $\psi$ is feasible.
Lemma~\ref{nearidenticalii} implies that $\ext(\psi,\bb) \leq Q(\bm{k})-\eta'\delta$, a contradiction.

Thus there exists $N \in \mathbb{N}$ such that the required sum of tuples is always at least $1/N$. We can now take $\eta$ to be the minimum of $\eta'/2$ and $1/N$.
\end{proof}

\section{Stability of optimal solutions}\label{staboptsols}

The aim of this section is to prove the following lemma, which forms the core of our proof of Theorem~\ref{stabilitysimp}.
Roughly speaking, it says that Problem $Q_0$ is stable, in the sense that both the vertex-weighting and the colour pattern of an almost $Q_0$-optimal solution are close to that of a $Q_0$-optimal solution (which can in turn be described in terms of a $Q_2$-optimal solution by Lemma~\ref{char}).
To prove Theorem~\ref{stabilitysimp}, we will later `transfer' this result to an almost optimal graph~$G$.

\begin{lemma}[Stability of optimal solutions]\label{bulk}
Let $s \in \mathbb{N}$ and let $\bm{k}=(k_1,\ldots,k_s) \in \mathbb{N}^s$ have the extension property, where $k_1 \geq \ldots \geq k_s$. Let $\nu>0$.
Then there exists $\eps>0$ such that
for every $(r,\phi,\ba) \in \feas_0(\bm{k})$ with
$$
q(\phi,\ba) > Q(\bm{k}) - 2\eps,
$$
there is $(r^*,\phi^*,\ba^*) \in \opt^*(\bm{k})$
and a partition $[r]=Y_0 \cup \ldots \cup Y_{r^*}$ such that, defining $\beta_i := \sum_{i' \in Y_i}\alpha_{i'}$ for all $i \in [r^*]$, the following holds.
\begin{enumerate}[label=(\roman*),ref=(\roman*)]
\item\label{bulki} $\|\bb-\ba^*\|_1 < \nu$ (and in particular, $\sum_{i' \in Y_0}\alpha_{i'} < \nu$).
\item\label{bulkii} For all $ij \in \binom{[r^*]}{2}$, $j' \in Y_j$ and $i' \in Y_i$, we have $\phi(i'j') \subseteq \phi^*(ij)$.
\item\label{bulkiii} For all $i \in [r^*]$ and every $i'j' \in \binom{Y_{i}}{2}$, we have $\phi(i'j') \subseteq \{ 1 \}$. 
\end{enumerate}
\end{lemma}

Note that the density of pairs (that is, the sum of the $\alpha_i\alpha_j$) where the inclusion $\phi(i'j') \subseteq \phi^*(ij)$ in~\ref{bulkii} is strict is $O(\nu)$.

\begin{proof}
We will apply a version of symmetrisation to the graph $([r],E)$, where $E$ is the set of pairs $ij$ on which $\phi$ has size at least two.
That is, at each step, we will consider two vertices $j,j'$ with $|\phi(jj')| \leq 1$ and replace one of them with a clone of the other, where the cloned vertex is the one which contributes the most to $q$.
In the first part of the proof we will perform this `forward symmetrisation'. 

Let $\mu$ be the output of Lemma~\ref{solutions} applied to $\bm{k}$, so $\alpha_j^* > \mu$ for all $(r^*,\phi^*,\ba^*) \in \opt^*(\bm{k})$ and $j \in [r^*]$.
Choose additional constants $\eps,\gamma,\eta,\delta$ such that
$
\eps \ll \gamma \ll \eta \ll \delta \ll \mu,\nu
$
where 
$\sqrt{\delta}$ is at most the output of Lemma~\ref{littlec}, 
$\eta$ is at most the output of Lemma~\ref{nearidenticaliii},
$\eps^{1/4}$ is at most the output of Lemma~\ref{Sfnear} with $\gamma$ playing the role of $\delta$.

Now let $(r,\phi,\ba) \in \feas_0(\bm{k})$ satisfy $q(\phi,\ba)>Q(\bm{k})-2\eps$.
Add all $i$ with $\alpha_i=0$ to $Y_0$.
We can take $\tau \ll \eps, \min_{i \in [r]}\alpha_i \leq 1/r$ so that, subdividing each part, we can remove subparts of total size at most $r\tau$ so that every other subpart has size exactly $\tau$.
Since we can put the removed parts into $Y_0$, we may assume without loss of generality that all $\alpha_i$ are equal to each other (and $\phi$ takes the value $\emptyset$ between parts obtained from subdividing a single original part).
So we may assume that
$\alpha_i = 1/r \ll \eps$ for all $i \in [r]$.
Altogether we have
$$
\alpha_1=\ldots = \alpha_r = 1/r \ll \eps \ll \gamma \ll \eta \ll \delta \ll \mu,\nu.
$$

\medskip
\noindent
\underline{\textit{The forward symmetrisation procedure.}}

\medskip
\noindent
Let $\mathcal{X}_0 := \{\{1\},\ldots,\{r\}\}$ and $\phi_0 := \phi$.

\begin{claim}\label{cl-1}
There is $f \in \mathbb{N}$ such that, after relabelling $[r]$, for all $i=0,\ldots,f$, there is a partition $\mathcal{X}_i$ of $[r]$ and colour pattern $\phi_i \in \Phi_0(r;\bm{k})$ such that 
the following hold.
\begin{enumerate}[label=(\roman*),ref=Claim~\ref{cl-1}(\roman*)]
\item\label{F2} There is a single $x_{i} \in [r]$ such that $\mathcal{X}_{i}$ consists of the same elements as $\mathcal{X}_{i-1}$, except that $x_i$ has moved from one part to another.
\item\label{F3} $\phi_i =_{\mathcal{X}_{i}} \psi_i$ where $\psi_i = \phi|_{\binom{[r_i]}{2}}$,
and $\psi_f \in \Phi_2(r_f;\bm{k})$ for some $2 \leq r_f \leq \ldots \leq r_0 = r$.
\item\label{F4} $q(\phi_{i},\ba) - q(\phi_{i-1},\ba) \geq 0$ (where $\phi_{-1}=\phi_0$). 
\end{enumerate}
\end{claim}

\bcpf
Let $W_0 := [r]$, 
$\phi_0 = \psi_0 := \phi$ and let $V_{0,x} := X_{0,x} := \{ x \}$ for all $x \in [r]$.
Let $\mathcal{V}_0 := \{\{x\}: x \in [r]\}$ and $\bb_0 := \ba$.
Initialise $i_1:=0$.

Inductively for $j \geq 0$, perform \emph{forward superstep $j+1$} by defining $W_{j+1}$, $\psi_{j+1}$, $\mathcal{V}_{j+1} := \{V_{j+1,x} : x \in W_{j+1}\}$, $\bb_{j+1}$ as follows.
Choose a pair $p_jt_j \in \binom{W_j}{2}$ with
$|\psi_j(p_jt_j)| \leq 1$, 
labelled so that $t_j$ has larger attachment under $\phi_j$; that is
\begin{equation}\label{goodchoice}
\sum_{\substack{yt_j \in \binom{W_j}{2}:\\ \psi_j(yt_j) \neq \emptyset}}\beta_{j,y}\log|\psi_j(yt_j)| \geq \sum_{\substack{yp_j \in \binom{W_j}{2}:\\ \psi_j(yp_j) \neq \emptyset}}\beta_{j,y}\log|\psi_j(yp_j)|.
\end{equation}
If there is no such pair, terminate the iteration.
Otherwise, let $W_{j+1} := W_j \setminus \{ p_j \}$.
Obtain $\mathcal{V}_{j+1}$ from $\mathcal{V}_j$ by replacing $V_{j,p_j},V_{j,t_j}$ with their union,
so $V_{j+1,t_j} := V_{j,t_j} \cup V_{j,p_j}$ and $V_{j+1,x} := V_{j,x}$ for all $x \in W_{j+1}\setminus\{t_j\}$.
For all $x \in W_{j+1}$, let $\beta_{j+1,x} := |V_{j+1,x}|/r$.
Let $\psi_{j+1} := \psi_j|_{\binom{W_{j+1}}{2}}$.
Note that
\begin{align}
\nonumber &\phantom{=} \sum_{\substack{xy \in \binom{W_{j+1}}{2}:\\ \psi_{j+1}(xy) \neq \emptyset}}\beta_{j+1,x}\beta_{j+1,y}\log|\psi_{j+1}(xy)| - \sum_{\substack{xy \in \binom{W_j}{2}:\\ \psi_j(xy) \neq \emptyset}}\beta_{j,x}\beta_{j,y}\log|\psi_{j}(xy)|\\
\label{diffeq} &= \beta_{j,p_j} \left( \sum_{\substack{yt_j \in \binom{W_j}{2}:\\ \psi_j(yt_j) \neq \emptyset}}\beta_{j,y}\log|\psi_j(yt_j)| - \sum_{\substack{yp_j \in \binom{W_j}{2}:\\ \psi_j(yp_j) \neq \emptyset}}\beta_{j,y}\log|\psi_j(yp_j)| \right) \stackrel{(\ref{goodchoice})}{\geq} 0.
\end{align}

Now we will symmetrise each part in $V_{j,p_j}$ one by one, defining $\phi_{i+1} : \binom{[r]}{2} \to 2^{[s]}$ for $i_j \leq i < i_{j+1}$ where $i_{j+1} := i_j+|V_{j,p_j}|$.
Let $y^* \in V_{j,t_j}$ be arbitrary.
Let $s_{j+1} := |V_{j,p_j}|$ and write $V_{j,p_j} := \{v_1,\ldots,v_{s_{j+1}}\}$.
We will perform $s_{j+1}$ \emph{forward steps}, as follows.
Inductively for $i \geq i_j$,
obtain $\mathcal{X}_{i+1}$ from $\mathcal{X}_{i}$ by moving $v_i$ from $X_{i,p_j}$ to $X_{i,t_j}$.
That is, for $i_j \leq i < i_{j+1}-1$, let $X_{i+1,t_j} := X_{i,t_j} \cup \{v_i\}$, $X_{i+1,p_j} := X_{i,p_j}\setminus \{v_i\}$ and $X_{i+1,x} := X_{i,x}$ for all $x \in W_{j+1}\setminus \{t_j\}$;
if $i=i_{j+1}-1$ we do the same but instead discard the (empty) $p_j$-th part, so $|\mathcal{X}_{i+1}|=|W_j|$ for $i_j \leq i < i_{j+1}$, while $|\mathcal{X}_{i_{j+1}}|=|W_{j+1}|$.
Let $v_i$ become a strong clone of $y^*$ in $\phi_{i+1}$; that is, for distinct $x,y \in [r]$, define
$$
\phi_{i+1}(xy) := \begin{cases} 
\phi_{i}(xy) &\mbox{if } x,y \neq v_i \\
\phi_{i}(y^*z) &\mbox{if } \{x,y\}=\{v_i,z\} \text{ and } z \neq y^* \\
\emptyset &\mbox{if } \{x,y\}=\{v_i,y^*\}.
\end{cases}
$$
After defining $\phi_{i+1}$ and $\mathcal{X}_{i+1}$ for all $i_j \leq i < i_{j+1}$, we proceed with superstep $j+2$.

The iteration will run until some forward step $i=f$ (for \emph{final}) when $|\phi_f(xy)| \geq 2$ for all $x,y$ in different parts of $\mathcal{X}_f$.
The process terminates in a finite number of steps since $|W_j|$ is strictly decreasing (so there are finitely many supersteps $j$), and there are finitely many steps $s_j$ at each superstep $j$.

Let $r_i := |\mathcal{X}_i|$.
By relabelling the elements of $[r]$, for all supersteps $j$, we can assume that
$W_j$ is always an initial segment of $[r]$ so we have $\psi_j = \phi|_{\binom{[|W_j|]}{2}}$.
Let $\ba_i := (\alpha_{i,1},\ldots,\alpha_{i,r_i}) := (|X_{i,1}|/r,\ldots,|X_{i,r_i}|/r) \in \Delta^{r_i}$.
We have shown that for each $i \in [f]$, we can obtain a function $\phi_i$, and sets $\mathcal{X}_i$ such that~\ref{F2} and~\ref{F3} hold.

We still need to prove~\ref{F4}. It is true by definition for $i=0$.
Equation~(\ref{goodchoice}) implies that $q_{t_j}(\psi_{j},\bb_{j}) - q_{p_j}(\psi_j,\bb_j) \geq 0$.
At step $i+1$ during the $j+1$ superstep, we change the attachment of a single vertex $v_i$, and we have $|\phi_i(v_iy)| \leq 1$ for all $y \in V_{j,p_j} \cup V_{j,t_j}$. Thus the only change to $q_{v_i}$ is for pairs $v_ix$ with $x \in V_{j,j'}$ for $j' \in W_{j+1}\setminus \{t_j\}$. Thus $q_{v_i}(\phi_{i+1},\ba) - q_{v_i}(\phi_i,\ba)$ 
is the difference of the left- and right-hand sides of~(\ref{goodchoice}).
The required statement follows since $q(\phi_{i+1},\ba)-q(\phi_i,\ba) = \alpha_{v_i}(q_{v_i}(\phi_{i+1},\ba) - q_{v_i}(\phi_i,\ba))$.
\ecpf

Since $\phi_f =_{\mathcal{X}_f} \psi_f \in \Phi_2(r_f;\bm{k})$ by definition, we also have that
\begin{equation}\label{feas}
(r_f,\psi_f,\ba_f) \in \textsc{feas}_2(\bm{k})
\quad\text{and}\quad
q(\phi_f,\ba) = q(\psi_f,\ba_f).
\end{equation}
Moreover, \ref{F4}~implies that
\begin{eqnarray}\label{interval}
Q(\bm{k}) \geq q(\psi_f,\ba_f) \geq q(\psi_{f-1},\ba_{f-1}) \geq \ldots \geq q(\psi_0,\ba_0)
\geq  Q(\bm{k}) - 2\eps.
\end{eqnarray}
Note that Lemma~\ref{Sfnear} implies there is some vertex weighting $\bb$ close to $\ba_f$ such that $(r_f,\psi_f,\bb)$ is optimal (but it could have zero parts).
So `forward symmetrisation' has allowed us to pass from our original feasible solution $(r,\phi,\ba)$ to a new feasible solution $(r_f,\psi_f,\ba_f)$, which is very close to a $Q_2$-optimal solution (both in terms of vertex weighting and colour pattern).
But our eventual aim is to show that $(r,\phi,\ba)$ itself is close to this optimal solution.
So we need to show that few `significant' changes were made during the forward symmetrisation procedure.
To this end, our next step will be to follow the procedure backwards, using the partitions $\mathcal{X}_i$ of $[r]$ at each step, to form a new partition $\mathcal{U}_i$ of $[r]$, which records how the solution at each step differs from $(r_f,\psi_f,\ba_f)$.

It will be convenient to define some normalised versions $\hat{q},\hat{q}_x$ of $q,q_x$ (for $x \in [r]$).
Here we recall that $\alpha_1=\ldots=\alpha_r=1/r$ which makes the normalisation simpler.
Given $(r,\phi,\bm{\ba}) \in \Phi_0(r;\bm{k})$ and $P \subseteq [r]$, write
\begin{align*}
&q(P,\phi) := 2\sum_{\substack{xy \in \binom{P}{2}:\\ \phi(xy) \neq \emptyset}}\alpha_x\alpha_y|\log\phi(xy)|
\quad\text{and}\quad
q_x(P,\phi) := \sum_{\substack{y \in P\setminus\{x\}:\\ \phi(xy) \neq \emptyset}}\alpha_y|\log\phi(xy)|,\quad\text{and}\\
&\hat{q}(P,\phi) := \left(\frac{r}{|P|}\right)^2 \cdot q(P,\phi)
\quad\text{and}\quad
\hat{q}_x(P,\phi) := \frac{r}{|P|} \cdot q_x(P,\phi),
\end{align*}
so that 
\begin{equation}\label{eq-norm}
\sum_{x \in P}\hat{q}_x(P,\phi) = |P|\hat{q}(P,\phi)\quad\text{and}\quad
\hat{q}(P,\phi) \leq Q(\bm{k})
\end{equation}
(if this inequality were not true, then setting $\beta_x:=1/|P|$ for $x \in P$ and $\beta_x:=0$ otherwise gives $q(\phi,\bb)>Q(\bm{k})$).
 
\medskip
\noindent
\underline{\textit{The backwards symmetrisation procedure.}}

The forwards symmetrisation procedure ended with $(r_f,\psi_f,\ba_f) \in \feas_2(\bm{k})$ which is very close to optimal. 
We now want to go backwards through each forwards step $i$ in turn, each time defining a partition of $[r]$ into $r_f$ sets $U^1,\ldots, U^{r_f}$ corresponding to the vertices of $\psi_f$ as well as a small exceptional set $U^0$.
The desired conclusion is that, at the end of this process, the final sets $U^1,\ldots,U^{r_f}$, that is, those corresponding to the original $\phi$ we started with, have sizes roughly $\alpha_1'r,\ldots, \alpha'_{r_f} r$ for some vertex weighting $\ba'$ where $(r_f,\psi_f,\ba') \in \opt_2(\bm{k})$ (so $\ba'$ could differ significantly from $\ba_f$ and could have zero parts, but is nevertheless optimal). Thus the exceptional set together with the `extra' parts outside of the support of $\ba'$ are small. This will mean that between parts, $\phi_i$ resembles $\psi_f$ throughout the process, but the sizes of the parts $U^1,\ldots,U^{r_f}$ could change during the process.
Thus $\phi$ resembles $\psi_f$ on the support of $\ba'$ in the required sense.

At each forwards step $i$, we modified the solution $\phi_{i-1}$ to obtain a new solution $\phi_i$ by changing the attachment at a single vertex $x_i$, so that $q$ did not decrease.
Now, in the corresponding backwards step, initially no vertex is exceptional. Then, we reconsider the attachment at $x_i$: if it was small in $\phi_i$ we remove it into the exceptional set $U^0$. If any other vertex $y$ also has small attachment in $\phi_i$ we also remove it to $U^0$.
If $x_i$ was not removed, we assign it to the part $U^j$ where $x_i$ looks most like a $\psi_f$-clone of $j$ in $\phi_i$, and similarly assign vertices which are no longer exceptional.

The extension property guarantees that any vertex which was not moved into the exceptional set, and therefore has large attachment, looks similar to a $\psi_f$-clone.
There cannot be too many exceptional vertices since they all have small attachment, whereas $q$ is large.

We now formally describe the $i$-th backwards step.
Define $\mathcal{U}_{f} := \{ U^0_{f},\ldots,U^{r_f}_{f} \}$ by setting $U^j_{f} := X_{f,j}$ for all $j \in [r_f]$ and $U^0_{f} := \emptyset$. 
For each $i=f-1,\ldots,0$, define $U_i$ and $\mathcal{U}_{i} := \{ U^0_{i},\ldots,U^{r_f}_{i} \}$ inductively as follows. Initially, $U_i=[r]$ and $U^0_i=\emptyset$.
If 
$\hat{q}_{x_i}(U_{i},\phi_{i}) < Q(\bm{k})-\sqrt{\eps}$, move $x_i$ from $U_i$ into $U^0_{i}$.
Next, if there is $y \in U_{i}$ such that $\hat{q}_y(U_{i},\phi_{i})<Q(\bm{k})-\sqrt{\eps}$, move $y$ into $U^0_{i}$ (we also include the special vertex $x_i$ here, if at some point its attachment becomes too small).
Update $U_{i}$ and repeat until there are no such vertices left in $U_{i}$.

Next, for each $j \in [r_f]$, let $U_i^j$ be the restriction of $U_{i+1}^j$ to $U_i\setminus \{x_i\}$. For each $z \in B_i := (U_{i+1}^0 \cup \{x_i\}) \cap U_i$, add $z$ to the part $U_i^j$ such that $z$ looks most like a $\psi_f$-clone of $j$ under $\phi_i|_{U_i}$; that is, choose the index $j\in[r_f]$ such that
$$
\sum_{j' \in [r_f]\setminus \{j\}}\left(|\{y \in U_i^{j'}: \phi_i(yz) \neq \psi_f(j'j)\}|\right)+|\{y \in U_i^j: |\phi_i(yz)| \geq 2\}|
$$
is minimal (breaking ties arbitrarily).
This completes \emph{backwards step} $i$; now move on to backwards step $i-1$.

We show that the exceptional set $U_i^0$ is always small.

\begin{claim}\label{cl-Ui0}
For all $i=f,\ldots,0$, we have $|U^0_{i}| \leq 2\sqrt{\eps}r$.
\end{claim}

\bcpf
Let $y_1,\ldots,y_\ell$ be the vertices which are moved into $U_{i}^0$ at step $i$, in this order.
So $|U_{i}^0| = \ell$.

Given distinct $x,y \in [r]$, write $d_i(xy) := \log|\phi_i(xy)|$ if $\phi_i(xy)\neq\emptyset$ and $d_i(xy):=0$ otherwise.
For $1 \leq k \leq \ell$, the vertex $y_k$ is moved to $U^0_{i}$
due to $\hat{q}_{y_k}(U_{i,k},\phi_{i}) < Q(\bm{k})-\sqrt{\eps}$, where $U_{i,k} := [r]\setminus \{y_1,\ldots,y_{k-1}\}$. 
Note that $\hat{q}_{y_k}(U_{i,k},\phi_i)|U_{i,k}|=\sum_{x \in U_{i,k+1}}d_i(xy_k)$.
We have
\begin{eqnarray*}
Q(\bm{k})\frac{(r-\ell)^2}{2} &\stackrel{(\ref{eq-norm})}{\geq}& \hat{q}(U_i,\phi_i)\frac{|U_i|^2}{2} = \sum_{xy \in \binom{U_i}{2}}d_i(xy) =  \sum_{xy \in \binom{[r]}{2}}d_i(xy) - \sum_{k \in [\ell]}\hat{q}_{y_k}(U_{i,k},\phi_i)|U_{i,k}|\\
&\geq& \sum_{xy \in \binom{[r]}{2}}d_i(xy) - \sum_{k \in [\ell]}(Q(\bm{k})-\sqrt{\eps})(r-k+1)\\
&=& q(\phi_i,\ba)\frac{r^2}{2} - (Q(\bm{k})-\sqrt{\eps})\left(r\ell - \binom{\ell}{2}\right)\\
&\stackrel{(\ref{interval})}{\geq}& (Q(\bm{k})-2\eps)\frac{r^2}{2}- (Q(\bm{k})-\sqrt{\eps})\left((r-\ell)\ell + \binom{\ell+1}{2}\right). 
\end{eqnarray*}
Rearranging, we have $\sqrt{\eps}((r-\ell)\ell+\binom{\ell+1}{2}) \leq \eps r^2 +Q(\bm{k})\ell/2 < 3\eps r^2/2$.
But if $2\sqrt{\eps}r \leq \ell \leq r/2$, we have $\sqrt{\eps}(r-\ell)\ell \geq \sqrt{\eps}(1-2\sqrt{\eps})\cdot 2\sqrt{\eps}r^2 > 3\eps r^2/2$.
Thus, at the moment when $2\sqrt{\eps}r$ vertices are added to $U^0_i$, we obtain a contradiction.
\ecpf

Every $x \in U_i$ satisfies $\hat{q}_x(U_i,\phi_i) \geq Q(\bm{k})-\sqrt{\eps}$.
Let $n := |U_i|$ and let $G_i$ be the complete graph with vertex set $U_i$ whose edges are coloured red (for missing), blue (for extra) or green (for perfect) as follows.
For each $x \in U_i$, let $j_x \in [r_f]$ be such that $x \in U^{j_x}_i$.
For each $xy \in E(G_i)$,
\begin{itemize}
\item $xy$ is \emph{red} if $j_x \neq j_y$ and $\phi_i(xy) \subsetneq \psi_f(j_xj_y)$, so there are missing colours.
\item $xy$ is \emph{blue} if either
$j_x \neq j_y$ and $\phi_i(xy) \setminus \psi_f(j_xj_y) \neq \emptyset$, or 
$j_x=j_y$ and $\phi_i(xy) \not\subseteq \{1\}$,
so there are extra colours.
\item $xy$ is \emph{green} otherwise. 
\end{itemize}
Recall that we defined
\begin{equation}\label{eq-B}
B_i = (U^0_{i+1} \cup \{x_i\}) \cap U_i,\quad\text{and that}\quad|B_i| \leq 3\sqrt{\eps}r
\end{equation}
by Claim~\ref{cl-Ui0}.
The colouring of $G_i-B_i$ depends only on the previous partition $\mathcal{U}_{i+1}$
and the colour pattern $\phi_{i+1}$ since every vertex in $U_i \setminus B_i$
lies in $U^j_{i+1} \cap U^j_i$ for some $j \in [r_f]$, and the colour patterns $\phi_i$ and $\phi_{i+1}$ only differ at $x_i$.
Also,
$$
n = r - |U^0_i| \geq (1-2\sqrt{\eps})r.
$$
Write $\bb_i := (|U_i^1|,\ldots,|U_i^{r_f}|)/|U_i| \in \Delta^{r_f}$.
\begin{claim}\label{cl-blue}
For all $i=f,\ldots,0$, $G_i$ has no blue edges. 
\end{claim}
\bcpf
We prove this by backwards induction for $i=f,\ldots,0$.
The claim is true for $i=f$, as every edge in $G_f$ is green.
Suppose it is true for all backwards steps $f,\ldots,i+1$.
The induction hypothesis and the fact that $\phi_i$ differs from $\phi_{i+1}$ only at $x_i$ implies that only vertices in $B_i$ can be incident with blue edges in $G_i$.

First we show that $G_i$ contains few red edges and $\bb_i$ is close to an optimal vertex weighting. Indeed,
\begin{align*}
Q(\bm{k})-\sqrt{\eps} &\leq \sum_{x \in U_i}\hat{q}_x(U_i,\phi_i)/|U_i| = \hat{q}(U_i,\phi_i)\\
&\leq\left(\frac{r}{|U_i|}\right)^2\left(\frac{2}{r^2}\sum_{\substack{xy \in E(G_i) \\ \psi_f(j_xj_y) \neq \emptyset}}\log|\psi_f(j_xj_y)| - \frac{2}{r^2}\sum_{xy \text{ red}}\log\left(\frac{|\psi_f(j_xj_y)|}{|\psi_f(j_xj_y)|-1}\right)+\frac{2\log s}{r}|B_i|\right)\\
&\leq q(\psi_f,\bb_i)-2\log\left(\frac{s}{s-1}\right)\frac{e_{\rm red}(G_i)}{n^2} + \frac{2r\log s|B_i|}{n^2}\\
&\leq Q(\bm{k})-2\log\left(\frac{s}{s-1}\right)\frac{e_{\rm red}(G_i)}{n^2} + 7\sqrt{\eps}\log s.
\end{align*}
Thus
\begin{equation}\label{eq-redG}
e_{\rm red}(G_i) \leq \eps^{1/3}n^2
\end{equation}
and additionally, from the penultimate inequality, $q(\psi_f,\bb_i) \geq Q(\bm{k}) - \eps^{1/4}$. 
Lemma~\ref{Sfnear} applied with parameters $s,\bm{k},\gamma$ implies that there exists
$\ba'_i \in \Delta^{r_f}$ such that $(r_f,\psi_f,\ba'_i) \in \opt_2(\bm{k})$ and
\begin{equation}\label{eq-Uij}
\|\bb_i-\ba'_i\|_1 < \gamma \ll \eta \delta.
\end{equation}
Without loss of generality, suppose the nonzero entries of $\ba'_i$ form an initial segment of length $\tilde{r}_i$, and let $\tilde{\ba}_i$ be this initial segment. 
Let $\tilde{\phi}_i := \psi_f|_{\binom{[\tilde{r}_i]}{2}}$. 
Then $(\tilde{r}_i,\tilde{\phi}_i,\tilde{\ba}_i) \in \opt^*(\bm{k})$.
We claim that
\begin{equation}\label{eq-alphaprime}
q_{j}(\psi_f,\ba'_i) \leq Q(\bm{k})\quad\text{for all}\quad j \in [r_f].
\end{equation}
Indeed, if $j \in [\tilde{r}_i]$, then Proposition~\ref{lagrange} implies that $q_{j}(\psi_f,\ba_i')=q_j(\tilde{\phi}_i,\tilde{\ba}_i) = Q(\bm{k})$.
If $\tilde{r}_i < j \leq r_f$, then Proposition~\ref{extendbd}
implies the bound $q_{j}(\psi_f,\ba'_i)=\ext(\psi_f|_{\binom{[\tilde{r}_i] \cup \{j\}}{2}},\tilde{\ba}_i) \leq Q(\bm{k})$.

Next we claim that every
$z \in B_i$ is $\delta$-close under $\phi_i|_{\binom{U_i}{2}}$ to being a $\psi_f$-clone of some $j \in [\tilde{r}_i] \subseteq [r_f]$, which will follow from an application of Lemma~\ref{nearidenticaliii}.
Suppose not. To apply the lemma, let 
$U'_i := \bigcup_{j \in [\tilde{r}_i]}U^j_i$ and
$U^z_i := \{z\} \cup (U_i' \setminus B_i)$ and $t:=|U_i' \setminus B_i|$.
Let $\tilde{\bb}_i := (|U_i^1\setminus B_i|,\ldots,|U_i^{\tilde{r}_i}\setminus B_i|)/|U_i' \setminus B_i|$, so $\|\tilde{\bb}_i-\tilde{\ba}_i\|_1 \leq 2\gamma$ by~(\ref{eq-Uij}).
Let $\phi': \binom{U^z_i}{2} \to 2^{[s]}$ be obtained from $\phi_i$ as follows. Let $\phi'(zy):=\phi_i(zy)$ for all $y \in U_i^z\setminus\{z\}$, and let $\phi'$ agree with $\phi_f$ elsewhere, that is, 
$\phi'(xy):=\psi_f(j_xj_y)$ whenever $j_x \neq j_y$, and $\phi'(xy):=\emptyset$ if $j_x=j_y$.
Let $\ba_t$ be the length-$t$ vector which is identically $1/t$.
We have
$$
\ext(\phi',\ba_t)=q_z(\phi_i|_{\binom{U^z_i}{2}},(\tfrac{1}{t},\ldots,\tfrac{1}{t},0)) \geq \frac{1}{t}\left(|U_i|\hat{q}_z(U_i,\phi_i)-(\log s)(|B_i|+\gamma r)\right)>Q(\bm{k})-r(\log s)\gamma.
$$
We can apply Lemma~\ref{nearidenticaliii} with $(\tilde{r}_i,\tilde{\phi}_i,\tilde{\ba}_i),(U^j_i \setminus B_i: j \in [\tilde{r}_i]),\ba_t,\tilde{\bb_i},z$ playing the roles of $(r^*,\phi^*,\ba^*),\mathcal{V},\ba,\bb,r+1$ to see that, writing $L$ for the set of sets $\{y_1,\ldots,y_{k_1-1}\} \in \binom{U_i'\setminus B_i}{k_1-1}$ (i.e.~$(k_1-1)$-subsets of vertices of $G_i[U_i']-B_i$) such that $(\phi')^{-1}(c)[\{z,y_1,\ldots,y_{k_1-1}\}] \supseteq K_{k_c}$ for some $c \in [s]$, we have $|L| \geq \eta n^{k_1-1}$.

For every $\{y_1,\ldots,y_{k_1-1}\} \in L$, there are $\ell,\ell' \in [k_1-1]$ such that $y_\ell y_{\ell'}$ is red (recalling that these edges are either red or green),
otherwise $\phi_i$ is identical to $\phi'$ on all pairs of these vertices, and thus $\phi_i^{-1}(c)$ contains a copy of $K_{k_c}$ for some $c$.
Each pair appears in at most $n^{k_1-3}$ sets in $L$, so the number of red edges in $G_i$ is at least
$
\eta n^{k_1-1}/n^{k_1-3} = \eta n^2 > 2\eps^{1/3} n^2
$,
a contradiction to~(\ref{eq-redG}).
Thus $z$ is $\delta$-close to being a $\psi_f$-clone of some $j \in [\tilde{r}_i]$ under $\phi_i|_{\binom{U_i^z}{2}}$.

During backwards symmetrisation we added $z$ to the part $U_i^{j_{z}}$ such that $z$ was closest to a $\psi_f$-clone of $j_{z}$ under $\phi_i$, 
so 
\begin{align*}
d_{\rm red}(z)+d_{\rm blue}(z)
&\leq \sum_{j' \in [r_f]\setminus\{j_z\}}|\{y \in U^{j'}_i: \phi_i(zy) \neq \psi_f(j_zj_y)\}|+|\{y \in U^{j_z}_i: |\phi_i(zy)| \geq 2\}|\\
&\leq \sum_{j' \in [r_f]\setminus\{j\}}|\{y \in U^{j'}_i: \phi_i(zy) \neq \psi_f(jj_y)\}|+|\{y \in U^{j}_i: |\phi_i(zy)| \geq 2\}|\\
&\leq \delta t + |B_i| \leq 2\delta n.
\end{align*}
Therefore the green degree of $z$ in $G_i$ is
\begin{equation}\label{eq-greenxi}
d_{\rm green}(z) \geq (1-2\delta)n \quad \text{for all }z \in B_i.
\end{equation}
Thus
\begin{eqnarray}
\nonumber Q(\bm{k})-\sqrt{\eps} &\leq& \hat{q}_z(U_i,\phi_i)\leq \frac{r}{|U_i|}\left( \frac{1}{r}\sum_{\substack{y \in U_i \\ \psi_f(j_zj_y) \neq \emptyset}}\log|\psi_f(j_zj_y)| + \frac{\log s}{r}(n-1-d_{\rm green}(z))\right)\\
\nonumber &\leq& q_{j_z}(\psi_f,\bb_i) + 3\delta\log s \leq q_{j_z}(\psi_f,\ba'_i) + 2(\log s)\|\ba'_i-\bb_i\|_1 + 3\delta\log s\\
\nonumber &\stackrel{(\ref{eq-Uij})}{\leq}& q_{j_z}(\psi_f,\ba'_i) + 4\delta\log s
\end{eqnarray}
and therefore
\begin{equation}\label{eq-Bz}
q_{j_z}(\psi_f,\ba_i') \geq Q(\bm{k})-\sqrt{\delta}\quad\text{for all }z \in B_i.
\end{equation}
(By the optimality of $(r_f,\psi_f,\ba_i')$ this is automatically true for nonzero parts.)
Next, we show the number of red edges incident to a vertex $x$ in $U_i\setminus B_i$ is
\begin{equation}\label{eq-greenx}
d_{\rm red}(x) \leq \sqrt{\gamma}n \quad\text{and hence}\quad
d_{\rm green}(x) \geq (1-2\sqrt{\gamma})n \quad\text{for all }x \in U_i\setminus B_i.
\end{equation}
Indeed, the second part follows from the first since every edge incident to $x$ in $G_i-B_i$ is either green or red.
To prove the first part, let $x \in U_i\setminus B_i$.
Since $x$ can have blue neighbours only in $B_i$, we have
\begin{align}
\nonumber Q(\bm{k})-\sqrt{\eps} &\leq \hat{q}_x(U_i,\phi_i)\\
\nonumber &\leq \frac{r}{|U_i|}\left( \frac{1}{r}\sum_{\substack{y \in U_i \setminus B_i \\ \psi_f(j_xj_y) \neq \emptyset}}\log|\psi_f(j_xj_y)| - \frac{1}{r}\sum_{y \in N_{\rm red}(x)}\log\left(\frac{|\psi_f(j_xj_y)|}{|\psi_f(j_xj_y)|-1}\right) + \frac{\log s}{r}|B_i|\right)\\
\nonumber &\leq q_{j_x}(\psi_f,\bb_i) - \log\left(\frac{s}{s-1}\right)\frac{d_{\rm red}(x)}{n} + 4\log s\sqrt{\eps}\\
\nonumber &\leq q_{j_x}(\psi_f,\ba'_i)  - \frac{d_{\rm red}(x)}{sn}+3\gamma\log s\\
\label{eq-attach} &\leq Q(\bm{k})  - \frac{d_{\rm red}(x)}{sn}+3\gamma\log s,
\end{align}
where the final inequality follows from~(\ref{eq-alphaprime}).
Therefore, $d_{\rm red}(x) \leq \sqrt{\gamma}n$, as required,
and the penultimate inequality implies that $q_{j_x}(\psi_f,\ba'_i) \geq Q(\bm{k})-\gamma^{1/3}$.

Combined with~(\ref{eq-Bz}), we have shown that $q_{j_y}(\psi_f,\ba'_i) \geq Q(\bm{k})-\sqrt{\delta}$ for all $y \in U_i$.
We will now show that this means that $\bb_i$ has the same support as $\ba'_i$;
that is, either $\tilde{r}_i=r_f$, or for all $\tilde{r}_i<j \leq r_f$, we have $\beta_{i,j}=0$.
Suppose not; then without loss of generality there is some $x \in U^{\tilde{r}_i+1}_i$ (so $j_x=\tilde{r}_i+1$). 
We have $\ext(\psi_f|_{\binom{[\tilde{r}_i+1]}{2}},\tilde{\ba}_i) = q_{\tilde{r}_i+1}(\psi_f,(\alpha'_{i,1},\ldots,\alpha'_{i,\tilde{r}_i},0)) \geq Q(\bm{k})-\sqrt{\delta}$.
Thus Lemma~\ref{littlec} implies that $\tilde{r}_i+1$ is a $\psi_f$-clone of some $j^* \in [\tilde{r}_i]$ under $\psi_f$,
which is a contradiction since $|\psi_f(\{\tilde{r}_i+1,j^*\})| \geq 2$ for all $j^* \in [\tilde{r}_i]$.
Thus $U^j_i=\emptyset$ for all $\tilde{r}_i < j \leq r_f$, so~(\ref{eq-Uij}) implies that
\begin{equation}\label{eq-Uij2}
\beta_{i,j} \geq \tilde{\alpha}_{i,j}-\gamma \geq \mu-\gamma \geq \mu/2\quad\text{for all }j \in [\tilde{r}_i],
\quad\text{and}\quad U_i = \bigcup_{j \in [\tilde{r}_i]}U^j_i.
\end{equation}
We can now complete the claim, comparing $\phi_i$ and the partition $\bigcup_{j \in [\tilde{r}_i]}U^j_i$ of $U_i$ to $(\tilde{r}_i,\tilde{\phi}_i,\tilde{\ba}_i) \in \opt^*(\bm{k})$.
Suppose for a contradiction that there is a blue edge $zy$, so $z \in B_i$ and $y \in U_i$
(where $y$ could also be in $B_i$).
Let $j_1 := j_{z}$ and $j_2 := j_y$, so $\{j_1,j_2\} \subseteq [\tilde{r}_i]$ by~(\ref{eq-Uij2}).
By definition, either 
\begin{itemize}
\item[(i)] $j_1 \neq j_2$ and there is $c \in \phi_i(zy)\setminus \tilde{\phi}_i(j_1j_2)$, or
\item[(ii)] $j_1=j_2$ and there is some $1 \neq c \in \phi_i(zy)$.
\end{itemize}
We claim that, in both cases, there exist $j_3,\ldots,j_{k_c} \in [\tilde{r}_i]\setminus\{j_1,j_2\}$ such that $c \in \tilde{\phi}_i(j_\ell j_{\ell'})$ for all pairs among $\{j_1,\ldots,j_{k_c}\}$ except $j_1j_2$ if they are distinct. 
Indeed, suppose~(i) holds. 
By Lemma~\ref{nocap}\ref{nocapi}, the graph $([\tilde{r}_i],(\tilde{\phi}_i)^{-1}(c))$ is maximally $K_{k_c}$-free. Since it is not complete, we are done.
Suppose~(ii) holds.
By Lemma~\ref{nocap}\ref{nocapii}, $(\bm{1}+\bm{e}_{j_1})(\tilde{\phi}_i)^{-1}(c)$ contains a copy of $K_{k_c}$ as $c \neq 1$.

We say that a subset $\{y_3,\ldots,y_{k_c}\}$ with $y_\ell \in U_i^{j_\ell}$ for $3 \leq \ell \leq k_c$ is \emph{bad} if $zy_\ell$, $yy_\ell$ and $y_{\ell}y_{\ell'}$ for every $\ell\ell'\in\binom{[k_c]}{2}$ are green. 
Since $c \in \phi_i(zy)$, there can be no bad subsets in $G_i$ since then $c \in \phi_i(xy)$ for every pair $xy$ among vertices in the subset, contradicting $K_{k_c} \notin \tilde{\phi}_i^{-1}(c)([\tilde{r}_i])$.
On the other hand, at least, say, $\prod_{3 \leq \ell \leq k_c}(|U_i^{j_\ell}|/2)$ subsets are bad.
Indeed,~(\ref{eq-Uij2}) implies that $|U_i^j| =\beta_{i,j}n \geq \mu n/2$ for all $j \in [\tilde{r}_i]$.
Also,~(\ref{eq-greenx}) and~(\ref{eq-greenxi}) imply that every vertex has at most $2\delta n$ nongreen neighbours.
Thus, choosing $y_3,\ldots,y_{k_c}$ sequentially, among the vertices in $U_i^{j_\ell}$, there are at most $2(\ell-1)\delta n < \mu n/4 < |U_i^{j_\ell}|/2$ vertices forbidden for $y_\ell$ due to not being a green neighbour of every $y_1,\ldots,y_{\ell-1}$, as required.
This contradiction implies that $z$ is not incident to any blue edges, and thus $G_i$ contains no blue edges.
This finishes the proof of Claim~\ref{cl-blue}.
\ecpf

\medskip
\noindent
As before, let us assume that nonzero entries of $\ba'_0$ are indexed by $[\tilde{r}_0]$.
So Lemma~\ref{bulk} holds when we set $Y_0 := U_0^0 \cup \bigcup_{\tilde{r}_0+1 \leq j \leq r_f}U_0^j$ and $Y_j := U_0^j$ for all $j \in [\tilde{r}_0]$ and $(\tilde{r}_0,\tilde{\phi}_0,\tilde{\ba}_0)$ plays the role of $(r^*,\phi^*,\ba^*)$.
This completes the proof of Lemma~\ref{bulk}.
\end{proof}

\section{Stability of asymptotically extremal graphs}\label{staboptgraphs}

\subsection{Tools for large graphs}\label{regtools}

One of our main tools is Szemer\'edi's regularity lemma, which allows us to discretise a large edge-coloured graph and thus approximate it by a feasible solution to Problem $Q_0$. 
We will need the following definitions relating to regularity.

\begin{definition}[Edge density, regularity of pairs and partitions]
\rm
Given a graph $G$ and disjoint non-empty sets $A,B \subseteq V(G)$, we define the \emph{edge density} between $A$ and $B$ to be
$$
d_G(A,B) := \frac{e_G(A,B)}{|A||B|}.
$$
Given $\eps,d > 0$, the pair $(A,B)$ is called
\begin{itemize}
\item \emph{$\eps$-regular} if for every $X \subseteq A$ and $Y \subseteq B$ with $|X| \geq \eps|A|$ and $|Y| \geq \eps|B|$, we have that $|d(X,Y) - d(A,B)| \leq \eps$.
\item \emph{$(\eps, d)$-regular} if $(A,B)$ is 
$\eps$-regular and $d_G(A,B) = d \pm \eps$.
\item \emph{$(\eps,\geq\! d)$-regular} if it is $\eps$-regular and has density at least $d-\eps$.
\end{itemize}
An \emph{equitable partition} of a set $V$ is a partition of $V$ into parts $V_1,\ldots,V_m$ such that $|\,|V_i| - |V_j|\,| \leq 1$ for all $i,j \in [m]$.
An equitable partition of $V(G)$ into parts $V_1,\ldots,V_m$ is called \emph{$\eps$-regular} if $|V_i| \leq \eps|V(G)|$ for every $i \in [m]$, and all but at most $\eps\binom{m}{2}$ of the pairs $(V_i,V_j)$ are $\eps$-regular.
\end{definition}

We use the following multicolour version of Szemer\'edi's regularity lemma~\cite{reg}.
This version can be deduced from the original; see for example Theorems~1.8 and~1.18 in Koml\'os and Simonovits~\cite{komsim}.

\begin{theorem}[Multicolour regularity lemma]\label{multicol}
For every $\eps > 0$ and $s \in \mathbb{N}$ , there exists $M \in \mathbb{N}$ such that for any graph $G$ on $n \geq M$ vertices and any edge $s$-colouring $\chi : E(G) \rightarrow [s]$, there is an equitable partition $V(G) = V_1 \cup \ldots \cup V_m$ with $1/\eps \leq m \leq M$, which is $\eps$-regular simultaneously with respect to all graphs $(V(G),\chi^{-1}(i))$, with $i \in [s]$.\qed
\end{theorem}

Our first tool states that a subgraph of a regular pair is still regular, provided both parts are not too small.

\begin{proposition}\cite[Proposition~9]{psy}\label{badrefine}
Let $\eps,\delta$ be such that $0 < 2\delta \leq \eps < 1$.
Suppose that $(X,Y)$ is a $\delta$-regular pair, and let $X' \subseteq X$ and $Y' \subseteq Y$.
If
$$
\min \left\{ \frac{|X'|}{|X|} , \frac{|Y'|}{|Y|} \right\} \geq \frac{\delta}{\eps},
$$
then $(X',Y')$ is $\eps$-regular.\qed
\end{proposition}

The next proposition states that, given a set of edge-disjoint subgraphs $G_1,\ldots,G_s$ of a bipartite graph, if at least one of the graphs $G_i$ is not regular of density $s^{-1}$, then there is a $G_j$ whose density on a pair of large sets is reduced.

\begin{proposition}\label{badred}
Let $A,B$ be disjoint sets of vertices, $s \in \mathbb{N}$ and let $\eps > 0$ be a constant with $1/|A|,1/|B| \ll \eps \ll 1/s$.
Let $G_1,\ldots,G_s$ be pairwise edge-disjoint subgraphs of $K[A,B]$.
Suppose that not all of $G_1,\ldots,G_s$ are $(\eps,s^{-1})$-regular graphs.
Then there exists $c \in [s]$ and $X \subseteq A$, $Y \subseteq B$ with $|X| = \lceil \eps|A|\rceil$ and $|Y| = \lceil \eps |B| \rceil$ such that 
$$
d_{G_c}(X,Y) \leq \frac{1}{s}\left( 1 - \frac{\eps}{2}\right).
$$
\end{proposition}

\begin{proof}
Given $c \in [s]$, $X \subseteq A, Y \subseteq B$, let
$$
\text{diff}_c(X,Y) := s^{-1}|X||Y|-e_{G_c}(X,Y).
$$
If $G_c$ is not $(\eps,s^{-1})$-regular, then either
\begin{enumerate}
\item[(i)] $|\text{diff}_c(A,B)| > \frac{\eps}{2}|A||B|$; or
\item[(ii)]
there is some $X \subseteq A$ and $Y \subseteq B$ with $|X| \geq \eps|A|$ and $|Y| \geq \eps|B|$ such that
$$
\bigg|\frac{e_{G_c}(X,Y)}{|X||Y|}-\frac{e_{G_c}(A,B)}{|A||B|}\bigg| > \eps;
$$
\end{enumerate}
or both. 
(The immediate implication from the definition of $(\eps,s^{-1})$-regular would have (i) replaced by $|\text{diff}_c(A,B)| > \eps|A||B|$, which is stronger than the statement of (i).)
To prove the proposition, it is enough to exhibit $c^* \in [s]$, $X' \subseteq A$ and $Y' \subseteq B$ with $|X'| \geq \eps|A|$ and $|Y'| \geq \eps|B|$ so that
\begin{equation}\label{diffaim}
\text{diff}_{c^*}(X',Y') \geq \frac{\eps}{2s}|X'||Y'|.
\end{equation}
Indeed, if we can find such $c^*,X',Y'$, then, setting $k_1 := \lceil \eps|A|\rceil$ and $k_2 := \lceil \eps|B|\rceil$, we have that
\begin{align*}
\sum_{\stackrel{X \subseteq X'}{|X|=k_1}}\sum_{\stackrel{Y \subseteq Y'}{|Y|=k_2}} \text{diff}_{c^*}(X,Y) &= \binom{|X'|}{k_1}\binom{|Y'|}{k_2}s^{-1}k_1k_2 - \binom{|X'|-1}{k_1-1}\binom{|Y'|-1}{k_2-1}e_{G_{c^*}}(X',Y')\\
&= \binom{|X'|-1}{k_1-1}\binom{|Y'|-1}{k_2-1} \text{diff}_{c^*}(X',Y');
\end{align*}
so, by averaging, there is some $X \subseteq X'$ and $Y \subseteq Y'$ with $|X|=k_1$, $|Y|=k_2$ such that
\begin{equation}
\nonumber\text{diff}_{c^*}(X,Y) \geq \text{diff}_{c^*}(X',Y') \cdot \frac{k_1k_2}{|X'||Y'|} = \frac{\eps}{2s}|X||Y|,
\end{equation}
as required.
So we will now concentrate on finding $c^*,X',Y'$ so that (\ref{diffaim}) holds.
Suppose first that (i) holds for some $c \in [s]$.
If $\text{diff}_c(A,B) > 0$, then we are done by setting $c^* := c$, $X' := A$ and $Y' := B$.
So we may assume that $\text{diff}_c(A,B) < 0$. 
Observe that
\begin{eqnarray*}
\sum_{i \in [s]}\text{diff}_i(A,B) \geq 0.
\end{eqnarray*}
So $\sum_{i \in [s]\setminus \{ c \}}\text{diff}_i(A,B) \geq \eps|A||B|/2$.
By averaging, there is some $c' \in [s]$ such that
$$
\text{diff}_{c'}(A,B) \geq \frac{\eps|A||B|}{2(s-1)} \geq \frac{\eps}{2s}|A||B|.
$$
So we are done by setting $c^* := c'$, $X' := A$ and $Y' := B$.

Suppose instead that (ii) holds for some $c \in [s]$.
So there are $X \subseteq A$, $Y \subseteq B$ with $|X| \geq \eps|A|$, $|Y| \geq \eps|B|$ such that
$$
\eps < \bigg| \frac{\text{diff}_c(X,Y)}{|X||Y|} - \frac{\text{diff}(A,B)}{|A||B|} \bigg| < \frac{| \text{diff}_c(X,Y)|}{|X||Y|} + \frac{\eps}{2}.
$$
Therefore $|\text{diff}_c(X,Y)| > \eps|X||Y|/2$.
Again, we may assume that $\text{diff}_c(X,Y) < 0$, or we are done.
Then an almost identical argument to the one above yields $c' \in [s]$ such that $\text{diff}_{c'}(X,Y) \geq \eps|A||B|/2s$.
This completes the proof.
\end{proof}

The next proposition states that regular pairs are robust under small perturbations; the version stated here is a slight variation of Proposition~8 in~\cite{bst}.

\begin{proposition}\label{adjust}
Let $(A,B)$ be an $(\eps, d)$-regular pair and let $(A',B')$ be a pair such that $|A' \bigtriangleup A| \leq \alpha|A'|$ and $|B' \bigtriangleup B| \leq \alpha|B'|$ for some $0 \leq \alpha \leq 1$.
Then $(A',B')$ is an $(\eps+7\sqrt{\alpha}, d)$-regular pair.\qed
\end{proposition}

We will also frequently use the following standard embedding lemma (see, for example, Theorem~2.1 in~\cite{komsim}).

\begin{lemma}[Embedding lemma]\label{embed}
For every $\eta >0$ and integer $k \geq 2$ there exist $\eps>0$ and $m_0 \in \mathbb{N}$ such that the following holds.
Suppose that $G$ is a graph with a partition $V(G) = V_1 \cup \ldots \cup V_k$ such that $|V_i| \geq m_0$ for all $i \in [k]$, and every pair $(V_i,V_j)$ for $1 \leq i < j \leq k$ is $(\eps,\geq\! \eta)$-regular.
Then $G$ contains $K_k$.\qed
\end{lemma}

\subsubsection{Binomial tails}

In order to prove Part~\ref{stabilitysimpii} of Theorem~\ref{stabilitysimp}, we will need Corollary~\ref{binbound} below, which is a simple consequence of the Chernoff inequality.
The combinatorial interpretation of this fact is that almost every partition of $[n]$ into $k$ parts is such that every part has size roughly $n/k$.
Write $X \sim \text{Bin}(n,p)$ if a random variable $X$ is binomially distributed with parameters $n \in \mathbb{N}$, $p \in (0,1)$.

\begin{proposition}\cite[Theorem 2.1]{randomgraphs} \label{hoeffding} 
Suppose $X \sim \text{Bin}(n,p)$ where $0<p<1$. 
Let $k \leq np$.
Then
$$
\mathbb{P}(X \leq k) \leq \exp \left( \frac{-(np-k)^2}{2np}\right).
$$
\end{proposition}

\begin{corollary}\label{binbound}
Let $n,k \in \mathbb{N}$  and $\delta \in \mathbb{R}$, where $0 < 1/n \ll \delta \ll 1/k$.
Then
$$
\sum_{i=0}^{\lfloor (k^{-1}-\delta) n \rfloor}\binom{n}{i} (k-1)^{n-i} \leq e^{-\delta^2 kn/3} \cdot k^n.
$$
\end{corollary}

\begin{proof}
Let $X \sim \text{Bin}(n,k^{-1})$ be a binomial random variable.
Then Proposition~\ref{hoeffding} implies that
$$
\mathbb{P}(X \leq \lfloor(k^{-1}-\delta)n\rfloor) \leq \exp \left( \frac{-(n/k-\lfloor (k^{-1}-\delta)n \rfloor)^2}{2n/k} \right) \leq e^{-\delta^2kn/3}.
$$
But
$$
\mathbb{P}(X \leq \lfloor (k^{-1}-\delta)n\rfloor) = \sum_{i=0}^{\lfloor(k^{-1}-\delta)n\rfloor}\binom{n}{i}\left(\frac{1}{k}\right)^i \left(1-\frac{1}{k}\right)^{n-i} = k^{-n} \cdot \sum_{i=0}^{(k^{-1}-\delta)n\rfloor}\binom{n}{i}(k-1)^{n-i},
$$
as required.
\end{proof}

We will also need the following simple bound, which we state without proof.
Let $n \in \mathbb{N}$ and $\eps,\delta > 0$ such that $0 < 1/n \ll \eps \ll \delta < 1$.
Then
\begin{equation}\label{bin}
\binom{n}{\leq \eps n} \leq 2^{\delta n}, 
\end{equation}
where for integers $m \geq t$ we write $\binom{m}{\leq t}:=\sum_{0 \leq i \leq t}\binom{m}{i}$.

\subsection{Preparation for the proof of Theorem~\ref{stabilitysimp}}

We define a hierarchy of constants and assume that these relations hold throughout the remainder of this section.
Let $s \in \mathbb{N}$ and $\bm{k} \in \mathbb{N}^s$, and let $\delta > 0$.
In what follows, whenever we assume that a constant is sufficiently small, it is because a larger constant gives a weaker conclusion.
Let $\mu > 0$ be such that $\alpha^*_i \geq \mu$ for all $(r^*,\phi^*,\ba^*) \in \opt^*(\bm{k})$ and $i \in [r^*]$ (which exists by Lemma~\ref{solutions}).
We may assume that $\delta \ll \mu, 1/s, 1/R(\bm{k})$.
Choose $\gamma_3, \ldots, \gamma_6 \in \mathbb{R}$ such that $0 < \gamma_3 \ll \ldots \ll \gamma_6 \ll \delta$.
In particular, we may assume that $\gamma_4 \leq \gamma_5/2M_{\gamma_5}$, where $M_{\gamma_5}$ is the integer output of Theorem~\ref{multicol} applied with parameter $\gamma_5$; and $\gamma_5$ is at most the output of Lemma~\ref{embed} applied with parameter $\gamma_6$.
Let $0 < \nu \ll \gamma_3$.
Let $\eps > 0$ be such that $8\eps$ is the output of Lemma~\ref{Sfnear} applied with $2\nu$;
we may assume that $\eps \ll \nu$.
Choose $\gamma_2 \in \mathbb{R}$ such that $0 < \gamma_2 \ll \eps$.
Let $\gamma_1 > 0$ be the minimum constant obtained when Lemma~\ref{embed} is applied with $\gamma_2$ playing the role of $\eta$, and with $k_1,\ldots,k_c$ playing the role of $k$.
We may assume that $0 < \gamma_1 \ll \gamma_2$.
Apply Theorem~\ref{multicol} with parameter $\gamma_1$ to obtain $M \in \mathbb{N}$ such that the conclusions of the theorem hold.
We may assume that $1/M \ll \gamma_1$.
Now let $n_0 \in \mathbb{N}$ be such that $1/n_0 \ll 1/M$.
We have the hierarchy
\begin{equation}\label{stabhier}
0 < \frac{1}{n_0} \ll \frac{1}{M} \ll \gamma_1 \ll \gamma_2 \ll \eps \ll \nu \ll \gamma_3 \ll \gamma_4 \ll \gamma_5 \ll \gamma_6 \ll \delta \ll \mu\ll\frac{1}{R(\bm{k})}.
\end{equation}

We need the following somewhat technical definition of `popular vectors' from~\cite{PY}, which allows us to choose colourings $\chi$ of $G$ whose coloured regularity partition is a witness of many other valid colourings of $G$.

\begin{definition}[Popular vectors]\label{popular}
\rm
Let $G$ be a graph on $n \geq n_0$ vertices, and $\chi : E(G) \rightarrow [s]$ be an $s$-edge colouring of $G$ which is $\bm{k}$-valid.
Apply Theorem~\ref{multicol} to the pair $(G,\chi)$ with parameter $\gamma_1$ to obtain
an equitable partition $V(G) = U_1 \cup \ldots \cup U_r$ with $1/\gamma_1 \leq r \leq M$, which is $\gamma_1$-regular simultaneously with respect to all graphs $(V(G),\chi^{-1}(c))$, with $c \in [s]$.
Let
$$
\phi(ij) := \{ c \in [s] : \chi^{-1}(c)[U_i,U_j] \text{ is } (\gamma_1, \geq\! \gamma_2)\text{-regular} \}.
$$
Let $\mathcal{U} := \{ U_i : i \in [r] \}$.
We define the function $\RL$ by setting
$$
\RL(\chi) := (r,\phi,\mathcal{U})
$$
(where we arbitrarily fix a single output if there is more than one choice of $(r,\phi,\mathcal{U})$).
We say that $(r,\phi,\mathcal{U})$ is \emph{popular} if 
$$
|\textstyle{\RL^{-1}}((r,\phi,\mathcal{U}))| \geq F(G;\bm{k}) \cdot 2^{-3\eps n^2},
$$
and \emph{unpopular} otherwise.
Let $\mathrm{Pop}(G)$ be the set of popular $(r,\phi,\mathcal{U})$ and let $\mathrm{Col}(G)$ be the set of $\bm{k}$-valid colourings $\chi$ of $G$ such that $\RL(\chi) \in \mathrm{Pop}(G)$.
\end{definition}

As the following proposition shows, almost every colouring $\chi$ maps to a popular vector.

\begin{proposition}\label{mainlypop}
For all graphs $G$ on $n \geq n_0$ vertices,
$$
|\mathrm{Col}(G)| \geq (1-2^{-2\eps n^2}) \cdot F(G;\bm{k}).
$$
\end{proposition}

\begin{proof}
Let $M$ be the integer output of Theorem~\ref{multicol} applied with parameter $\gamma_1$.
Let $n \geq M$ and let $G$ be a graph on $n$ vertices.
The function $\RL$ is well-defined.
Then the number of outputs $(r,\phi,\mathcal{U})$ is at most
$$
M \cdot \left( 2^{\binom{M}{2}} \right)^s \cdot n^{M} = 2^{O(\log n)}.
$$
Now,
\begin{align*}
\sum_{(r,\phi,\mathcal{U}) \in \mathrm{Pop}(G)}|\textstyle{\RL^{-1}}((r,\phi,\mathcal{U}))| &= F(G;\bm{k}) - \sum_{(r,\phi,\mathcal{U}) \notin \mathrm{Pop}(G)}|\textstyle{\RL^{-1}}((r,\phi,\mathcal{U}))|\\
&\geq \left(1 - 2^{O(\log n)} \cdot 2^{-3\eps n^2} \right) F(G;\bm{k}) \geq (1-2^{-2\eps n^2})F(G;\bm{k}),
\end{align*}
as required.
\end{proof}

\subsection{The proof of Theorem~\ref{stabilitysimp}}

Using Lemma~\ref{bulk}, we can now prove Theorem~\ref{stabilitysimp}.
Although this lemma is really the heart of the proof, there are still many steps required to `transfer' its conclusion to the graph setting.
For this reason, we split the proof into a series of claims, and continue to use the constants defined in~(\ref{stabhier}).

\begin{proof}[Proof of Theorem~\ref{stabilitysimp}]
Suppose that $G$ is a graph on $n \geq n_0$ vertices, and 
\begin{equation}\label{assumption}
\frac{\log F(G;\bm{k})}{\binom{n}{2}} \geq Q(\bm{k})-\eps.
\end{equation}
We will show that the conclusion of Theorem~\ref{stabilitysimp} holds with parameter $\delta$.
If we decrease $\delta$, then the conclusion of Theorem~\ref{stabilitysimp} becomes only stronger, so we can assume that $\delta$ satisfies~(\ref{stabhier}).
Let $(r,\phi,\mathcal{U}) \in \mathrm{Pop}(G)$.
That is,
\begin{equation}\label{popularuse}
|\textstyle{\RL^{-1}}((r,\phi,\mathcal{U}))| \geq 2^{-3\eps n^2} \cdot F(G;\bm{k}).
\end{equation}
We will (for now) suppress the dependence of what follows on $(r,\phi,\mathcal{U})$.
Thus there is an equitable partition $V(G) = U_1 \cup \ldots \cup U_r$, where $1/\gamma_1 \leq r \leq M$, which is, for all $\chi \in \RL^{-1}((r,\phi,\mathcal{U}))$, $\gamma_1$-regular simultaneously with respect to all graphs $(V(G),\chi^{-1}(c))$, with $c \in [s]$.
Furthermore, for each $ij \in \binom{[r]}{2}$ and $c \in [s]$, we have that $c \in \phi(ij)$ if and only if $\chi^{-1}(c)[U_i,U_j]$ is $(\gamma_1,\geq\! \gamma_2)$-regular.
Lemma~\ref{embed} and our choice of parameters in~(\ref{stabhier}) imply that $\phi^{-1}(c)$ is $K_{k_c}$-free for all $c \in [s]$.

The next claim shows  that $G$ gives rise to a feasible solution $(r,\phi,\ba)$ of Problem~$Q_0$ which is almost optimal.
Moreover, $\ba$ is a good approximation of the structure of $G$, and because $(r,\phi,\mathcal{U})$ is popular, $\phi$ is a good approximation of many valid colourings of $G$.

\begin{claim2}\label{largecon}
Let $\ba:= ( |U_1|/n, \ldots, |U_r|/n )$.
Then
$(r,\phi,\ba) \in \feas_0(\bm{k})$ and $$
q(\phi,\ba) \geq Q(\bm{k})-8\eps + 2\sum_{\stackrel{ij \in \binom{[r]}{2}}{|\phi(ij)| \geq 2}} \left( \frac{e(\overline{G}[U_i,U_j])}{n^2} \right).
$$
Moreover,
for every $\chi \in \RL^{-1}((r,\phi,\mathcal{U}))$, we have
\begin{equation}\label{almostfull}
\sum_{\stackrel{ij \in \binom{[r]}{2}}{|\phi(ij)|\geq2}}\frac{e(\overline{G}[U_i,U_j])}{n^2} \leq 4\eps
\end{equation}
and there are at most $s\gamma_2 n^2$ edges $xy \in E(G)$ where $x \in U_i$ and $y \in U_j$, such that either $i=j$, or $i \neq j$ and $\chi(xy) \notin \phi(ij)$.
\end{claim2}

\bcpf
Consider the following procedure for producing colourings of $G$ whose image under $\RL$ is $(r,\phi,\mathcal{U})$.

\medskip
\noindent
\underline{\textit{Standard colouring procedure.}}

\begin{itemize}
\item[1.] Colour `atypical' edges as follows:
\begin{itemize}
\item[(i)] Assign arbitrary colours to all edges of $G$ that lie inside some part $U_i$.
\item[(ii)] Select at most $s\gamma_1\binom{r}{2}$ elements of $\binom{[r]}{2}$ and, for each selected pair $ij$, assign colours to $G[U_i,U_j]$ arbitrarily.
\item[(iii)] For every colour $c \in [s]$ and every $ij \in \binom{[r]}{2}$, colour an arbitrary subset of edges of $G[U_i,U_j]$ of size at most $\gamma_2|U_i||U_j|$ by colour $c$.
\end{itemize}
\item[2.] Colour most edges according to $\phi$: for every edge $ij \in \binom{[s]}{2}$ and $x \in U_i$, $y \in U_j$ where $xy \in E(G)$ and $xy$ is not yet coloured, pick an arbitrary colour from the set $\phi(ij)$. If $\phi(ij) = \emptyset$, colour $xy$ with colour $1$.
\end{itemize}
This procedure will generate every $\chi \in \RL^{-1}((r,\phi,\mathcal{U}))$ (as well as some further colourings, which may not even be $\bm{k}$-valid).
Indeed, this follows from Theorem~\ref{multicol}, and the statement that $\phi^{-1}(c)$ is $K_{k_c}$-free for all $c \in [s]$.

Let $S_1$ be the number of choices in Step~1.
We will call those edges which are not coloured according to $\phi$ (i.e.~not coloured in Step~2) the \emph{atypical} edges.
The number of these is at most
$$
r\left\lceil \frac{n}{r}\right\rceil^2 + s\gamma_1\binom{r}{2}\left\lceil \frac{n}{r} \right\rceil^2  + s\cdot \binom{r}{2} \cdot \gamma_2\left\lceil \frac{n}{r}\right\rceil^2 < s\gamma_2 n^2,
$$
proving the second part of the claim.
This also implies that
$$
S_1 \leq \binom{\binom{n}{2}}{\leq s\gamma_2 n^2}s^{s\gamma_2 n^2} \stackrel{(\ref{bin})}{<} 2^{\eps n^2/3}.
$$
Let $S_2$ be the number of choices in Step~2 given a fixed choice at Step~1.
Since $(r,\phi,\mathcal{U})$ is popular, we have that
\begin{equation}\label{logs2}
\log S_2 \stackrel{(\ref{popularuse})}{\geq} \log \left( 2^{-3\eps n^2} \cdot F(G;\bm{k}) \right) - \log S_1 \stackrel{(\ref{assumption})}{\geq} \binom{n}{2}\left( Q(\bm{k})-\eps \right) - \frac{\eps n^2}{3} - 3\eps n^2.
\end{equation}
We would now like to bound $S_2$ from above.
For each $ij \in \binom{[r]}{2}$, define $\delta_{ij}$ by setting
$$
\delta_{ij}n^2 := e(\overline{G}[U_i,U_j]) = |U_i||U_j| - e(G[U_i,U_j]).
$$
Now,
\begin{equation}\label{S2}
S_2 \leq \prod_{ij \in \binom{[r]}{2}} \left( \max \{ 1, |\phi(ij)| \} \right)^{|U_i||U_j|-\delta_{ij}n^2}.
\end{equation}
So
\begin{align}
\nonumber \log S_2 &\leq \sum_{\stackrel{ij \in \binom{[r]}{2}}{\phi(ij) \neq \emptyset}}\left(|U_i||U_j|-\delta_{ij}n^2 \right)\log|\phi(ij)| = n^2 \sum_{\stackrel{ij \in \binom{[r]}{2}}{\phi(ij) \neq \emptyset}}\left(\alpha_i\alpha_j-\delta_{ij} \right)\log|\phi(ij)|\\
\label{logS2upper} &= \frac{n^2}{2} q(\phi,\ba) - n^2\sum_{\stackrel{ij \in \binom{[r]}{2}}{|\phi(ij)| \geq 2}}\delta_{ij}\log|\phi(ij)|. 
\end{align}
Combining this with (\ref{logs2}), we have that
$$
q(\phi,\ba) \geq Q(\bm{k}) - 8\eps  + 2\sum_{\stackrel{ij \in \binom{[r]}{2}}{|\phi(ij)|\geq 2}} \delta_{ij},
$$
proving the first part of the claim. Every edge that the second part of the claim counts is atypical and by construction there are at most $s\gamma_2 n^2$ of these.
The final part of the claim follows from~(\ref{logs2}),~(\ref{S2}) and~(\ref{logS2upper}).
\ecpf

Apply Lemma~\ref{bulk} with parameter $2\nu$ to $(r,\phi,\ba)$ to obtain $(r^*,\phi^*,\ba^*) \in \opt^*(\bm{k})$ such that the following hold: there is a partition
$
[r] = V_0 \cup \ldots \cup V_{r^*}
$
where for all $i \in [r^*]$ we have
\begin{equation}\label{u}
\|\bm{x}-\ba^*\|_1 < 2\nu \text{ where }x_i := \sum_{j \in V_i}\alpha_j;
\end{equation}
for all $ij \in \binom{[r^*]}{2}$, $i' \in V_i$ and $j' \in V_j$ we have that $\phi(i'j')\subseteq\phi^*(ij)$;
and for all $i \in [r^*]$ and every $i'j' \in \binom{V_i}{2}$, we have $\phi(i'j') \subseteq \{ 1 \}$.

We would like to transfer this partition to $G$ itself. So for all $0 \leq i \leq r^*$, let
\begin{equation}\label{X}
X_i := \bigcup_{j \in V_i}U_j,\ \text{ so } \ V(G) = X_0 \cup \ldots \cup X_{r^*}.
\end{equation}
Then it is easy to see that
\begin{equation}\label{boldx}
\bm{x} = \left( \frac{|X_1|}{n},\ldots,\frac{|X_{r^*}|}{n} \right).
\end{equation}
Now~(\ref{u}) and~(\ref{stabhier}) imply that, for all $i \in [r^*]$,
\begin{equation}\label{mu}
|X_i| \geq (\alpha^*_i-2\nu)n \geq  \mu n/2\quad\text{and}
\end{equation}
\begin{equation}\label{X0}
2\nu > \|\ba^*-\bm{x}\|_1 \geq \bigg| \|\ba^*\|_1 - \|\bm{x}\|_1 \bigg| = \left| 1 - \left(1-\frac{|X_0|}{n}\right) \right| = \frac{|X_0|}{n}.
\end{equation}
Note that $\bm{x}$, $(r^*,\phi^*)$, $[r]=V_0 \cup \ldots \cup V_{r^*}$, $X_0, \ldots, X_{r^*}$ are fixed for every $\chi \in \RL^{-1}((r,\phi,\mathcal{U}))$.
Claim~\ref{largecon} implies that
\begin{equation}\label{eq-innergood}
\sum_{i \in [r^*]}(e_G(X_i) - |\chi^{-1}(1)[X_i]|) < 3\nu n^2\quad\text{for all }\chi \in \textstyle{\RL^{-1}}((r,\phi,\mathcal{U})).
\end{equation}
Say that $\chi \in \RL^{-1}((r,\phi,\mathcal{U}))$ is \emph{good} if
\begin{itemize}
\item $\chi^{-1}(c)[X_i,X_j]$ is $(\gamma_3,|\phi^*(ij)|^{-1})$-regular for all $ij \in \binom{[r^*]}{2}$ and $c \in \phi^*(ij)$.
\end{itemize}
Say that $\chi$ is \emph{bad} otherwise.
Let $\mathcal{G}=\mathcal{G}(r,\phi,\mathcal{U})$ be the set of good colourings $\chi \in \RL^{-1}((r,\phi,\mathcal{U}))$.
We will show that almost every $\chi$ is good.
The idea here is that, in every bad colouring $\chi$, there is a pair $(X_{i^*},X_{j^*})$ in which some colour graph of $\chi$ is not regular of the correct density.
Lemma~\ref{badred} implies that there must be some colour $c$ and large sets $X \subseteq X_{i^*}$ and $Y \subseteq X_{j^*}$ between which $\chi^{-1}(c)$ has density which is significantly smaller than expected.
So there are significantly fewer choices for colouring the edges between this pair, a loss which is quantified by Corollary~\ref{binbound}.

Consider the following procedure for generating a set of colourings of $G$ which (as we will show) includes every bad colouring.

\medskip
\noindent
\emph{\underline{Bad colouring procedure.}}

\begin{enumerate}
\item Choose at most $3\nu n^2$ edges of $G$ and colour them arbitrarily. For each $i \in [r^*]$, colour every remaining edge in $G[X_i]$ with colour $1$.
\item Pick $i^*j^* \in \binom{[r^*]}{2}$; $c^* \in \phi^*(i^*j^*)$ and subsets $X \subseteq X_{i^*}$ and $Y \subseteq X_{j^*}$ of size $\lceil \gamma_3|X_{i^*}| \rceil, \lceil \gamma_3|X_{j^*}|\rceil$ respectively.
\item Choose at most $(|\phi^*(i^*j^*)|^{-1}-\gamma_3/2s)|X||Y|$ edges in $G[X,Y]$ and colour them with colour~$c^*$. Arbitrarily colour the remaining edges in $G[X,Y]$ with colours from $\phi^*(i^*j^*) \setminus \{ c^* \}$.
\item Arbitrarily colour the remaining edges in $G[X_{i^*},X_{j^*}]$ with colours from $\phi^*(i^*j^*)$.
\item For all $ij \in \binom{[r^*]}{2} \setminus \{ i^*j^* \}$, arbitrarily colour all remaining edges in $G[X_i,X_j]$ using colours from $\phi^*(ij)$.
\end{enumerate}
Let $S_{p_1\ldots p_2}$ be the number of choices in Steps~$p_1$--$p_2$, having fixed choices in previous steps, where $[p_1,p_2] \subseteq [5]$.

\begin{claim2}\label{badprimer}
The number of bad $\chi \in \RL^{-1}((r,\phi,\mathcal{U}))$ is at most $S_{1\ldots 5}$.
\end{claim2}

\bcpf
If suffices to show that for any bad $\chi \in \RL^{-1}((r,\phi,\mathcal{U}))$, there is a set of choices in the bad colouring procedure which generates it.
So fix such a $\chi$.
Say that an edge $xy$ is \emph{contrary} if one of the following holds:
\begin{itemize}
\item[(a)] at least one of $x,y$ is in $X_0$;
\item[(b)] $\chi(xy) \notin \phi^*(ij)$, where $i \neq j$ and $x \in X_i$ and $y \in X_j$;
\item[(c)] $x,y \in X_i$ and $xy$ is not coloured with colour $1$.
\end{itemize}
By~(\ref{X0}), the number of edges of type (a) is at most $|X_0|n \leq 2\nu n^2$.
By Claim~\ref{largecon}, there at most $s\gamma_2 n^2$ edges $xy$ with $x \in U_i$, $y \in U_j$, such that either $i=j$, or $i \neq j$ and $\chi(xy) \notin \phi(ij)$.
Combining this with Lemma~\ref{bulk}\ref{bulkiii}, we see that the number of edges of types (b) and (c) is at most $s\gamma_2 n^2$.
Therefore there are at most $3\nu n^2$ contrary edges in $G$.
We colour these edges in Step~1.

Since $\chi$ is bad,  
there is some $i^*j^* \in \binom{[r^*]}{2}$ and $c \in [s]$ such that $\chi^{-1}(c)[X_{i^*},X_{j^*}]$ is not $(\gamma_3,|\phi^*(i^*j^*)|^{-1})$-regular.
Proposition~\ref{badred}
applied with $|\phi^*(i^*j^*)|,X_{i^*},X_{j^*}, \chi^{-1}(c)[X_{i^*},X_{j^*}],\gamma_3$ playing the roles of $r,A,B,G_c,\eps$
 implies that there exists $c^* \in [s]$ such that there are $X \subseteq X_{i^*}$, $Y \subseteq X_{j^*}$ with $|X| = \lceil \gamma_3 |X_{i^*}| \rceil$, $|Y| = \lceil \gamma_3|X_{j^*}|\rceil$ where 
$$
d(\chi^{-1}(c)(X,Y)) \leq |\phi^*(i^*j^*)|^{-1}\left(1-\frac{\gamma_3}{2}\right) \leq |\phi^*(i^*j^*)|^{-1} - \frac{\gamma_3}{2s}.
$$
So in Step~2 we can take $i^*j^* \in \binom{[r^*]}{2}$, $c^* \in [s]$ and $X \subseteq X_{i^*}$, $Y \subseteq X_{j^*}$, and choose a suitable colouring in Steps~3 and~4, which will generate $\chi[X_{i^*},X_{j^*}]$.
The only uncoloured edges are non-contrary edges in $(X_i,X_j)$ for $ij \in \binom{[r^*]}{2} \setminus \{ i^*j^* \}$, which can only use colours allowed by $\phi$, which form a subset of the colours allowed by $\phi^*$, by Lemma~\ref{bulk}.
So we can colour them as in $\chi$ in Step~5.
This proves that the bad colouring procedure will generate $\chi$.
Since $\chi$ was an arbitrary bad colouring, the claim is proved.
\ecpf

\medskip
\noindent
Therefore we can give an upper bound for the number of bad colourings by counting the number of steps in the bad colouring procedure.

\begin{claim2}\label{bad}
$
|\mathcal{G}| \geq (1-2^{-\gamma_3^5 n^2})|\RL^{-1}((r,\phi,\mathcal{U}))|.
$
\end{claim2}

\bcpf
By the previous claim, it suffices to bound $S_{1\ldots 5}$ from above.
Then
$$
S_{12} \leq \binom{n^2}{\leq 3\nu n^2} \cdot s^{3\nu n^2} \cdot \binom{r^*}{2} \cdot s \cdot \binom{|X_{i^*}|}{\leq\lceil \gamma_3|X_{i^*}|\rceil}\binom{|X_{j^*}|}{\leq\lceil \gamma_3|X_{j^*}|\rceil} \stackrel{(\ref{bin})}{\leq} 2^{\gamma_3^6 n^2}.
$$
Let $X \subseteq X_{i^*}$ and $Y \subseteq X_{j^*}$ be chosen at Step~2.
Now~(\ref{mu}) implies that $|X|,|Y| \geq \gamma_3 \mu n/2$.
Using $e(G[X,Y]) \leq |X||Y|$, we have 
\begin{align*}
S_{3} &\leq \sum_{i=0}^{\lfloor(|\phi^*(i^*j^*)|^{-1}-\gamma_3/2s)|X||Y|\rfloor}\binom{|X||Y|}{i} \left( |\phi^*(i^*j^*)|-1\right)^{|X||Y|-i}\\
&\leq e^{-\gamma_3^2|X||Y|/12s^2} \cdot |\phi^*(i^*j^*)|^{|X||Y|} \leq e^{-\gamma_3^4\mu^2n^2/48s^2} \cdot |\phi^*(i^*j^*)|^{|X||Y|}.
\end{align*}
where, in the second inequality, we used Corollary~\ref{binbound} with $|X||Y|,|\phi^*(i^*j^*)|,\gamma_3/2s$ playing the roles of $n,k,\delta$.
Therefore
\begin{eqnarray*}
S_{34} &\leq& S_3 \cdot |\phi^*(i^*j^*)|^{|X_{i^*}||X_{j^*}| - e(G[X,Y])} \stackrel{(\ref{almostfull})}{\leq} e^{-\gamma_3^4\mu^2 n^2/48s^2} \cdot s^{4\eps n^2} |\phi^*(i^*j^*)|^{|X_{i^*}||X_{j^*}|}.
\end{eqnarray*}

Let $B$ be the number of bad $\chi$.
Then, by Claim~\ref{badprimer},
\begin{eqnarray*}
\log B &\leq& \log S_{1\ldots 5} \leq \log \left( 2^{\gamma_3^6n^2} \cdot e^{-\gamma_3^4\mu^2n^2/48s^2} \cdot s^{4\eps n^2} \prod_{ij \in \binom{[r^*]}{2}} |\phi^*(ij)|^{|X_i||X_j|} \right)\\
&\leq& \gamma_3^6 n^2 - \frac{\log e \cdot \gamma_3^4 \mu^2 n^2}{48s^2} + \log s \cdot 4\eps n^2 + \sum_{ij \in \binom{[r^*]}{2}} |X_i||X_j|\log|\phi^*(ij)|\\
&\stackrel{(\ref{boldx})}{\leq}& -4\gamma_3^5 n^2 + \frac{q(\phi^*,\bm{x})n^2}{2} \leq -4\gamma_3^5 n^2 + \left(q(\phi^*,\ba^*) + 2\log s\|\bm{x}-\ba^*\|_1 \right) \frac{n^2}{2}\\
&\stackrel{(\ref{u})}{\leq}& Q(\bm{k})\binom{n}{2} - 3\gamma_3^5 n^2,
\end{eqnarray*}
where the penultimate inequality follows from Proposition~\ref{continuity}.
(Here, we also define $q(\phi^*,\bm{x})$ as in~(\ref{q}) even though $x_1+\ldots+x_{r^*} \leq 1$, as opposed to equal to $1$.)
Therefore
\begin{eqnarray*}
\log B &\leq& Q(\bm{k})\binom{n}{2}-3\gamma_3^5 n^2 \stackrel{(\ref{assumption})}{\leq}  \left( \log F(G;\bm{k}) + \eps\binom{n}{2} \right) - 3\gamma_3^5 n^2\\
&\leq& \log F(G;\bm{k}) - 2\gamma_3^5 n^2
\stackrel{(\ref{popularuse})}{\leq} \log |\textstyle{\RL^{-1}}((r,\phi,\mathcal{U}))| - \gamma_3^5 n^2.
\end{eqnarray*}
The claim now follows.
\ecpf

\medskip
\noindent
We would now like to adjust our partition $V(G) = X_0 \cup \ldots \cup X_{r^*}$ so that $X_0 = \emptyset$ and $|\,|X_i|-\alpha^*_i n\,| \leq 1$ for all $i \in [r^*]$, and the other properties we have proved are maintained (with slightly weaker parameters).
Clearly
\begin{equation}\label{astaru}
\|\bm{\alpha^*}-\bm{x}\|_1 = \sum_{x_i < \alpha^*_i}(\alpha^*_i-x_i) + \sum_{x_i \geq \alpha^*_i}(x_i-\alpha^*_i).
\end{equation}
For each $i \in [r^*]$, let $w_i := \min \{ \alpha^*_i,x_i \}$.
(Recall from~(\ref{boldx}) that $x_i = |X_i|/n$.)
Choose $|X_i| - \lfloor w_i n \rfloor$ vertices from each $X_i$ with $i \in [r^*]$ and choose every vertex in $X_0$.
Distribute them among the remainders of the $X_j$, $j \in [r^*]$, to create a new partition $V(G) = Y_1 \cup \ldots \cup Y_{r^*}$ such that $|\,|Y_i|-\alpha^*_i n\,| \leq 1$ for all $i \in [r^*]$.
This partition satisfies Theorem~\ref{stabilitysimp}\ref{stabilitysimpi}.

Recall Definition~\ref{popular}.
For every $(r',\phi',\mathcal{U}') \in \mathrm{Pop}(G)$, define $\mathcal{G}(r',\phi',\mathcal{U}')$ in analogy with $\mathcal{G}$ (defined with respect to $(r,\phi,\mathcal{U})$). 
Since $(r,\phi,\mathcal{U}) \in \mathrm{Pop}(G)$ chosen at the beginning of the proof was arbitrary,
\begin{align*}
\sum_{(r',\phi',\mathcal{U}') \in \mathrm{Pop}(G)} |\mathcal{G}(r',\phi',\mathcal{U}')| &\geq \left(1 - 2^{-\gamma_3^5 n^2} \right) \sum_{(r',\phi',\mathcal{U}') \in \mathrm{Pop}(G)} |\textstyle{\RL^{-1}}((r',\phi',\mathcal{U}'))|\\
&= \left(1 - 2^{-\gamma_3^5 n^2} \right) |\mathrm{Col}(G)| \geq \left(1 - 2^{-\gamma_3^5 n^2} \right)(1-2^{-2\eps n^2}) \cdot F(G;\bm{k})\\
&\geq \left(1 - 2^{-\eps n^2} \right) \cdot F(G;\bm{k}),
\end{align*}
where we used Claim~\ref{bad} and Proposition~\ref{mainlypop} in the first and third inequalities respectively.
Therefore, to prove the remainder of Theorem~\ref{stabilitysimp}, it suffices to show that every $\chi \in \mathcal{G}$ satisfies~\ref{stabilitysimpii} and~\ref{stabilitysimpiii}.

So we will now fix $\chi \in \mathcal{G}$.
Then the number of vertices which do not lie in $X_i \cap Y_i$ for any $i \in [r^*]$ is
\begin{eqnarray*}
n - \sum_{i \in [r^*]}|X_i \cap Y_i| &\leq& n - \sum_{x_i < \alpha^*_i}|X_i| - \sum_{x_i \geq \alpha^*_i}\lfloor \alpha^*_i n \rfloor\\
&\leq& n + \sum_{x_i < \alpha^*_i} (\alpha^*_i n - |X_i|) + \sum_{x_i \geq \alpha^*_i}(|X_i|-\alpha^*_i n) - \sum_{x_i < \alpha^*_i}\alpha^*_in - \sum_{x_i \geq \alpha^*_i}|X_i| + R(\bm{k})\\
&\stackrel{(\ref{astaru})}{\leq}& n + \|\ba^*-\bm{x}\|_1 n - \sum_{i \in [r^*]}\alpha^*_i n + R(\bm{k}) \stackrel{(\ref{u})}{\leq} 2\nu n + R(\bm{k}) \leq 3\nu n.
\end{eqnarray*}
Therefore
\begin{equation}\label{symdiff}
\sum_{i \in [r^*]}|X_i \bigtriangleup Y_i| = \sum_{i \in [r^*]}(|X_i|+|Y_i|-2|X_i \cap Y_i|) \leq 6\nu n.
\end{equation}
So, for all $i \in [r^*]$, our choice of $\mu$ in (\ref{stabhier}) implies that $|X_i \bigtriangleup Y_i| \leq 6\nu n \leq 7\nu|Y_i|/\mu$.
Now Proposition~\ref{adjust} implies that
$\chi^{-1}(c)[Y_i,Y_j]$ is $(\gamma_4,|\phi^*(ij)|^{-1})$-regular.
So $\chi$ satisfies Theorem~\ref{stabilitysimp}\ref{stabilitysimpii}.

We will now show that $\chi$ satisfies Theorem~\ref{stabilitysimp}\ref{stabilitysimpiii}.
Assume that
\begin{equation}
\label{gamma7}
\sum_{i \in [r^*]}e_G(Y_i) > \delta n^2> \sqrt{\gamma_6} n^2.
\end{equation}
We have that
\begin{eqnarray}
\label{ediff} \sum_{i \in [r^*]}(e_G(Y_i)-|\chi^{-1}(1)[Y_i]|) &=&  \sum_{i \in [r^*]}(e_G(Y_i)-e_G(X_i))\\
\nonumber &\phantom{+}& \hspace{-4cm} + \sum_{i \in [r^*]}(e_G(X_i)-|\chi^{-1}(1)[X_i]|) + \sum_{i \in [r^*]}(|\chi^{-1}(1)[X_i]| -|\chi^{-1}(1)[Y_i]|)\\
\nonumber &\stackrel{(\ref{eq-innergood}),(\ref{symdiff})}{\leq}& 2 \cdot 6\nu n\cdot n + 3\nu n^2 = 15\nu n^2.
\end{eqnarray}
For each $i \in [r^*]$, do the following (independently).
Let $M_{\gamma_5}$ be the integer output of Theorem~\ref{multicol} applied with $\gamma_5,1,1$ playing the roles of $\eps,s,M'$.
Recall that $|Y_i| \geq \mu n/2 > M_{\gamma_5}$.
Apply Theorem~\ref{multicol} to the monochromatic graph $\chi^{-1}(1)[Y_i]$, with parameter $\gamma_5$.
Thus obtain an equitable partition $Y_i = Z_{i,1}\cup\ldots\cup Z_{i,n_i}$ with $1/\gamma_5 \leq n_i \leq M_{\gamma_5}$ which is $\gamma_5$-regular with respect to $\chi^{-1}(1)[Y_i]$.
For each $i \in [r^*]$ and $j \in [n_i]$, we have that $|Z_{i,j}|/|Y_i| \geq M_{\gamma_5}^{-1}/2 \geq \gamma_4/\gamma_5$. 
Proposition~\ref{badrefine} now implies that, whenever $1 \in \phi^*(ii')$, we have that $\chi^{-1}(1)[Z_{i,j},Z_{i',j'}]$ is $(\gamma_5,|\phi^*(ii')|^{-1})$-regular.
Now, for each $i \in [r^*]$, we will remove any edge $xy$ from $\chi^{-1}(1)[Y_i]$ with $x \in Z_{i,j}$ and $y \in Z_{i,j'}$ such that either
$\chi^{-1}(1)[Z_{i,j},Z_{i,j'}]$ is not $(\gamma_5,\geq\! \gamma_6)$-regular; or $j=j'$.
Let $G'[Y_i]$ be the graph obtained after these removals.
Now,
\begin{eqnarray}
\nonumber \sum_{i \in [r^*]} (e_G(Y_i) - e_{G'}(Y_i)) &\stackrel{(\ref{ediff})}{\leq}& 15 \nu n^2 + \sum_{i \in [r^*]} \gamma_5\binom{n_i}{2}\left\lceil \frac{|Y_i|}{n_i} \right\rceil^2 + \sum_{i \in [r^*]}\gamma_6 \binom{n_i}{2} \left\lceil \frac{|Y_i|}{n_i} \right\rceil^2\\
\nonumber &\phantom{+}& \hspace{2cm} + \sum_{i \in [r^*]}\left\lceil \frac{|Y_i|}{n_i} \right\rceil^2\\
\label{edgecount} &\leq& 15\nu n^2 + \gamma_5 n^2/2 + \gamma_6 n^2/2 + \gamma_5 n^2 \leq \gamma_6 n^2.
\end{eqnarray}
Observe that for every $i \in [r^*]$, every edge in $G'[Y_i]$ is coloured with colour $1$ by $\chi$, and lies in a $(\gamma_5, \geq\! \gamma_6)$-regular pair.
Let $J_i$ be the graph on vertex set $[n_i]$ in which $jj'$ is an edge if and only if $\chi^{-1}(1)[Z_{i,j},Z_{i,j'}]$ is a $(\gamma_5,\geq\!\gamma_6)$-regular pair.
Let $\omega_i := \omega(J_i)$ be the size of a maximal clique in $J_i$ and let $\bm{\omega} := (\omega_1,\ldots,\omega_{r^*})$.

\begin{claim2}\label{capclaim}
$
\bm{\omega} \in \{ \bm{1} \} \cup \{ \bm{\ell} \in \mathbb{N}^{r^*} : \|\bm{\ell}\|_1 \leq k_1-1 \}.
$
\end{claim2}

\bcpf
Without loss of generality, we may suppose that, for each $i \in [r^*]$, $Z_{i,1},\ldots,Z_{i,\omega_{i}}$ span a clique in $J_i$.
Let $H$ be the graph with vertex set $\{ (i,j) : i \in [r^*], j \in [\omega_{i}] \}$ in which $\{(i,j),(i',j')\}$ is an edge if $i=i'$; or $i \neq i'$ and $1 \in \phi^*(ii')$.
Then (recalling Definition~\ref{cap})
$$
H = (\omega_{1},\ldots,\omega_{r^*})(\phi^*)^{-1}(1).
$$ 
Suppose that $H$ contains a copy of $K_{k_1}$.
Observe that, for every $\{(i,j),(i',j')\} \in E(H)$, we have that $\chi^{-1}(1)[Z_{i,j},Z_{i',j'}]$ is $(\gamma_5,\geq\! \gamma_6)$-regular.
Lemma~\ref{embed} and our choice of parameters implies that $G'$ contains a $K_{k_1}$ of colour $1$, a contradiction.
Therefore $(\omega_1,\ldots,\omega_{r^*}) \in \Capa((\phi^*)^{-1}(1),k_1)$, and Lemma~\ref{nocap} proves the claim.
\ecpf

\medskip
\noindent
For each $i \in [r^*]$, let $\ell_i$ be such that $G[Y_i]$ is $\delta$-far from being $K_{\ell_i}$-free.
Then, by~(\ref{edgecount}), $G'[Y_i]$ is $(\delta/2)$-far from being $K_{\ell_i}$-free.
So we can remove $\delta|Y_i|^2/3$ edges from $G'[Y_i]$ and there will still be a copy $T$ of $K_{\ell_i}$.
But, by the definition of $G'[Y_i]$, every edge in $T$ lies in a pair $(Z_{i,j},Z_{i,j'})$ which is $(\gamma_5,\geq\! \gamma_6)$-regular.
Thus $J_i$ contains a copy of $K_{\ell_i}$, and so $\ell_i \leq \omega_i$.
Therefore $\bm{\ell} = \bm{1}$, or $\|\bm{\ell}\|_1 \leq k_1-1$.

We claim that our assumption~(\ref{gamma7}) means that the first alternative cannot hold.
Indeed, (\ref{gamma7}) and~(\ref{edgecount}) imply that $\sum_{i \in [r^*]}e_{G'}(Y_i) \geq (\sqrt{\gamma_6}-\gamma_6) n^2$.
So there is some $i \in [r^*]$ with $e_{G'}(Y_i) \geq (\sqrt{\gamma_6}-\gamma_6)n^2/R(\bm{k}) > 0$.
Thus $J_i$ contains at least one edge, and so $\omega_i \geq 2$.
We have proved that $\|\bm{\ell}\|_1 \leq k_1-1$ as required.
This together with Lemma~\ref{strong} further implies that $\bm{k}$ does not have the strong extension property.
This completes the proof that $\chi$ satisfies Theorem~\ref{stabilitysimp}\ref{stabilitysimpiii}.
\end{proof}

We end this section with a proof of Corollary~\ref{uniform}, a stability theorem for $\bm{k}$ with the strong extension property.

\begin{proof}[Proof of Corollary~\ref{uniform}]
Let $\delta > 0$.
Let $\eps$ be the output of Theorem~\ref{stabilitysimp} applied with parameter $\delta' \leq \delta/(5s),\mu/10$, where $\mu$ is the output of Lemma~\ref{solutions}.
Now let $G$ be a graph on $n \geq n_0$ vertices such that $\log F(G;\bm{k})/\binom{n}{2} \geq Q(\bm{k})-\eps$.

Let $(r^*,\phi^*,\ba^*) \in \opt^*(\bm{k})$ be such that
at least one of the specified $(1-2^{-\eps n^2}) \cdot F(n;\bm{k})$ colourings is associated with this triple by Theorem~\ref{stabilitysimp}.
Let $Y_1, \ldots, Y_{r^*}$ be the partition of $V(G)$ given by~\ref{stabilitysimpi}.
Writing $K_{\ba^*}(n)$ for the $n$-vertex complete partite graph whose $i$th part has size $\alpha_i^* n \pm 1$, we have
\begin{equation}\label{edit}
d_{\text{edit}}(G,K_{\ba^*}(n)) \leq \sum_{ij \in \binom{[r^*]}{2}}e(\overline{G}[V_i,V_j]) + \sum_{i \in [r^*]}e(G[V_i]).
\end{equation}
Now, Part~\ref{stabilitysimpii} of Theorem~\ref{stabilitysimp} implies that, for all $ij \in \binom{[r^*]}{2}$, we have that
$$
e(G[V_i,V_j]) \geq \sum_{c \in \phi^*(ij)}|\chi^{-1}(c)[V_i,V_j]| \geq \sum_{c \in \phi^*(ij)}\left(\phi^*(ij)^{-1}-\delta'\right)|V_i||V_j| \geq \left(1-s\delta'\right)|V_i||V_j|.
$$
So
$$
\sum_{ij \in \binom{[r^*]}{2}}e(\overline{G}[V_i,V_j]) \leq \frac{\delta}{5}\cdot\sum_{ij \in \binom{[r^*]}{2}}|V_i||V_j| \leq \frac{\delta n^2}{10}.
$$
Finally, by Part~\ref{stabilitysimpiii} of Theorem~\ref{stabilitysimp},
$
\sum_{i \in [r^*]}e(G[V_i]) \leq \delta n^2/(5s)
$.
Together with (\ref{edit}), we have $d_{\text{edit}}(G,K_{\ba^*}(n)) \leq \delta n^2/5$.
Suppose $(r,\phi,\ba)$ is another triple associated with one of the specified colourings.
Then
$\|\ba-\ba^*\|_1 \cdot n^2/2 + o(n^2) \leq d_{\text{edit}}(K_{\ba}(n),K_{\ba^*}(n)) \leq 2\delta n^2/5$.
Moreover, each entry of $\ba$ and $\ba^*$ is at least $\mu$ by Lemma~\ref{solutions}, and the above inequality can only be satisfied if $r=r^*$. This completes the proof.
\end{proof}

\section{Applications}\label{applicationsec}

\subsection{Recovering some previous results}
Previous works~\cite{abks,PY} have (implicitly) solved the optimisation problem by solving a linear program with real variables $\bm{x}=(x_1,\ldots,x_t)$ such that any $(r,\phi,\ba) \in \feas^*(\bm{k})$ corresponds to some feasible $\bm{x}$ (but not necessarily vice versa).
If, for every optimal $\bm{x}$, there is some $(r,\phi,\ba)\in\feas^*(\bm{k})$ which corresponds to it, then this triple is a basic optimal solution.
Unfortunately, for all but a few small cases, the optimal solutions of the linear program do not correspond to a feasible triple.

We define a `basic' linear program, to which we will then add extra constraints.

\medskip
\noindent
\textbf{Problem $L$:}
\it
Given a sequence $\bm{k}:=(k_1,\ldots,k_s)\in \mathbb{N}^s$ of natural numbers,
determine
$
\ell^{\rm max}(\bm{k}) := \max_{\bm{d} \in D(\bm{k})}\ell(\bm{d})
$,
the maximum value of
$$
\ell(\bm{d}) := \sum_{2 \leq t \leq s}\log t \cdot d_t
$$
over the set $D(\bm{k})$ of tuples $\bm{d}=(d_2,\ldots,d_s)$
such that $0 \leq d_t \leq 1$ for all $2 \leq t \leq s$, and $\sum_{2 \leq t \leq s}td_t \leq \sum_{c \in [s]}\left(1-(k_c-1)^{-1}\right)$.
\rm

Say that $\bm{d}$ which is feasible for Problem $L$ is \emph{realisable} if there is some $(r,\phi,\ba)\in\feas^*(\bm{k})$ with
\begin{equation}\label{realisable}
d_t=2\sum_{ij \in \binom{[r]}{2}:|\phi(ij)|=t}\alpha_i \alpha_j\quad \text{for all }2 \leq t \leq s
\end{equation}
and call such a feasible triple a \emph{realisation (of $\bm{d}$)}.

\begin{lemma}\label{transfer}
Let $s \in \mathbb{N}$ and $\bm{k} \in \mathbb{N}$.
Then $Q(\bm{k}) \leq \max_{\bm{d} \in D(\bm{k})}\ell(\bm{d})$.
Moreover, the following is true.
Suppose that at least one optimal solution $\bm{d}$ to Problem $L$ is realisable. Then $\max_{\bm{d}\in D(\bm{k})}\ell(\bm{d})=Q(\bm{k})$
and $\opt^*(\bm{k})$ is the set of all $(r,\phi,\ba) \in \feas^*(\bm{k})$ which are realisations of some optimal (realisable) $\bm{d}$.
\end{lemma}

\begin{proof}
Let $(r,\phi,\ba) \in \opt^*(\bm{k})$.
For all $L \subseteq [s]$, let
$f_L := 2\sum_{ij \in \binom{[r]}{2} :\phi(ij) = L} \alpha_i\alpha_j$ and
for all $2 \leq t \leq s$ let $d_t:=\sum_{|L|=t}f_L$.
Then
$q(\phi,\bm{\alpha}) = \ell(\bm{d})$.
We have
$$
\sum_{2 \leq t \leq s}td_t = \sum_{L \subseteq [s]}|L|f_L = \sum_{c \in [s]}\sum_{L\subseteq [s]\setminus c}f_{L\cup\{c\}},
$$
so it suffices to show that $\sum_{L \subseteq [s]\setminus\{c\}}f_{L\cup\{c\}} \leq 1-(k_c-1)^{-1}$ for all $c \in [s]$.
For this, let $n \in \mathbb{N}$ and let $H_c:=H^n_c(\phi,\ba)$ be the graph on $n$ vertices with vertex classes $X_1,\ldots,X_r$ where $|\,|X_i|-\alpha_i n\,| \leq 1$ for all $i \in [r]$ and $xy \in E(H_c)$ for $x \in X_i$, $y \in X_j$ if and only if $c \in \phi(ij)$.
Then $H_c$ is $K_{k_c}$-free since $\phi^{-1}(c)$ is.
Therefore Tur\'an's Theorem~\cite{turan} implies that
$
e(H_c) \leq (1 - (k_c-1)^{-1}) n^2/2
$.
Let $c \in [s]$.
So
\begin{align*}
\frac{n^2}{2}\sum_{L \subseteq [s]\setminus \{ c \}} f_{L\cup \{ c \}} 
&= \sum_{L \subseteq [s]\setminus\{c\}} \sum_{\stackrel{ij \in \binom{[r]}{2}}{\phi(ij)=L\cup \{ c \}}}\alpha_i n \cdot \alpha_j n \leq \sum_{L \subseteq [s]\setminus \{ c \}}\sum_{\stackrel{ij \in \binom{[r]}{2}}{\phi(ij)=L\cup \{ c \}}}|X_i||X_j| + 2s^2 n\\
&= \sum_{\stackrel{ij \in \binom{[r]}{2} }{c \in \phi(ij)}} |X_i||X_j| + 2s^2n= e(H_c) + 2s^2 n \leq \left(1-\frac{1}{k_c-1}\right)\frac{n^2}{2} + 2s^2n.
\end{align*}
Dividing through by $n^2/2$ and taking the limit as $n \rightarrow \infty$ gives
the required inequality.
\end{proof}

We wish to add more constraints to Problem $L$.
Indeed, without additional constraints, Problem $L$ only yields realisable solutions in some very special cases, for example $\bm{k}=(k,k)$ or $\bm{k}=(k,k,k)$.
A constraint is \emph{valid} if every $\bm{d}$ which has a realisation $(r,\phi,\ba)\in\opt^*(\bm{k})$ must satisfy the constraint. 
We use $I$ for a set of constraints, each of the type $\sum_{t \in T}d_t \leq 1-\frac{1}{k-1}$ for some $T \subseteq \{2,\ldots,s\}$ and integer $k \geq 3$.
We call this constraint a \emph{$(T,k)$-constraint}.
Let Problem $(L,I)$ be Problem $L$ with the constraints in $I$ added to it, and 
let $\ell^{\rm max}_I(\bm{k})$ be the optimal solution of Problem $(L,I)$.
We will still discuss \emph{realisable} solutions $\bm{d}$ and \emph{realisations} of $\bm{d}$ for Problem $(L,I)$ without referring to $I$ when it is clear from the context.

For our purposes, it suffices to consider constraints as follows.
Let $T \subseteq \{2,\ldots,s\}$.
Next, given $(r,\phi,\ba)\in\feas^*(\bm{k})$, let
\begin{equation}\label{HphiT}
H_\phi(T):=\left\{ij\in\binom{[r]}{2}: |\phi(ij)|\in T\right\}.
\end{equation}
Suppose that $H_\phi(T)$ is $K_k$-free for all $(r,\phi,\ba)\in\feas^*(\bm{k})$. Then
$$
\sum_{t \in T}d_t \leq 1-\frac{1}{k-1}
$$
is a valid constraint.
This follows as in the proof of Lemma~\ref{transfer} from defining $H^n_T(\phi,\ba)$ to be the $n$-vertex $\ba$-blow-up of $H_\phi(T)$ (in analogy with $H^n_c(\phi,\ba)$) and using the observation that
$$
e(H^{n}_T(\phi,\ba)) = \left(\sum_{t \in T}d_t\right)\frac{n^2}{2}+O(n).
$$

\begin{lemma}\label{transfer2}
Let $s \in \mathbb{N}$ and $\bm{k} \in \mathbb{N}$.
Let $I$ be a set of valid $(T,k)$-constraints where each $T \subseteq \{2,\ldots,s\}$ and $k \geq 3$ is an integer.
Then $Q(\bm{k}) \leq \ell^{\rm max}_I(\bm{k})$.
Moreover, the following is true.
Suppose that at least one optimal solution $\bm{d}$ to Problem $(L,I)$ is realisable. Then $\ell^{\rm max}_I(\bm{k})=Q(\bm{k})$
and $\opt^*(\bm{k})$ is the set of all $(r,\phi,\ba) \in \feas^*(\bm{k})$ which are realisations of some optimal (realisable) $\bm{d}$.
\qed\end{lemma}

The following lemma will enable us to prove that many of the sequences $\bm{k}$ for which Problem~$Q_2$ has been solved do indeed have the strong extension property.

\begin{lemma}\label{numcheck}
Let $s \in \mathbb{N}$ and $\bm{k} \in \mathbb{N}^s$.
Suppose that, for all $(r^*,\phi^*,\bm{\alpha}^*) \in \opt^*(\bm{k})$, we have that
\begin{enumerate}[label=(\roman*),ref=(\roman*)]
\item\label{numchecki} $(\phi^*)^{-1}(c) \cong T_{k_c-1}(r^*)$ and $(k_c-1)|r^*$ for all $c \in [s]$;
\item\label{numcheckii} $|\,|\phi^*(ij)|-|\phi^*(i'j')|\,| \leq 1$ for all $ij,i'j' \in \binom{[r^*]}{2}$ and $\bm{\alpha}^*$ is uniform;
\item\label{numcheckiii} every solution $\bm{t} := (t_1,\ldots,t_{r^*}) \in [s]^{r^*}$ of
\begin{equation}\label{num}
\prod_{i \in [r^*]} t_i^{\alpha_i^*} = 2^{Q(\bm{k})}
\end{equation}
is such that $t_i = 1$ for exactly one value $i \in [r^*]$.
\end{enumerate}
Then $\bm{k}$ has the strong extension property.
\end{lemma}

\begin{proof}
Let $r^*+1$ be a new vertex and let $\phi : \binom{[r^*+1]}{2} \rightarrow 2^{[s]}$ be such that $\phi|_{\binom{[r^*]}{2}} = \phi^*$ and
\begin{equation}\label{C}
\ext(\phi,\ba^*) = Q(\bm{k}). 
\end{equation}
Since $(\phi^*)^{-1}(c) \cong T_{k_c-1}(r^*)$ for each $c \in [s]$, we have equally-sized sets $P^c_1,\ldots,P^c_{k_c-1}$ which partition $[r^*]$ and which are the vertex classes of $(\phi^*)^{-1}(c)$.
Let $\phi'$ be
obtained from $\phi$ by maximally enlarging the values on the pairs
that contain $r^*+1$ so that $(\phi')^{-1}(c)$ still does not contain
a clique on $k_c$ vertices. Clearly, for each colour $c$ this can be
done independently of the other colours and every maximal
$K_{k_c}$-free attachment of a new vertex to $(\phi^*)^{-1}(c)\cong
T_{k_c-1}(r^*)$ is to connect the vertex to all but one parts of the
Tur\'an graph. Thus for each $c\in [s]$ there is $j_c \in [k_c-1]$ such that $c \notin \phi'(\{ x,r^*+1\})$ if and only if $x \in P^c_{j_c}$.

For each $x \in [r^*]$, let $i_c(x)$ be the unique member of $[k_c-1]$ such that~$x \in P^c_{i_c(x)}$.
So $c \notin \phi'(\{ y,r^*+1 \})$ if and only if $i_c(y) = j_c$.
Then $\ext(\phi',\ba^*) \geq \ext(\phi,\ba^*)=Q(\bm{k})$ so
by Proposition~\ref{extendbd}, $\ext(\phi',\ba^*) = \ext(\phi,\ba^*) = Q(\bm{k})$, so $\phi(xy) \neq \phi'(xy)$ only if $|\phi'(xy)| = 1$. 
Observe that $\phi'$ is determined completely by $\phi^*$ and $\{ j_1,\ldots, j_s \}$.

Define $\bm{t} \in \mathbb{N}^{r^*}$ by setting $t_i := \max \{ |\phi'(\{ i,r^*+1 \})|,1\}$.
Exponentiating~(\ref{C}) implies that $\prod_{i \in [r^*]}t_i^{\alpha^*_i} = 2^{Q(\bm{k})}$. 
So, by our hypothesis~\ref{numcheckiii}, there exists $x^* \in [r^*]$ such that $|\phi'(\{ x^*,r^*+1 \})| \leq 1$; and $|\phi'(\{ i,r^*+1 \})| \geq 2$ for all $i \in [r^*] \setminus \{ x^* \}$.
Suppose first that $\phi'(\{ x^*,r^*+1\}) = \emptyset$.
Then $j_c = i_c(x^*)$ for all $c \in [s]$, and so $r^*+1$ is a twin of $x^*$, as required.

Therefore we may assume that $\phi'(\{ x^*,r^*+1\}) = \{ c^* \}$ for some $c^* \in [s]$.
Note that $j_c = i_c(x^*)$ for all $c \in [s]\setminus \{ c^* \}$.
So the attachment of $r^*+1$ is almost the same as that of $x^*$, and we will compare them to obtain a contradiction.
Without loss of generality, assume that $i_{c^*}(x^*) = 1$ and $j_{c^*} = 2$.
Now, for $i \in [k_{c^*}-1]$ and $y \in P^{c^*}_i\setminus \{x^*\}$, we have that
$$
\phi'(\{ y,r^*+1 \}) = \begin{cases} 
\phi^*(x^*y) \cup \{ c^*\} &\mbox{if } i=1. \\
\phi^*(x^*y) \setminus \{ c^*\} &\mbox{if } i=2. \\
\phi^*(x^*y) &\mbox{if }  3 \leq i \leq k_{c^*}-1. 
\end{cases}
$$
Since $\ext(\phi',\ba^*) = Q(\bm{k}) = q_{x^*}(\phi^*,\ba^*) = \sum_{x \in [r^*]}\alpha^*_x \log|\phi^*(x^*x)|$ by Lemma~\ref{lagrange}, we have that
\begin{equation}\label{sumsum}
\sum_{y \in P^{c^*}_1 \cup P^{c^*}_2 \setminus \{ x^*\}} \alpha^*_i\log|\phi^*(x^*y)| = \sum_{y \in P^{c^*}_1}\alpha^*_i\log \left(|\phi^*(x^*y)|-1 \right) + \sum_{y \in P^{c^*}_2 \setminus \{ x^*\}} \alpha^*_i\log \left(|\phi^*(x^*y)|+1\right).
\end{equation}
Let $p \in \mathbb{N}$ be such that $|\phi^*(xy)| \in \{ p,p+1 \}$ for all $xy \in \binom{[r^*]}{2}$ (which exists by~\ref{numcheckii}).
Note that $p \geq 2$.
Since $(k_{c^*}-1)|r^*$ we may write $|P^{c^*}_1| =|P^{c^*}_2| = r^*/(k_c-1) =:r$.
Suppose $\ell \leq r-1$ and $k \leq r$ are such that $|\phi^*(x^*y)| = p$ for $\ell$ elements $y$ in $P^{c^*}_1$ and $|\phi^*(x^*y)| = p$ for $k$ elements $y$ in $P^{c^*}_2$.
Then, since $\bm{\alpha}^*$ is uniform, exponentiating~(\ref{sumsum}) gives
$$
p^\ell(p+1)^{r-1-\ell}p^k(p+1)^{r-k} = (p+1)^\ell(p+2)^{r-1-\ell}(p-1)^k p^{r-k},
$$
i.e. $p^{\ell+2k-r}(p+1)^{2r-1-k-\ell} = (p+2)^{r-1-\ell}(p-1)^k$.
But $p,p-1$ are coprime, and so are $p+1,p+2$.
So $p|(p+2)$ and $(p-1)|(p+1)$.
Therefore $p = 2$, giving
$$
2^{\ell+2k-r}3^{2r-1-k-2\ell} = 2^{2r-2-2\ell}.
$$
So $\ell+2k-r=2r-2-2\ell$ and $2r-1-k-2\ell=0$, and hence $3(k+1) = 6(r-\ell) = 4(k+1)$, which implies $k=-1$, a contradiction.

Therefore $\phi'(\{ x^*,r^*+1 \}) = \emptyset$, and $\phi'(\{ x,r^*+1 \}) = \phi^*(x^*x)$ for all $x \in [r^*]\setminus \{ x^*\}$.
By our earlier observation, $\phi \equiv \phi'$.
Therefore $r^*+1$ is a twin of $x^*$, as required.
\end{proof}

In all these results, every basic optimal $(r,\phi,\ba)$ has $\phi^{-1}(c) \cong T_{k-1}(r)$ for all $c \in [s]$ and $\ba$ is the uniform vector of length $r$. The figure for $k=4$, $s=4$ is the complement of the optimal solution.

\begin{table}[h]
\centering
\caption{Basic optimal solutions}
In all these results, every basic optimal $(r,\phi,\ba)$ has $\phi^{-1}(c) \cong T_{k-1}(r)$ for all $c \in [s]$ and $\ba$ is the uniform vector of length $r$. The figure for $k=4$, $s=4$ is the complement of the optimal solution.

\vspace{0.5cm}
\label{knownsolutions2}
\begin{tabular}{rl|l|l|l}
 $\bm{k}=$ & $(k;s)$              &                     $F(\bm{k})$       & basic optimal $(r,\phi,\ba)$ & \\ 
\hline
any $k$ & $s=2$   & $1-\frac{1}{k-1}$                                 &         $r=k-1$, $|\phi|=s$       &    every $\phi^{-1}(c) \cong K_r$     \\
\hline
 & $s=3$  & $(1-\frac{1}{k-1})\log 3$                                      &     $r=k-1$, $|\phi|=s$  &    every $\phi^{-1}(c) \cong K_r$ \\
\hline
$k=3$ & $s=4$           &  $\frac{1}{4}+\frac{1}{2}\log 3$                         & $r=4$, $|\phi| \in \{ 2,3 \}$  &       \includegraphics[scale=1]{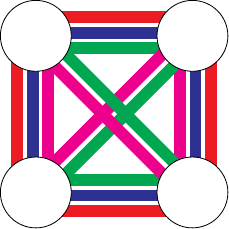}  
\\
\hline    
$k=4$ & $s=4$           &      $\frac{8}{9}\log 3$                              &   $r=9$, $|\phi| \equiv 3$ & \includegraphics[scale=0.8]{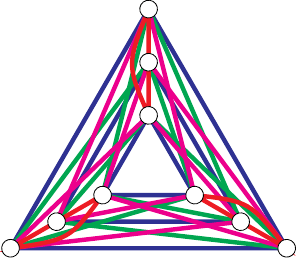}                  
\end{tabular}
\end{table}

\begin{proof}[Proof of Theorem~\ref{recover}]
First we must solve Problem $Q_2$ for all specified $\bm{k}$.
This was implicitly done in~\cite{abks,PY} but we repeat the arguments here for completeness, and to demonstrate that the arguments are much cleaner and shorter when one is working with optimal solutions rather than regularity partitions of large graphs.
We will solve Problem $Q_2$ by solving Problem $L$, sometimes with some additional valid constraints $I$, and then applying Lemma~\ref{transfer2}.
First we make some general observations.
Suppose $\bm{d}$ is a feasible solution of Problem $L$ with additional constraints $I$, each constraint corresponding to some $(T,k)$, and $\bm{d}$ has realisation $(r,\phi,\ba)$.
\begin{itemize}
\item Let $T \subseteq \{2,\ldots,s\}$ and $k \geq 3$ be such that the $(T,k)$-constraint is valid and in $I$. Suppose further that $\sum_{t \in T}d_{t}=1-\frac{1}{k-1}$ (that is, there is equality in the $(T,k)$-constraint). Then there is a partition of $[r]$ into parts $A_1,\ldots,A_{k-1}$ such that $\sum_{i \in A_{i'}}\alpha_{i}=\frac{1}{k-1}$ for all $i' \in [k-1]$, and $ij \in H_\phi(T)$ if and only if $i,j$ lie in different parts $A_{i'},A_{j'}$. 
(Recall that $H_\phi(T)$ was defined in~(\ref{HphiT}).)
\item If $S \subseteq [r]$ has $|S| \leq k$, then $2\sum_{ij \in \binom{S}{2}}\alpha_i\alpha_j \leq \sum_{i \in S}\alpha_i\left(1-\frac{1}{k-1}\right)$.
\end{itemize}

These follow as in the proof of Lemma~\ref{transfer} by taking $\ba$-weighted blow-ups of $H_\phi(T)$ and $\phi|_{\binom{S}{2}}$ respectively. For the first assertion, apply the stability theorem of Erd\H{o}s~\cite{Estab} and Simonovits~\cite{simonovits} for the Tur\'an problem, which states that any large $n$-vertex $K_k$-free graph with density close to $1-\frac{1}{k-1}$ must be close in edit distance to $T_{k-1}(n)$. For the second, apply Tur\'an's theorem.

For ease of notation we will write $H_\phi(t_1,\ldots,t_\ell)$ for $H_\phi(\{t_1,\ldots,t_\ell\})$ below.

\medskip
\noindent
\underline{\textit{The cases $\bm{k}=(k,k)$ and $\bm{k}=(k,k,k)$.}}

\medskip
\noindent
We omit $\bm{k}=(k,k)$ since it is similar to $\bm{k}=(k,k,k)$.
Problem $L$ for $\bm{k}=(k,k,k)$ 
is to maximise $d_2+\log 3\cdot d_3$ subject to $\bm{d}\geq\bm{0}$ and
$2d_2+3d_3 \leq 3(1-\frac{1}{k-1})$.
It is easy to see that the maximum is $\frac{k-2}{k-1}\log 3$ with unique optimal solution $(d_2,d_3)=(0,1-\frac{1}{k-1})$.
Now, if $(r,\phi,\ba) \in \opt^*(\bm{k})$ is a realisation of $\bm{d}$, then $H_\phi(3) \cong \phi^{-1}(c)$ for all colours $c$, so $H_\phi(3)$ is $K_k$-free.
Thus $H_\phi(3)$ is a complete $(k-1)$-partite graph and the sum of $\alpha_{i'}$ over $i'$ in a single part is $\frac{1}{k-1}$, and in fact each part is a singleton.
So $r=k-1$ and $\alpha_i=\frac{1}{k-1}$ for all $i \in [r]$, and $\phi^{-1}(c) = H_\phi(3) \cong K_{k-1}$ for all colours $c$.

\medskip
\noindent
\underline{\textit{The case $\bm{k}=(3,3,3,3)$.}}

\medskip
\noindent
We use the argument from~\cite{PY}, which requires an additional constraint.
Let $T := \{3,4\}$.
We claim that $H_\phi(T)$ is $K_3$-free for all $(r,\phi,\ba)\in\feas^*(\bm{k})$.
Indeed, if it contained a triangle $i_1i_2i_3$, then there is at most one colour in $[4]$ missing from each $\phi(i_si_t)$, and thus there is one colour in $[4]$ which appears on every edge, a contradiction.
Thus the $(\{3,4\},3)$-constraint is valid. So adding this constraint to Problem $L$, we seek to maximise $d_2 + \log 3\cdot d_3 + 2d_4$ subject to $\bm{d}\geq\bm{0}$, 
$2d_2+3d_3+4d_4 \leq 2$ and $d_3+d_4 \leq \frac{1}{2}$.
This has maximum $\frac{1}{4}+\frac{1}{2}\log 3$ with unique optimal solution $(d_2,d_3,d_4)=(\frac{1}{4},\frac{1}{2},0)$.
Thus if $(r,\phi,\ba)$ is a realisation of $\bm{d}$,
there is a partition of $[r]$ into $A,B$ such that $H_\phi(3,4)=H_\phi(3)$ is a complete bipartite graph with parts $A,B$, and $\sum_{i \in A}\alpha_{i}=\sum_{i\in B}\alpha_i=\frac{1}{2}$.
Since $d_4=0$, for distinct $i,j \in A$ or $i,j \in B$ we have $|\phi(ij)|=2$,
and, because $H_\phi(2)$ is disjoint from $H_\phi(3)$,
\begin{equation}\label{f2}
\frac{1}{4}=d_2\stackrel{(\ref{realisable})}{=}2\sum_{ij\in\binom{A}{2}}\alpha_{i}\alpha_{j}+2\sum_{ij\in\binom{B}{2}}\alpha_i\alpha_j.
\end{equation}

Without loss of generality suppose that $|A| \leq |B|$.
Next we show that $|A|=|B|=2$ via a series of claims.
Note that $|A|+|B| \geq 4$, otherwise $|A|=1$ and $|B|\leq 2$ and the second bullet point above implies that $2d_2 \leq \frac{1}{2}\cdot\frac{1}{2}$, a contradiction.
The first claim is that $|A| \leq |B| \leq 4$.
If not, then there are $a \in A$ and $b_1,\ldots,b_5 \in B$.
Since $|\phi(ab_j)|=3$ for all $j\in[5]$, we may assume that $\phi(ab_1)=\phi(ab_2)$.
But then $\phi(b_1b_2)$, of size $2$, has non-empty intersection with this set,
so $\{a,b_1,b_2\}$ span a monochromatic triangle, a contradiction which proves the claim.
The second claim is that if $|A| \geq 2$, then $|A|=|B|=2$.
If not, then there are $a_1,a_2 \in A$ and $b_1,b_2,b_3 \in B$.
Let $S$ be the multiset obtained by collecting all $\phi(a_1a_2),\phi(a_ib_j), \phi(b_jb_{j'})$ for $i \in [2]$, $j,j'\in[3]$. Then $|S|= 6\cdot 3 + 4\cdot 2 = 26$.
So there is $c \in [4]$ which appears in $\phi$ on $\lceil\frac{26}{4}\rceil=7$ pairs among $5$ vertices, so $\phi^{-1}(c)$ contains a triangle by Tur\'an's theorem.
It remains to rule out the case $|A|=1$ and $|B| \geq 3$.
Since $|B| \leq 4$, we have $2\sum_{ij \in \binom{B}{2}}\alpha_i\alpha_j \leq \frac{1}{2}\cdot\frac{3}{4}$, so $d_2 \leq \frac{3}{16}$, a contradiction.
This completes the proof that $|A|=|B|=2$.

So $r=4$, and~(\ref{f2}) holds if and only if $\alpha_i=\frac{1}{4}$ for all $i\in[4]$.
We have $\sum_{c \in [4]}|\phi^{-1}(c)| = \sum_{ij\in\binom{[r]}{2}}|\phi(ij)|=2\cdot 2+ 4\cdot 3=16$ and every $|\phi^{-1}(c)| \leq 4$ otherwise there would be a triangle in colour $c$.
Thus every $|\phi^{-1}(c)|=4$ and moreover $\phi^{-1}(c) \cong K_{2,2}=T_2(4)$.
One can check that, up to relabelling, there is a unique way to choose the $\phi^{-1}(c)$ to attain the given multiplicities (as in Table~\ref{knownsolutions2}).

\medskip
\noindent
\underline{\textit{The case $\bm{k}=(4,4,4,4)$.}}

\medskip
\noindent
No additional constraints are necessary in this case.
Problem $L$ is to maximise $d_2+\log 3\cdot d_3 + 2\cdot d_4$ subject to $\bm{d}\geq \bm{0}$ and $2d_2+3d_3+4d_4 \leq \frac{8}{3}$. This has maximum $\frac{8}{9}\cdot \log 3$, attained uniquely by $(d_2,d_3,d_4)=(0,\frac{8}{9},0)$.
Suppose $(r,\phi,\ba)\in\opt^*(\bm{k})$ is a realisation of $\bm{d}$.
Since $d_3$ is the only non-zero entry in $\bm{d}$, we have $H_\phi(3) \cong K_r$.
We claim $r \leq 9$.
If not, then
$$
\sum_{c\in[4]}|\phi^{-1}(c)[[10]]|=\sum_{ij\in\binom{[10]}{2}}|\phi(ij)|= 3\cdot\binom{10}{2}=135,
$$
so there is some $c \in [4]$ such that $\phi^{-1}(c)$ has at least $\lceil\frac{135}{4}\rceil = 34$ edges among $10$ vertices. But by Tur\'an's theorem, $\phi^{-1}(c)$ contains a $K_4$, a contradiction.
So $H_\phi(3)$ is $K_{10}$-free and $d_3=\frac{8}{9}$, so the first bullet point implies that $H_\phi(3)$ is a complete $9$-partite graph and the sum of $\alpha_i$ over all $i$ in a single part is $\frac{1}{9}$.
But $ij \in H_\phi(3)$ if and only if $\phi(ij)\neq 0$, so $r=9$ and $\alpha_i=\frac{1}{9}$ for all $i \in [9]$.
Again we must have $\phi^{-1}(c) \cong T_3(9)$ for all $c \in [4]$, and one can check that there is a unique way, up to relabelling, so choose the $\phi^{-1}(c)$ to attain the given multiplicities (see Table~\ref{knownsolutions2}, where the complement of $(r,\phi,\ba)$ is drawn, i.e.~there is an edge of colour $c$ drawn between $i$ and $j$ if and only if $c \notin \phi(ij)$).

\noindent
\underline{\textit{The strong extension property.}}
\medskip

\noindent
Given $\bm{k} \in \mathbb{N}^s$, and $(r,\phi,\bm{\alpha}) \in \opt^*(\bm{k})$ and let $\bm{t} \in [s]^{r}$ be such that
\begin{equation}\label{tprod}
\prod_{i \in [r]} t_i^{\alpha_i} = 2^{Q(\bm{k})}.
\end{equation}
Using what we have just proved about basic optimal solutions, summarised in Table~\ref{knownsolutions2}, we have the following.

\begin{table}
\centering
\vspace{0.5cm}
\label{products}
\begin{tabular}{ll|l|ll}
 $\bm{k}=$ & $(k;s)$  & $r$           &       $2^{Q(\bm{k})}$    &        \\ 
\hline
any $k$ & $s=2$  & $k-1$ & $2^{k-2}$                       &                       \\
 & $s=3$  & $k-1$ & $3^{k-2}$                               &           \\
$k=3$ & $s=4$    & $4$       &  $2 \cdot 3^2$            &      \\   
$k=4$ & $s=4$  & $9$         &     $3^8$              &                    
\end{tabular}
\end{table}
We can easily solve all of these using Lemma~\ref{numcheck}.
Indeed, in every case, $2^{Q(\bm{k})}$ is a product $p_1\ldots p_{r-1}$ of $r-1$ primes each larger than $\sqrt{s}$.
If $t_1\ldots t_{r}=2^{Q(\bm{k})}$ for positive integers $t_1,\ldots,t_{r}$,
since the $p_i$ are prime, each $t_i$ is a product of $k_i$ elements of $p_1,\ldots,p_{r-1}$ for some $k_i$.
But $p_jp_k > s$ for any $jk \in \binom{[r-1]}{2}$, so $k_i \in \{ 0,1\}$.
By the pigeonhole principle, there is exactly one $i \in [r]$ with $t_i = 1$. 
Now every $\bm{k}$ in the table satisfies the hypotheses of Lemma~\ref{numcheck}.
So each of these $\bm{k}$ have the strong extension property.
\end{proof}

\subsection{The two colour case}

We will now compute $Q(\bm{k})$ in the case when $s =2$.
When $k \geq \ell$ and $\bm{k} = (k,\ell)$ we will show that $\opt^*(\bm{k})$ depends only on $\ell$, but $\opt_0(\bm{k})$ depends on both $k,\ell$.

\begin{proof}[Proof of Lemma~\ref{2col}]
Let $(r^*,\phi^*,\ba^*) \in \opt^*(\bm{k})$.
Since $2 = s \geq |\phi^*(ij)| \geq 2$ for all $ij \in \binom{[r^*]}{2}$, we must have that $(\phi^*)^{-1}(c) \cong K_{r^*}$ for $c=1,2$.
Lemma~\ref{nocap}\ref{nocapiii} implies that $r^* \geq \ell-1$.
Therefore $r^* = \ell-1$.
So we have that
$$
q(\phi^*,\ba^*)  = 2\sum_{ij \in \binom{[\ell-1]}{2}}\alpha^*_i\alpha^*_j = 1 - \sum_{i \in [\ell-1]}(\alpha^*_i)^2 \leq 1 - \frac{1}{\ell-1},
$$
with equality if and only if $\alpha^*_i = 1/(\ell-1)$ for all $i \in [\ell-1]$.

Next we show that $\bm{k}$ has the extension property.
So suppose we can attach a vertex $\ell$ and extend $\phi^*$ to $\phi$ as in Definition~\ref{extprop}.
Then
$$
1 - \frac{1}{\ell-1} = \ext(\phi,\ba^*) = \sum_{i \in [\ell-1] : \phi(i\ell) \neq \emptyset}\frac{\log|\phi(i\ell)|}{\ell-1}
$$
so
$$
\prod_{i \in [\ell-1] : \phi(i\ell) \neq \emptyset}|\phi(i\ell)| = 2^{\ell-2}.
$$
The left-hand side is a product of at most $\ell-1$ $1$-s and $2$-s.
So there is some $j \in [\ell-1]$ such that $|\phi(i\ell)| = 2$ for all $i \in [\ell-1] \setminus \{ j \}$ and $|\phi(\ell j)| \leq 1$.
This proves that $\bm{k}$ has the extension property.
If $k=\ell$, then we must have $\phi(\ell j) = \emptyset$.
But if $k > \ell$ we can set $\phi(\ell j) = \{ 1 \}$; then $\phi^{-1}(1) \cong K_{\ell}$ and so $\phi \in \Phi_1(\ell;\bm{k})$.
So $\bm{k}$ has the strong extension property if and only if $k=\ell$. 
\end{proof}

Theorem~\ref{2colcor} follows from combining Lemma~\ref{2col} with Theorem~\ref{stabilitysimp}.

\section{Concluding remarks}\label{conclude}

In this paper we have proved a stability theorem which roughly says that all almost optimal graphs for the Erd\H{o}s-Rothschild problem are similar in structure to the blow-up of a basic optimal solution with graphs of controlled clique number added inside parts.
From this, one can systematically recover almost all known stability results.
Unfortunately Problem $Q_2$ is difficult to solve in general.
It would be very interesting to see it solved in further cases. Currently, all known solutions have been obtained by relaxing it to a linear program (which is easy to solve), whose variables are graph densities and whose constraints essentially replace combinatorial constraints such as some graph being $K_k$-free, with the linear constraint that its density must be at most $1-\frac{1}{k-1}$, by Tur\'an's theorem. 
For some few cases, solutions of this linear program correspond to feasible solutions of Problem $Q_2$, but in general they do not. 
So one possible avenue to solve it in more cases is to add more sophisticated constraints to decrease the feasible set of the linear program, which is typically much larger than that of Problem $Q_2$.

In~\cite{ps3} we apply our stability theorem to prove an exact result for every $\bm{k}$ with the strong extension property, proving a part of Conjecture~\ref{conj}.
Given Theorem~\ref{recover}, this will systematically recover most existing exact results (see Table~\ref{knownsolutions}).
For the weak extension property, it is harder to obtain an exact result as there is typically a large family of asymptotically extremal graphs, with similar structures, and these graphs could have small parts.

\section{Acknowledgments}
We are grateful to Zelealem Yilma, our collaborator on an earlier paper on this topic, for helpful discussions at the beginning of this project;
and to Emil Powierski who found and helped us correct an error in the proof of Lemma~\ref{bulk}.
We also thank an anonymous referee for their careful reading of our paper.

Oleg Pikhurko was supported by ERC Advanced Grant 101020255 and Leverhulme Research Project Grant RPG-2018-424. Katherine Staden was supported by EPSRC Fellowship EP/V025953/1.

\end{document}